\documentclass[11pt]{article}
\title{Tensor categories: A selective guided tour}
\author{Michael M\"uger \\ Institute for Mathematics, Astrophysics and Particle Physics \\ 
Radboud University Nijmegen \\ The Netherlands}

\usepackage{amssymb,amsmath, theorem,latexsym,diagrams}
\usepackage{tan-gle,color}

\usepackage{epic}
\usepackage{eepicemu}

\usepackage{epsfig}

\newlength{\dinwidth}
\newlength{\dinmargin}
\def\1#1{{\bf #1}}
\def\2#1{{\cal #1}}
\def\3#1{{\sl #1}}
\def\4#1{{\tt #1}}
\def\5#1{{\sf #1}}
\def\6#1{{\mathfrak #1}}
\def\7#1{{\mathbb #1}}

\newcommand{\be}{\begin{equation}}
\newcommand{\ee}{\end{equation}}
\newcommand{\ba}{\begin{array}}
\newcommand{\ea}{\end{array}}
\newcommand{\bea}{\begin{eqnarray}}
\newcommand{\eea}{\end{eqnarray}}
 \newcommand{\bean}{\begin{eqnarray*}}
\newcommand{\eean}{\end{eqnarray*}}

\newcommand{\ve}{\varepsilon}

\newcommand{\impl}{\Rightarrow}
\newcommand{\rarr}{\rightarrow}

\newcommand{\ol}{\overline}

\newcommand{\del}{\partial}
\newcommand{\id}{\mathrm{id}}
\newcommand{\obj}{\mathrm{Obj}}
\newcommand{\mcirc}{\,\circ\,}

\newcommand{\Hom}{\mathrm{Hom}}
\newcommand{\End}{\mathrm{End}}
\newcommand{\Aut}{\mathrm{Aut}}
\newcommand{\Rep}{\mathrm{Rep}}
\newcommand{\Vect}{\mathrm{Vect}}
\newcommand{\Tr}{\mathrm{Tr}}
\newcommand{\op}{{\mathrm{\scriptsize op}}}

\renewcommand{\mod}{\mathrm{mod}}
\newcommand{\Mod}{\mathrm{Mod}}
\newcommand{\Obj}{\mathrm{Obj}\,}

\newcommand{\DS}{\displaystyle}

\def\endexem{\hfill{$\Box$}\medskip}

\theoremstyle{change}
\theoremheaderfont{\scshape}
\newtheorem{defin}{Definition}[section]
\newtheorem{defprop}{Definition/Proposition}[section]
\newtheorem{lemma}[defin]{Lemma}
\newtheorem{prop}[defin]{Proposition}
\newtheorem{theorem}[defin]{Theorem}
\newtheorem{coro}[defin]{Corollary}
\newtheorem{conj}[defin]{Conjecture}
\theorembodyfont{\rmfamily}
\newtheorem{remark}[defin]{Remark}

\newcommand{\bdefin}{\begin{defin}}
\newcommand{\blemma}{\begin{lemma}}
\newcommand{\bprop}{\begin{prop}}
\newcommand{\btheor}{\begin{theorem}}
\newcommand{\bcoro}{\begin{coro}}
\newcommand{\bdefprop}{\begin{defprop}}
\newcommand{\edefprop}{\end{defprop}}
\newcommand{\edefin}{\end{defin}}
\newcommand{\elemma}{\end{lemma}}
\newcommand{\eprop}{\end{prop}}
\newcommand{\etheor}{\end{theorem}}
\newcommand{\ecoro}{\end{coro}}
\newcommand{\bconj}{\begin{conj}}
\newcommand{\econj}{\end{conj}}
\newcommand{\brem}{\begin{remark}}
\newcommand{\erem}{\endexem\end{remark}}

\def\mobj#1{\raise .4\unitlens\hbox{\put(0,0){$#1$}}}
\def\mychi{\raise 2pt\hbox{$\chi$}}

\newarrow{Congruent} 33333

\numberwithin{equation}{section}

\begin{document}
\maketitle

\abstract{These are the, somewhat polished and updated, lecture notes for a three hour course on
tensor  categories, given at the CIRM, Marseille, in April 2008.
The coverage in these notes is relatively  non-technical, focusing on the essential 
  ideas. They are meant to be accessible for beginners, but it is hoped that also some of the
  experts will find something interesting in them.
 
Once the basic definitions are given, the focus is mainly on categories that are linear over a field
$k$ and have finite dimensional hom-spaces. Connections with quantum groups and low dimensional topology
  are pointed out, but these notes have no pretension to cover the latter subjects to any depth.
Essentially, these notes should be considered as annotations to the extensive bibliography.
We also recommend the recent review \cite{CE}, which covers less ground in a deeper way.
}


\section{Tensor categories}\label{sec-1}
These informal notes are an outgrowth of the three hours of lectures that I gave at the Centre
International de Rencontres Mathematiques, Marseille, in April 2008. The original version of text
was projected to the screen and therefore kept maximally concise. For this publication, I have
corrected the language where needed, but no serious attempt has been made to make these notes
conform with the highest standards of exposition. I still believe that publishing them in this form
has a purpose, even if only providing some pointers to the literature.

\subsection{Strict tensor categories}
We begin with strict tensor categories, despite their limited immediate applicability.
\begin{itemize}
\item We assume that the reader has a working knowledge of categories, functors and natural
transformations. Cf.\ the standard reference \cite{cwm}. Instead of $s\in\Hom(X,Y)$ we will
occasionally write $s:X\rarr Y$.
\item We are interested in ``categories with multiplication''. (This was the title of a paper
\cite{benab1} by B\'enabou 1963, cf.\ also Mac Lane \cite{macl2} from the same year). 
This term was soon replaced by `monoidal categories' or `tensor categories'. (We use these
synonymously.) It is mysterious to this author why the explicit formalization of tensor categories
took twenty years to arrive after that of categories, in particular since monoidal categories
appear in protean form, e.g., in Tannaka's work \cite{tannaka}.  

\item A {\bf strict tensor category} (strict monoidal category)
is a triple $(\2C,\otimes,\11)$, where $\2C$ is a category, $\11$ a
  distinguished object and $\otimes:\2C\times\2C\rarr\2C$ is a functor, satisfying
\[ (X\otimes Y)\otimes Z=X\otimes(Y\otimes Z) \quad\quad\mbox{and}\quad\quad
    X\otimes\11=X=\11\otimes X \ \ \forall X,Y,Z.\]

If $(\2C,\otimes,\11),(\2C',\otimes',\11')$ are strict tensor categories, a {\bf strict tensor
functor} $\2C\rarr\2C'$ is a functor $F:\2C\rarr\2C'$ such that 
\[ F(X\otimes Y)=F(X)\otimes'F(Y), \quad\quad F(\11)=\11'. \]

If $F,F':\2C\rarr\2C'$ are strict tensor functors, a natural transformation $\alpha:F\rarr F'$ is
{\bf monoidal} if and only if $\alpha_\11=\id_{\11'}$ and
\[ \alpha_{X\otimes Y}=\alpha_X\otimes\alpha_Y\quad \forall X,Y\in\2C. \]
(Both sides live in $\Hom(F(X\otimes Y),F'(X\otimes Y))=\Hom(F(X)\otimes' F(Y),F'(X)\otimes' F'(Y))$.)
\item WARNING: The coherence theorems, to be discussed in a bit more detail in Subsection
  \ref{ss-nonstrict}, will imply that, in a sense, strict tensor categories are sufficient for all
  purposes. However, even when dealing with strict tensor categories, one needs non-strict tensor functors! 
\item Basic examples:
\begin{itemize}
\item Let $\2C$ be any category and let $\End\,\2C$ be the category of functors $\2C\rarr\2C$ and
their natural transformations. Then $\End\,\2C$ is a strict $\otimes$-category, with composition of
functors as tensor product. It is also denoted as the `center' $Z_0(\2C)$. (The subscript is needed
since various other centers will be encountered.)
\item To every group $G$, we associate the discrete tensor category $\2C(G)$:
\[ \obj\,\2C(G)=G, \ \ \Hom(g,h)=\left\{ \begin{array}{cc} \{ \id_g \} &  g=h \\ \emptyset & g\ne
   h\end{array}\right., 
 \quad g\otimes h=gh.\]

\item The {\bf symmetric category} $\7S$: 
\[ \obj\,\2S=\7Z_+, \ \ \Hom(n,m)=\left\{ \begin{array}{cc} S_n &  n=m \\ \emptyset & n\ne m\end{array}\right.,
 \quad n\otimes m=n+m\]
with tensor product of morphisms given by the obvious map $S_n\times S_m\rarr S_{n+m}$. 

Remark: 1. $\7S$ is the free symmetric tensor category on one monoidal generator.

2. $\7S$ is equivalent to the category of finite sets and bijective maps.

2. This construction works with any family $(G_i)$ of groups with an associative composition
$G_i\times G_j\rarr G_{i+j}$.

\item Let $\2A$ be a unital associative algebra with unit over some field. We define $\End\,\2A$ to
have as objects the unital algebra homomorphisms $\rho:\2A\rarr\2A$. The morphisms are defined by
\[ \Hom(\rho,\sigma)=\{ x\in\2A\ | \ x\rho(y)=\sigma(y)x\ \forall y\in\2A\} \]
with $s\circ t=st$ and $s\otimes t=s\rho(t)=\rho'(t)s$ for $s\in\Hom(\rho,\rho'), t\in\Hom(\sigma,\sigma')$.
This construction has important applications in in subfactor theory \cite{lo1} and (algebraic)
quantum field theory \cite{dhr3,frs1}. Yamagami \cite{yama2} proved that every countably generated
$C^*$-tensor category with conjugates (cf.\ below) embeds fully into $\End\,\2A$ for some
von Neumann-algebra $\2A=\2A(\2C)$. (See the final section for a conjecture concerning an algebra
that should work for all such categories.) 

\item The {\bf Temperley-Lieb categories} $\2T\2L(\tau)$. (Cf.\ e.g.\ \cite{gw}.) Let $k$ be a
field and $\tau\in k^*$. We define
\[ \obj\,\2T\2L(\tau)=\7Z_+,\quad n\otimes m=n+m,\]
as for the free symmetric category $\2S$. But now
\[ \Hom(n,m)=\mathrm{span}_k\{ \mbox{Isotopy classes of } (n,m)\mbox{-TL diagrams}\}.\]
Here, an $(n,m)$-diagram is a planar diagram where $n$ points on a line and $m$ points
on a parallel line are connected by lines without crossings. The following example of a (7,5)-TL
diagram will explain this sufficiently:
\begin{center}
\setlength{\unitlength}{0.00087489in}
\begingroup\makeatletter\ifx\SetFigFont\undefined%
\gdef\SetFigFont#1#2#3#4#5{%
  \reset@font\fontsize{#1}{#2pt}%
  \fontfamily{#3}\fontseries{#4}\fontshape{#5}%
  \selectfont}%
\fi\endgroup%
\renewcommand{\dashlinestretch}{30}
\begin{picture}(1824,2084)(0,-10)
\path(12,2047)(1812,2047)(1812,22)
    (12,22)(12,2047)
\thicklines
\path(462,2047)(463,2046)(464,2044)
    (466,2040)(470,2034)(475,2025)
    (482,2013)(490,1998)(500,1981)
    (512,1960)(526,1936)(541,1909)
    (557,1880)(575,1849)(594,1816)
    (613,1781)(634,1745)(654,1707)
    (676,1668)(698,1628)(720,1587)
    (742,1545)(765,1502)(788,1457)
    (812,1411)(836,1364)(861,1315)
    (886,1264)(912,1211)(939,1156)
    (966,1099)(993,1041)(1020,982)
    (1047,922)(1076,855)(1104,791)
    (1129,730)(1153,672)(1174,618)
    (1194,567)(1212,520)(1228,476)
    (1242,434)(1256,395)(1268,358)
    (1280,322)(1290,289)(1300,256)
    (1309,226)(1317,197)(1324,170)
    (1331,145)(1337,121)(1343,101)
    (1347,82)(1351,66)(1355,53)
    (1357,42)(1359,34)(1360,28)
    (1361,25)(1362,23)(1362,22)
\path(237,2047)(237,22)
\path(462,22)(463,25)(465,32)
    (470,45)(476,63)(484,86)
    (494,114)(505,144)(517,175)
    (529,205)(541,234)(552,261)
    (564,286)(575,308)(586,328)
    (596,345)(607,360)(618,374)
    (630,386)(642,397)(655,407)
    (668,416)(682,424)(697,431)
    (713,437)(729,442)(746,446)
    (764,449)(781,451)(800,451)
    (818,451)(835,449)(853,446)
    (870,442)(886,437)(902,431)
    (917,424)(931,416)(944,407)
    (957,397)(969,386)(981,374)
    (992,360)(1003,345)(1013,328)
    (1024,308)(1035,286)(1047,261)
    (1058,234)(1070,205)(1082,175)
    (1094,144)(1105,114)(1115,86)
    (1123,63)(1129,45)(1134,32)
    (1136,25)(1137,22)
\path(1137,2047)(1138,2045)(1139,2041)
    (1141,2033)(1144,2021)(1149,2004)
    (1155,1981)(1163,1952)(1173,1917)
    (1184,1877)(1196,1832)(1210,1783)
    (1224,1729)(1240,1673)(1255,1614)
    (1272,1555)(1288,1495)(1304,1435)
    (1320,1376)(1336,1318)(1351,1262)
    (1365,1208)(1379,1157)(1392,1107)
    (1405,1060)(1416,1016)(1427,974)
    (1438,934)(1448,896)(1457,860)
    (1465,826)(1473,794)(1481,764)
    (1488,734)(1494,707)(1501,680)
    (1506,654)(1512,629)(1520,592)
    (1528,556)(1535,521)(1541,488)
    (1546,456)(1552,425)(1556,393)
    (1561,362)(1564,330)(1568,298)
    (1571,266)(1574,233)(1577,201)
    (1579,170)(1581,139)(1583,111)
    (1584,87)(1585,66)(1586,49)
    (1586,36)(1587,28)(1587,24)(1587,22)
\path(687,2047)(687,2046)(689,2041)
    (693,2029)(699,2009)(708,1984)
    (717,1956)(728,1928)(737,1903)
    (746,1881)(755,1865)(763,1852)
    (770,1843)(777,1839)(785,1837)
    (792,1839)(800,1843)(809,1852)
    (818,1865)(829,1881)(841,1903)
    (854,1928)(868,1956)(882,1984)
    (894,2009)(904,2029)(909,2041)
    (912,2046)(912,2047)
\path(687,22)(687,23)(689,28)
    (693,40)(699,60)(708,85)
    (717,113)(728,141)(737,166)
    (746,188)(755,204)(763,217)
    (770,226)(777,230)(785,232)
    (792,230)(800,226)(809,217)
    (818,204)(829,188)(841,166)
    (854,141)(868,113)(882,85)
    (894,60)(904,40)(909,28)
    (912,23)(912,22)
\end{picture}
\end{center}
The tensor product of morphisms is given by horizontal juxtaposition, whereas composition of
morphisms is defined by vertical juxtaposition, followed by removal all newly formed closed circles
and multiplication by a factor $\tau$ for each circle. (This makes sense  since the category is
$k$-linear.) 

Remark: 1. The Temperley-Lieb algebras TL$(n,\tau)=\End_{\2T\2L(\tau)}(n)$ first appeared in the theory of
exactly soluble lattice models of statistical mechanics. They, as well as $\2T\2L(\tau)$ are closely
related to the Jones polynomial \cite{vfr1} and the quantum group $SL_q(2)$. Cf.\ \cite[Chapter
XII]{turaev}. 

2. The Temperley-Lieb algebras, as well as the categories $\2T\2L(\tau)$ can be defined purely
algebraically in terms of generators and relations.

\item In dealing with (strict) tensor categories, it is often convenient to adopt a graphical
notation for morphisms:
\[ s: X\rarr Y \quad\Leftrightarrow \quad\quad
    \begin{tangle}\object{Y}\\ \hh\id\\ \O{s}\\ \hh\id\\ \object{X}\end{tangle} \]

If $s:X\rarr Y,\ t:Y\rarr Z,\ u: Z\rarr W$ then we write
\[ t\circ s: X\rarr Z \quad\Leftrightarrow\quad\quad
  \begin{tangle}\object{Z}\\ \hh\id\\ \O{t}\\  \O{s}\\ \hh\id\\ \object{X}\end{tangle}
\quad\quad\quad\quad
  s\otimes u:X\otimes Z\rarr Y\otimes W \quad\Leftrightarrow\quad\quad\quad
  \begin{tangle}\object{Y}\Step\object{W}\\ \hh\id\Step\id\\ \O{s}\Step\O{u}\\ \hh\id\Step\id\\
  \object{X}\Step\object{Z}\end{tangle}\]
\end{itemize}
The usefulness of this notation becomes apparent when there are morphisms with `different numbers of
in- and outputs': Let, e.g., 
$a:X\rarr S\otimes T,\ b:\11\rarr U\otimes Z,\ c: S\rarr\11,\ d: T\otimes U\rarr V,\ e:Z\otimes Y\rarr W$
and consider the composite morphism
\be\label{e} c\otimes d\otimes e\mcirc a\otimes b\otimes\id_Y: X\otimes Y\rarr V\otimes W. \ee
This formula is almost unintelligible. (In order to economize on brackets, we follow the majority of
authors and declare $\otimes$ to bind stronger than $\circ$, i.e.\ 
$a\circ b\otimes c\equiv a\circ(b\otimes c)$. Notice that inserting brackets in (\ref{e}) does
nothing to render the formula noticeably more intelligible.)  
It is not even clear whether it represents a morphism in the
category. This is immediately obvious from the diagram:
\[ \begin{tangle}
\step[3]\object{V}\step[2]\object{W}\\
\hh\hstep\step[2.5]\id\step[2]\id\\
\hh\frabox{c}\step[1.5]\frabox{d}\step[1.5]\frabox{e}\\
\id\obj{S}\step[2]\step[-.5]\obj{T}\hstep\ne2\hstep\id\obj{U}\step[1.5]\obj{Z}\step[-.5]\dd\hstep\id\\
\hh\frabox{a}\step[2.5]\frabox{b}\step[1.5]\id\\
\hh\step\id\step[5]\id\\
\step\object{X}\step[5]\object{Y}
\end{tangle}\]
Often, there is more than one way to translate a diagram into a formula, e.g.\ 
\[ \begin{tangle} \object{Z}\step[2]\object{Z'}\\
  \O{t}\Step\O{t'} \\  \O{t}\Step\O{t'} \\ \object{X}\step[2]\object{X'}\end{tangle}\]
can be read as $t\otimes t'\mcirc s\otimes s'$ or as $(t\circ s)\otimes(t'\otimes s')$. But by the
interchange law (which is just the functoriality of $\otimes$), these two morphisms
coincide. For proofs of consistency of the formalism, cf.\ \cite{js6,FY2} or \cite{kassel}. 

\subsection{Non-strict tensor categories}\label{ss-nonstrict}
\item For almost all situations where tensor categories arise, strict tensor categories are not
general enough, the main reasons being:

\begin{itemize}
\item Requiring equality of objects as in $(X\otimes Y)\otimes Z=X\otimes(Y\otimes Z)$
is highly unnatural from a categorical point of view.
\item Many would-be tensor categories are not strict; in particular this is the case for Vect$_k$,
as well as for representation categories of groups (irrespective of the class of groups and
representations under consideration).
\end{itemize}

\item The obvious minimal modification, namely to require only existence of isomorphisms 
$(X\otimes Y)\otimes Z\cong  X\otimes(Y\otimes Z)$ for all $X,Y,Z$ and 
$\11\otimes X\cong X\cong X\otimes\11$ for all $X$, turns out to be too weak to be useful.
\item The correct definition of not-necessarily-strict {\bf tensor categories} was given in
\cite{benab1}: It is a sextuplet $(\2C,\otimes,\11,\alpha,\lambda,\rho)$, where $\2C$ is a 
category, $\otimes:\2C\times\2C\rarr\2C$ a functor, $\11$ an object, and
$\alpha:\otimes\circ(\otimes\times\id)\rarr\otimes\circ(\id\times\otimes)$,
$\lambda:\11\otimes -\rarr\id$, $\rho: -\otimes\11\rarr\id$ are natural isomorphisms 
 (i.e., for all $X,Y,Z$ we have isomorphisms 
$\alpha_{X,Y,Z}:(X\otimes Y)\otimes Z\rarr X\otimes(Y\otimes Z)$ and $\lambda_X:\11\otimes X\rarr X$,
$\rho_X:X\otimes\11\rarr X$) such that all morphisms between the same pair of objects that can be
built from $\alpha,\lambda,\rho$ coincide. (Examples of what this means are given by the
commutativity of the following two diagrams.)

\item There are two versions of the coherence theorem for tensor categories:\\
Version I (Mac Lane \cite{macl2,cwm}): All morphisms built from $\alpha,\lambda,\rho$
are unique provided $\alpha$ satisfies the pentagon identity, i.e.\ commutativity of
\begin{diagram}
((X\otimes Y)\otimes Z)\otimes T & \rTo^{\alpha_{X,Y,Z}\otimes\id_T} & (X\otimes(Y\otimes Z))\otimes T & 
  \rTo^{\alpha_{X,Y\otimes Z,T}} &  X\otimes((Y\otimes Z)\otimes T) \\
\dTo^{\alpha_{X\otimes Y,Z,T}} &&&& \dTo_{\id_X\otimes\alpha_{Y,Z,T}} \\
(X\otimes Y)\otimes(Z\otimes T) &&\rTo_{\alpha_{X,Y,Z\otimes T}}&& X\otimes (Y\otimes(Z\otimes T))
\end{diagram}
and $\lambda,\rho$ satisfy the unit identity
\begin{diagram} (X\otimes\11)\otimes Y & \rTo^{\alpha_{X,\11,Y}} &X\otimes(\11\otimes Y) \\
   \dTo^{\rho_X\otimes\id_Y} && \dTo_{\id_X\otimes\lambda_Y} \\
   X\otimes Y & \rCongruent &  X\otimes Y
\end{diagram}
For modern expositions of the coherence theorem see \cite{cwm,kassel}. (Notice that the original
definition of non-strict tensor categories given in \cite{macl2} was modified in slightly
\cite{kelly1,kelly2}.) 

\item Examples of non-strict tensor categories:

\begin{itemize}
\item Let $\2C$ be a category with products and terminal object $T$. Define $X\otimes Y=X\prod Y$
(for each pair $X,Y$ choose a product, non-uniquely) and $\11=T$. Then $(\2C,\otimes,\11)$ is
non-strict tensor category. (Existence of associator and unit isomorphisms follows from the
universal properties of product and terminal object). An analogous construction works with coproduct
and initial object. 
\item Vect$_k$ with $\alpha_{U,V,W}$ defined on simple tensors by 
$(u\otimes v)\otimes w\mapsto u\otimes (v\otimes w)$. Note: This trivially satisfies the pentagon
identity, but the other choice $(u\otimes v)\otimes w\mapsto - u\otimes (v\otimes w)$ does not!
\item Let $G$ be a group, $A$ an abelian group (written multiplicatively) and $\omega\in Z^3(G,A)$, i.e.\ 
\[ \omega(h,k,l)\omega(g,hk,l)\omega(g,h,k)=\omega(gh,k,l)\omega(g,h,kl) \quad\forall g,h,k,l\in G.\]
Define $\2C(G,\omega)$ by 
\[ \obj\,\2C=G, \ \ \Hom(g,h)=\left\{ \begin{array}{cc} A &  g=h \\ \emptyset & g\ne h\end{array}\right.,
 \quad g\otimes h=gh.\]
with associator $\alpha=\omega$, cf.\ \cite{sinh}. If $k$ is a field, $A=k^*$, one
has a $k$-linear version where 
$\Hom(g,h)=\left\{ \begin{array}{cc} k &  g=h \\ \{0\} & g\ne h\end{array}\right.$.
I denote this by $\2C_k(G,\omega)$, but also $\mathrm{Vect}^G_\omega$ appears in the literature.

The importance of this example lies in its showing relations between categories and cohomology,
which are reinforced by `higher category theory', cf.\ e.g.\ \cite{baezS}. 
But also the concrete example is relevant for the classification of fusion categories, at least 
the large class of `group theoretical  categories'. (Cf.\ Ostrik et al.\ \cite{ostrik2,eno}.) See
Section \ref{sec-3}.

\item A  categorical group is a tensor category that is a groupoid (all morphisms are invertible)
  and where every object has a tensor-inverse, i.e.\ for every $X$ there is an object $\ol{X}$ such
 that $X\otimes\ol{X}\cong\11$. The categories $\2C(G,\omega)$ are just the skeletal categorical groups.
\end{itemize}

\item Now we can give the general definition of a tensor functor (between non-strict tensor
categories or non-strict tensor functors between strict tensor categories): A tensor functor between
tensor categories $(\2C,\otimes,\11,\alpha,\lambda,\rho),(\2C',\otimes',\11',\alpha',\lambda',\rho')$
consists of a functor $F:\2C\rarr\2C'$, an isomorphism $e^F:F(\11)\rarr\11'$ and a family of natural
isomorphisms $d^F_{X,Y}: F(X)\otimes F(Y)\rarr F(X\otimes Y)$ satisfying commutativity of
\begin{diagram}
(F(X)\otimes' F(Y))\otimes' F(Z) & \rTo^{d_{X,Y}\otimes\id} & F(X\otimes Y)\otimes' F(Z) 
   & \rTo^{d_{X\otimes Y,Z}} & F((X\otimes Y)\otimes Z) \\
\dTo^{\alpha'_{F(X),F(Y),F(Z)}} &&&& \dTo_{F(\alpha_{X,Y,Z})} \\
 F(X)\otimes'(F(Y)\otimes' F(Z)) & \rTo_{\id\otimes d_{Y,Z}} & F(X)\otimes'F(Y\otimes Z) 
   &\rTo_{d_{X,Y\otimes Z}} & F(X\otimes(Y\otimes Z))
\end{diagram}
(notice that this is a 2-cocycle condition, in particular when $\alpha\equiv\id$) and
\[\begin{array}{cc} 
\begin{diagram} F(X)\otimes F(\11) & \rTo^{\id\otimes e^F} & F(X)\otimes \11' \\
\dTo_{d^F_{X,\11}} && \dTo^{\rho'_{F(X)}}\\ F(X\otimes\11) & \rTo_{F(\rho_X)} & F(X)\end{diagram}
& \quad\quad\quad\quad\mbox{(and similar for $\lambda_X$)} \end{array}\]

Remark: Occasionally, functors as defined above are called {\bf strong} tensor functors in order to
distinguish them from the lax variant, where the $d^F_{X,Y}$ and $e^F$ are not required to be
isomorphisms. (In this case it also makes sense to consider $d^F, e^F$ with source and target
exchanged.)

\item Let $(\2C,\otimes,\11,\alpha,\lambda,\rho),(\2C',\otimes',\11',\alpha',\lambda',\rho')$ be
  tensor categories and $(F,d,e),(F',d',e'):\2C\rarr\2C'$ tensor functors. Then a natural transformation
  $\alpha:F\rarr F'$ is monoidal if
\begin{diagram} F(X)\otimes' F(Y) & \rTo^{d_{X,Y}} & F(X\otimes Y) \\
\dTo^{\alpha_X\otimes\alpha_Y} && \dTo_{\alpha_{X\otimes Y}} \\
F'(X)\otimes' F'(Y) & \rTo^{d'_{X,Y}} & F'(X\otimes Y) \end{diagram}
For strict tensor functors, we have $d\equiv\id\equiv d'$, and we obtain the earlier condition.

\item A tensor functor
$F: (\2C,\otimes,\11,\alpha,\lambda,\rho)\rarr(\2C',\otimes',\11',\alpha',\lambda',\rho')$ is called an
equivalence if there exist a tensor functor $G:\2C'\rarr\2C$ and natural monoidal
isomorphisms $\alpha:G\circ F\rarr\id_\2C$ and $\beta:F\circ G\rarr\id_{\2C'}$. For the existence of
such a $G$ it is necessary and sufficient that $F$ be full, faithful and essentially surjective (and
of course monoidal), cf.\ \cite{sr}. (We follow the practice of not worrying too much about size
issues and assuming a sufficiently strong version of the axiom of choice for classes. On this
matter, cf.\ the different discussions of foundational issues given in the two editions of \cite{cwm}.)

\item Given a group $G$ and $\omega,\omega'\in Z^3(G,A)$, the identity functor is part of a
  monoidal equivalence $\2C(G,\omega)\rarr\2C(G,\omega')$ if and only if $[\omega]=[\omega']$ in
  $H^3(G,A)$. Cf.\ e.g. \cite[Chapter 2]{dav1}. Since categorical groups form a 2-category $\2C\2G$,
  they are best classified by providing a 2-equivalence between $\2C\2G$ and a 2-category $\2H^3$
  defined in terms of cohomology groups $H^3(G,A)$. The details are too involved to give here;
  cf.\ \cite{js0}. 
(Unfortunately, the theory of categorical groups is marred by the fact that important works
  \cite{sinh,js0} were never formally published. For a comprehensive recent treatment cf.\ \cite{baezL}.)

\item Version II of the Coherence theorem (equivalent to Version I): Every tensor category is
  monoidally equivalent to a strict one. \cite{cwm,kassel}. As mentioned earlier, this allows us to
  pretend that all tensor categories are strict. (But we cannot restrict ourselves to strict tensor functors!) 

\item One may ask what the strictification of $\2C(G,\omega)$ looks like. The answer is somewhat
  complicated, cf.\ \cite{js0}: It involves the free group on the set underlying $G$. (This shows
  that sometimes it is actually more convenient to work with non-strict categories!)

\item As shown in \cite{schau3}, many non-strict tensor categories can be turned into equivalent
strict ones by changing only the tensor functor $\otimes$, but {\it leaving the underlying
category unchanged}.

\item We recall the ``Eckmann-Hilton argument'': If a set has two monoid structures $\star_1,\star_2$  
  satisfying $(a\star_2b)\star_1(c\star_2d)=(a\star_1c)\star_2(b\star_1d)$ with the same unit, the
two products coincide and are commutative. If $\2C$ is a tensor category and we consider $\End\,\11$
with $\star_1=\circ, \star_2=\otimes$ we find that $\End\,\11$ is commutative, cf.\ \cite{kel3}.
In the Ab- ($k$-linear) case, defined in Subsection \ref{ss-ab}, $\End\,\11$ is a commutative unital ring
($k$-algebra). (Another classical application of the Eckmann-Hilton argument is the abelianness of
the higher homotopy groups $\pi_n(X), n\ge 2$ and of $\pi_1(M)$ for a topological monoid $M$.)

\subsection{Generalization: 2-categories and bicategories}
\item Tensor categories have a very natural and useful generalization. We begin with `2-categories',
  which generalize strict tensor categories: A 2-category $\2E$ consists of a set (class) of objects
and, for every $X,Y\in\Obj\,\2E$, a category HOM$(X,Y)$. The objects (morphisms) in HOM$(X,Y)$ are
called 1-morphisms (2-morphisms) of $\2E$. For the detailed axioms we refer to the references given
below. In particular, we have functors 
$\circ: \mathrm{HOM}(\6A,\6B)\times\mathrm{HOM}(\6B,\6C)\rarr\mathrm{HOM}(\6A,\6C)$, and $\circ$ is
associative (on the nose). 
\item The prototypical example of a 2-category is the 2-category $\2C\2A\2T$. Its objects are the
  small categories, its 1-morphisms are functors and the 2-morphisms are natural transformations.
\item We notice that if $\2E$ is a 2-category and $X\in\Obj\,\2E$, then
END$(X)=\mathrm{HOM}(X,X)$ is a strict tensor category. 
This leads to the non-strict version of 2-categories called bicategories: We replace the
associativity of the composition $\circ$ of 1-morphisms by the existence of invertible 2-morphisms 
$(X\circ Y)\circ Z\rarr X\circ(Y\circ Z)$ satisfying axioms generalizing those of a tensor
category. Now, if $\2E$ is a bicategory and 
$X\in\Obj\,\2E$, then END$(X)=\mathrm{HOM}(X,X)$ is a (non-strict) tensor category. 
Bicategories are a very important generalization of tensor categories, and we'll meet them
again. Also the relation between bicategories and tensor categories is prototypical for `higher
category theory'. 

References: \cite{ks} for 2-categories and \cite{benab3} for bicategories, as well as the very
recent review by Lack \cite{lack}.

\subsection{Categorification of monoids}
Tensor categories (or monoidal categories) can be considered as the categorification of the notion
of a monoid. This has interesting consequences:

\item Monoids in monoidal categories: 
Let $(\2C,\otimes,\11)$ be a strict $\otimes$-category. A monoid in $\2C$ (B\'enabou \cite{benab2})
is a triple $(A,m,\eta)$ with $A\in\2C,\ m:A\otimes A\rarr A,\ \eta:\11\rarr A$ satisfying
\[ m\mcirc m\otimes\id_A=m\mcirc \id_A\otimes m,\quad\quad
m\mcirc\eta\otimes\id_A=\id_A=m\mcirc\id_A\otimes\eta. \]
(In the non-strict case, insert an associator at the obvious place.) A monoid in Ab ($\Vect_k$)
is a ring ($k$-algebra). Therefore, in the recent literature monoids are often called `algebras'.

Monoids in monoidal categories are a prototypical example of the `microcosm principle' of Baez and 
Dolan \cite{baez4} affirming that ``certain algebraic structures can be   defined in any category
equipped  with a categorified version of the same structure''. 

\item If $\2C$ is any category, monoids in the tensor category $\End\,\2C$ are known as
`monads'. As such they are older than tensor categories! Cf.\ \cite{cwm}.

\item If $(A,m,\eta)$ is a monoid in the strict tensor category $\2C$, a left A-module is a pair
$(X,\mu)$, where $X\in\2C$ and $\mu:A\otimes X\rarr X$ satisfies
\[ \mu\mcirc m\otimes\id_X=\mu\mcirc \id_A\otimes \mu, \quad\quad \mu\mcirc\eta\otimes\id_X=\id_X. \]
Together with the obvious notion of A-module morphism
\[ \Hom_{A-\Mod}((X,\mu),(X',\mu'))=\{ s\in\Hom_\2C(X,X')\ | \ s\mcirc \mu=\mu'\mcirc\id_A\otimes s\},\]
$A$-modules form a category. Right A-modules and $A-A$ bimodules are defined analogously.

The free A-module of rank 1 is just $(A,m)$.
\item If $\2C$ is abelian, then $A-\Mod_\2C$ is abelian under weak assumptions, cf.\
\cite{ardi}. (The latter are satisfied when $A$ has duals, as e.g.\ when it is a strongly separable
Frobenius algebra \cite{FS}. All this could also be deduced from \cite{EM}.)

\item Every monoid $(A,m,\eta)$ in $\2C$ gives rise to a monoid $\Gamma_A=\Hom(\11,A)$ in the category
$\2S\2E\2T$ of sets. We call it the {\bf elements of $A$}. ($\Gamma_A$ is related to the
  endomorphisms of the unit object in the tensor categories of $A-A$-bimodules and $A$-modules (in
the braided case), when the latter exist.) 

\item Let $\2C$ be abelian and $(A,m,\eta)$ an algebra in $\2C$. An {\bf ideal} in $A$ is an A-module
  $(X,\mu)$ together with a monic morphism $(X,\mu)\hookrightarrow(A,m)$. Much as in ordinary algebra, one
  can define a quotient algebra $A/I$. Furthermore, every ideal is contained in a maximal ideal, and
an ideal $I\subset A$ in a commutative monoid is maximal if and only if the ring $\Gamma_{A/I}$ is a
field. (For the last claim, cf.\ \cite{mue16}.)

\item {\bf Coalgebras} and their comodules are defined analogously. In a tensor category equipped
  with a symmetry or braiding $c$ (cf.\ below), it makes sense to say that an (co)algebra is
  (co){\bf commutative}. For an algebra $(A,m,\eta)$ this means that $m\circ c_{A,A}=m$.

\item (B) Just as monoids can act on sets, tensor categories can act on categories:

Let $\2C$ be a tensor category. A left {\bf $\2C$-module category} is a pair ($\2M,F)$ where $\2M$
  is a category and $F:\2C\rarr\End\,\2M$ is a tensor functor. (Here, $\End\,\2M$ is as in our first
  example of a tensor category.) This is equivalent to having a functor   $F':\2C\times\2M\rarr\2M$ and
  natural isomorphisms $\beta_{X,Y,A}: F'(X\otimes Y,A)\rarr F(X,F(Y,A))$ satisfying a pentagon-type
  coherence law, unit constraints, etc. Now one can define indecomposable module categories,
  etc. (Ostrik \cite{ostrik}) 

\item There is a close connection between module categories and categories of modules:

If $(A,m,\eta)$ is an   algebra in $\2C$, then there is an natural right $\2C$-module
structure on the category $A-\Mod_\2C$ of left A-modules:
\[ F': A-\Mod_\2C\times\2C, \ \ (M,\mu)\times X\mapsto (M\otimes X,\mu\otimes\id_X). \]
(In the case where $(M,\mu)$ is the free rank-one module $(A,m)$, this gives the free $A$-modules
$F'((A,m),X)=(A\otimes X,m\otimes\id_X)$.)
For a fusion category (cf.\ below), one can show that every semisimple indecomposable left
$\2C$-module category arises in this way from an algebra in $\2C$, cf.\ \cite{ostrik}.


\subsection{Duality in tensor categories I}
\item If $G$ is a group and $\pi$ a representation on a finite dimensional vector space $V$, we
  define the `dual' or `conjugate' representation $\ol{\pi}$ on the dual vector space $V^*$ by 
$\langle\ol{\pi}(g)\phi,x\rangle=\langle\phi,\pi(g)x\rangle$. Denoting by $\pi_0$ the trivial
  representation, one finds $\Hom_{\Rep\,G}(\pi\otimes\ol{\pi},\pi_0)\cong\Hom_{\Rep\,G}(\pi,\pi)$,
  implying $\pi\otimes\ol{\pi}\succ\pi_0$. If $\pi$ is irreducible, then so is $\ol{\pi}$ and the
  multiplicity of $\pi_0$ in $\pi\otimes\ol{\pi}$ is one by Schur's lemma.

Since the above discussion is quite specific to the group situation, it clearly needs to be generalized.

\item Let $(\2C,\otimes,\11)$ be a strict tensor category and $X,Y\in\2C$. 
We say that $Y$ is a {\bf left dual} of $X$ if there are morphisms 
$e: Y\otimes X\rarr\11$ and $d:\11\rarr X\otimes Y$ satisfying
\[ \id_X\otimes e\mcirc d\otimes\id_X=\id_X, \quad\quad 
   e\otimes\id_{Y}\mcirc\id_{Y}\otimes d=\id_{Y}, \]
or, representing $e: Y\otimes X\rarr\11$ and $d:\11\rarr X\otimes Y$ by 
$\begin{tangle}\hh\coev\end{tangle}$ and $\begin{tangle}\hh\ev\end{tangle}$, respectively,
\[ \begin{tangle}
\object{X}\\
\hh\id\step\hcoev\obj{e}\\
\hh\id\hstep\obj{Y}\hstep\id\step\id\\
\hh\step[-.5]\obj{d}\hstep\hev\step\id\\
\Step\object{X}
\end{tangle}
\quad\quad=\quad\begin{tangle}\object{X}\\ \id\\ \object{X}\end{tangle}
\quad\quad\quad\quad\quad\quad
\begin{tangle}
\Step\object{Y}\\
\hh\step[-.8]\obj{e}\step[.8]\hcoev\step\id\\
\hh\id\hstep\obj{X}\hstep\id\step\id\\
\hh\id\step\hev\obj{d}\\
\object{Y}
\end{tangle}
\quad=\quad\begin{tangle}\object{Y}\\ \id\\ \object{Y}\end{tangle}
\]
($e$ stands for `evaluation' and $d$ for `dual'.). In this situation, $X$ is called a {\bf right dual}
of $Y$. 

Example: $\2C=\Vect_k^{\mathrm{fin}},\ X\in\2C$. Let $Y=X^*$, the dual vector space. Then
$e:Y\otimes X\rarr\11$ is the usual pairing. With the canonical isomorphism 
$f: X^*\otimes X\stackrel{\cong}{\longrightarrow}\End\,X$, we have $d=f^{-1}(\id_X)$.

We state some facts: 
\begin{enumerate}
\item Whether an object $X$ admits a left or right dual is not for us to choose. It is a
property of the tensor category.
\item If $Y,Y'$ are left (or right) duals of $X$ then $Y\cong Y'$. 
\item If ${}^\vee\! A,\ {}^\vee\! B$ are left duals of $A,B$, respectively, then ${}^\vee\! B\otimes{}^\vee\! A$ is
  a left   dual for $A\otimes B$, and similarly for right duals.
\item If $X$ has a left dual $Y$ and a right dual $Z$, we may or may not have $Y\cong Z$ ! (Again,
  that is a property of $X$.)
\end{enumerate}

While duals, if they exist, are unique up to isomorphisms, it is often convenient to make choices.
One therefore defines a {\bf left duality} of a strict tensor category $(\2C,\otimes,\11)$ to be a
map that assigns to each object $X$ a left dual ${}^\vee\! X$ and morphisms $e_X: {}^\vee\! X\otimes X\rarr\11$ and 
$d_X:\11\rarr X\otimes {}^\vee\! X$ satisfying the above identities.

Given a left duality and a morphism, $s:X\rarr Y$ we define 
\[ {}^\vee\! s= 
  e_Y\otimes\id_{{}^\vee\! X}\mcirc\id_{{}^\vee\! Y}\otimes s\otimes\id_{{}^\vee\! X}\mcirc\id_{{}^\vee\! Y}\otimes d_X
=\quad\quad\begin{tangle}
\Step\object{{}^\vee\! X}\\
\hh\step[-.8]\obj{e_Y}\step[.8]\hcoev\step\id\\
\id\step\O{s}\step\id\\
\hh\id\step\hev\obj{d_X}\\
\object{{}^\vee\! Y}
\end{tangle}
\]
Then $(X\mapsto{}^\vee\! X,\ \ s\mapsto{}^\vee\! s)$ is a contravariant functor. (We cannot recover the
e's and d's from the functor!) It can be equipped with a natural (anti-)monoidal isomorphism 
${}^\vee\!(A\otimes B)\rarr{}^\vee\! B\otimes{}^\vee\! A,\ {}^\vee\!\11\rarr\11$. Often, the duality functor
comes with a given anti-monoidal structure, e.g.\ in the case of pivotal categories, cf.\ Section \ref{sec-3}.

\item A chosen {\bf right duality}
$X\mapsto(X^\vee, e'_X:X\otimes X^\vee\rarr\11,\ d'_X:\11\rarr  X^\vee\otimes X)$
also give rise to a contravariant anti-monoidal functor $X\mapsto X^\vee$.

\item Categories equipped with a left (right) duality are called left (right) {\bf rigid} (or {\bf
  autonomous}). Categories with left and right duality are called rigid (or autonomous).

\item Examples: $\Vect_k^{\mathrm{fin}}, \Rep\,G$ are rigid.

\item Notice that ${}^{\vee\vee}X\cong X$ holds if and only if ${}^\vee\! X\cong X^\vee$, for which
  there is no general reason. 

\item If every object $X\in\2C$ admits a left dual ${}^\vee\! X$ and a right dual $X^\vee$, {\it and
  both are isomorphic}, we say that $\2C$ has {\bf two-sided duals} and write $\ol{X}$. 
We will only consider such categories, but we will  need stronger axioms.

\item Let $\2C$ be a $*$-category (cf.\ below) with left duality. If $({}^\vee\! X,e_X,d_X)$ is a left
dual of $X\in\2C$ then $(X^\vee={}^\vee\! X, d_X^*,e_X^*)$ is a right dual. Thus duals in
$*$-categories are automatically two-sided. For this reason, duals in $*$-category are often
axiomatized in a symmetric fashion by saying that a {\bf conjugate}, cf.\ \cite{DR,lro}, of an
object $X$ is a triple $(\ol{X},r,\ol{r})$, where 
$r:\11\rarr\ol{X}\otimes X,\ \ol{r}:\11\rarr X\otimes\ol{X}$ satisfy 
\[ \id_X\otimes r^*\mcirc \ol{r}\otimes\id_X=\id_X,\quad\quad\quad
   \id_{\ol{X}}\otimes\ol{r}^*\mcirc r\otimes\id_{\ol{X}}=\id_{\ol{X}}. \]
It is clear that then $(\ol{X},r^*,\ol{r})$ is a left dual and $(\ol{X},\ol{r}^*,r)$ a right dual.

\item Unfortunately, there is an almost Babylonian inflation of slightly different notions
  concerning duals, in particular when braidings are involved: A category can be rigid, 
autonomous, sovereign, pivotal, spherical, ribbon, tortile, balanced, closed, category with
conjugates, etc. To make   things worse, these terms are not always used in the same way!

\item Before we continue the discussion of  duality in tensor categories, we will discuss
  symmetries. For symmetric tensor categories, the discussion of duality is somewhat simpler than in
  the general case. Proceeding like this seems justified since symmetric (tensor) categories already
  appeared in the second paper (\cite{macl2} 1963) on tensor categories.


\subsection{Additive, linear and $*$-structure}\label{ss-ab}
\item The discussion so far is quite general, but often one encounters categories with more structure.
\item We begin with `Ab-categories' (=categories `enriched over abelian groups'): For such a
  category, each $\Hom(X,Y)$   is an abelian   group, and $\circ$ is bi-additive, cf.\ \cite[Section
    I.8]{cwm}. Example:  The category {\bf Ab} of abelian groups. In $\otimes$-categories, also
  $\otimes$ must be bi-additive on the morphisms. Functors of 
  Ab-tensor categories required to be additive on hom-sets.
\item If $X,Y,Z$ are objects in an Ab-category, $Z$ is called a direct sum of $X$ and $y$ if there
  are morphisms 
$X\stackrel{u}{\rightarrow}Z\stackrel{u'}{\rightarrow}X, Y\stackrel{v}{\rightarrow}Z\stackrel{v'}{\rightarrow}Y$
satisfying $u\circ u'+v\circ v'=\id_Z, u'\circ u=\id_X, v'\circ v=\id_Y$.  An {\bf additive}
category is an Ab-category having direct sums for all pairs of objects and a zero object.

\item An {\bf abelian category} is an additive category where every morphism has a kernel and a
cokernel and every monic (epic) is a kernel (cokernel). We do not have the space to go further into
this and must refer to the literature, e.g.\ \cite{cwm}. 

\item A category is said to have splitting idempotents (or is `Karoubian') if $p=p\circ p\in\End\,X$
implies the existence of an object $Y$ and of morphisms $u:Y\rarr X,\ u':X\rarr Y$ such that
$u'\circ u=\id_Y$ and $u\circ u'=p$. An additive category with splitting idempotents is called {\bf
pseudo-abelian}. Every abelian category is pseudo-abelian.

\item In an abelian category with duals, the functors $-\otimes X$ and $X\otimes -$ are
  automatically exact, cf.\ \cite[Proposition 1.16]{DM}. But without rigidity this is far from true.

\item A {\bf semisimple} category is an abelian category where every short exact sequence splits.

An alternative, and more pedestrian, way to define semisimple categories is as pseudo-abelian
categories admitting a family of simple objects $X_i,\ i\in I$ such  that every
$X\in\2C$ is a finite direct sum of $X_i$'s.

Standard examples: The category $\Rep\,G$ of finite dimensional representations of a compact group 
$G$, the category $H-\Mod$ of finite dimensional left modules for a finite dimensional semisimple
Hopf algebra $H$.  

\item In {\bf $k$-linear} categories, each $\Hom(X,Y)$ is $k$-vector space (often required finite
  dimensional), and $\circ$ (and $\otimes$ in the monoidal case) is bilinear. Functors must be
  $k$-linear. Example: $\Vect_k$. 
\item Pseudo-abelian categories that are $k$-linear with finite-dimensional hom-sets are called {\bf
Krull-Schmidt} categories. (This is slightly weaker than semisimplicity.)

\item A {\bf fusion category} is a semisimple $k$-linear category with finite dimensional hom-sets,
finitely many isomorphism classes of simple objects and $\End\,\11=k$. We also require that $\2C$
has 2-sided duals. 

\item A {\bf finite tensor category} (Etingof, Ostrik \cite{EO1}) is a $k$-linear tensor category with
  $\End\,\11=k$ that is equivalent (as a category) to the category of modules over a finite
dimensional k-algebra. (There is a more intrinsic characterization.) Notice that semisimplicity is
not assumed. 

\item Dropping the condition $\End\,\11=k\id_k$, one arrives at a {\bf multi-fusion category}
(Etingof, Nikshych, Ostrik \cite{eno}). 
\item Despite the recent work on generalizations, most of these lectures will be concerned with
  \underline{semisimple $k$-linear categories satisfying $\End\,\11=k\,\id_\11$}, including infinite
  ones! (But see the remarks at the end of this section.)

\item If $\2C$ is a semisimple tensor category, one can choose representers $\{X_i, i\in I\}$ of the
  simple isomorphism classes and define $N_{i,j}^k\in\7Z_+$ by
\[ X_i\otimes X_j\cong\bigoplus_{k\in I} N_{i,j}^k X_k. \]
There is a distinguished element $0\in I$ such that $X_0\cong I$, thus
$N_{i,0}^k=N_{0,i}^k=\delta_{i,k}$. By associativity of $\otimes$ (up to isomorphism)
\[ \sum_n  N_{i,j}^n N_{n,k}^l = \sum_m  N_{i,m}^l N_{j,k}^m\quad\forall i,j,k,l\in I.\] 
If $\2C$ has two-sided duals, there is an involution $i\mapsto\ol{\imath}$ such that 
$\ol{X_i}\cong X_{\ol{\imath}}$. One has $N_{i,j}^0=\delta_{i,\ol{\jmath}}$.
The quadruple $(I,\{N_{i,j}^k\},0,i\mapsto\ol{\imath})$ is called the {\bf fusion ring} or 
{\bf fusion hypergroup} of $\2C$. 

\item The above does not work when $\2C$ is not semisimple. But: In any abelian tensor category, one
  can consider the {\bf Grothendieck  ring} $R(\2C)$, the free abelian 
  group generated by the isomorphism classes $[X]$ of objects in $\2C$, with a relation
  $[X]+[Z]=[Y]$ for every short exact sequence $0\rarr X\rarr Y\rarr Z\rarr 0$ and
  $[X]\cdot[Y]=[X\otimes Y]$. 

In the semisimple case, the Grothendieck ring has $\{[X_i], i\in I\}$
  as $\7Z$-basis and $[X_i]\cdot[X_j]=\sum_k N_{i,j}^k[X_k]$. Obviously, an isomorphism of
  hypergroups gives rise to a ring isomorphism of Grothendieck rings, but the converse is not
  obvious. While the author is not aware of counterexamples, in order to rule out this annoying
  eventuality, some authors work with isomorphisms of the Grothendieck {\bf semi}ring or the
{\bf ordered} Grothendieck ring, cf e.g.\ \cite{hand}.

Back to hypergroups:

\item The hypergroup contains important information about a tensor category, but it misses that
  encoded in the associativity constraint. In fact, the hypergroup of $\Rep\,G$ for a finite group
  $G$ contains exactly the   same information as the character table of $G$, and it is well known
  that there are non-isomorphic finite groups with isomorphic character tables. (The simplest
  example is given by the dihedral group 
$D_8=\7Z_4\rtimes\7Z_2$ and the quaternion group $Q$, cf.\ any elementary textbook, e.g.\ \cite{JL}.)
Since $D_8$ and $Q$ have the same number of irreducible representations, the categories $\Rep\,D_8$
and $\Rep\,Q$ are equivalent (as categories). They are not equivalent as symmetric tensor
categories, since this would imply $D_8\cong Q$ by the duality theorems of Doplicher and Roberts
\cite{DR} or Deligne \cite{del} (which we will discuss in Section \ref{sec-3}). In fact, $D_8$ and
$Q$ are already inequivalent as 
tensor categories (i.e.\ they are not isocategorical in the sense discussed below). Cf.\ \cite{TY},
where fusion categories with the fusion hypergroup of $D_8$ are classified (among other things).

\item On the positive side: (1) If a finite group $G$ has the same fusion hypergroup (or character
  table) as a finite simple group $G'$, then $G\cong G'$, cf.\ \cite{chen}. (The proof uses the
classification of finite simple groups.) (2) Compact groups that are abelian or connected are
determined by their fusion rings (by Pontrjagin duality, respectively by a result of McMullen
\cite{mcm} and Handelman \cite{hand}. The latter is first proven for simple compact Lie groups and
then one deduces the general result via the structure theorem for connected compact groups.)

\item If all objects in a semisimple category $\2C$ are invertible, the fusion hypergroup becomes a
  group. Such fusion categories are called {\bf pointed} and are just the linear versions of the
  categorical groups encountered earlier. This situation is very special, but:

\item To each hypergroup $\{I,N,0,i\mapsto\ol{\imath}\}$ one can associate a group $G(I)$ as
  follows: Let $\sim$ be the smallest equivalence relation on $I$ such that
\[ i\sim j\quad\mbox{whenever}\quad \exists\, m,n\in I:\ 
     i\prec mn\succ j\ \  (\mbox{i.e.}\ N_{n,m}^i\ne 0\ne N_{n,m}^j).\]
Now let $G(I)=I/\!\sim$ and define
\[ [i]\cdot[j]=[k]\ \ \mbox{for any}\ k\prec ij, \quad\quad\quad[i]^{-1}=[\ol{\imath}],
\quad\quad e=[0]. \]
Then $G(I)$ is a group, and it has the universal property that every map $p:I\rarr K$, $K$ a group, 
satisfying $p(k)=p(i)p(j)$ when $k\prec ij$ factors through the map $I\rarr G(I),\ i\mapsto[i]$. 

In analogy to the abelianization of a non-abelian group, $G(I)$ should perhaps be called the {\bf
groupification} of the hypergroup $I$. But it was called the {\bf universal grading group} by
Gelaki/Nikshych \cite{GN1}, to which this is due in the above generality, since every group-grading on
the objects of a fusion category having fusion hypergroup $I$ factors through the map $I\rarr G(I)$.

\item In the symmetric case (where $I$ and $G(I)$ are abelian, but everything else as above)
this groupification is due to Baumg\"artel/Lled\'o \cite{BL}, who spoke of the `chain group'. They
stated the conjecture that if $K$ is a compact group, then the (discrete) universal grading group
$G(\Rep\,K)$ of $\Rep\,K$ is the Pontrjagin dual of the (compact) center $Z(K)$.  Thus: The {\it
  center of a compact group $K$ can be recovered from the fusion ring} of $K$, even if $K$ itself in
general cannot!  This conjecture was proven in \cite{mue14}, but the whole circle of ideas is
already implicit in \cite{mcm}.

Example: The representations of $K=SU(2)$ are labelled by $\7Z_+$ with 
\[ i\otimes j=|i-j|\ \oplus\ \cdots\ \oplus\ i+j-2 \ \oplus\  i+j.\] 
From this one easily sees that there are two $\sim$-equivalence classes, consisting of the even and
odd integers. This is compatible with $Z(SU(2))=\7Z/2\7Z$. Cf.\ \cite{BL}.

\item There is another application of $G(\2C)$: If $\2C$ is $k$-linear semisimple then group of
  natural monoidal isomorphisms of $\id_\2C$ is given by $\Aut_\otimes(\id_\2C)\cong\Hom(G(\2C),k^*)$. 

\item Given a fusion category $\2C$ (where we have two-sided duals $\ol{X}$), Gelaki/Nikshych
\cite{GN1} define the full subcategory $\2C_{ad}\subset\2C$ to be the generated by the objects
$X\otimes\ol{X}$ where $X$ runs through the simple objects.  
Notice that $\2C_{ad}$ is just the full subcategory of objects of universal grading zero.

Example: If $G$ is a compact group then $(\Rep\,G)_{ad}=\Rep(G/Z(G))$.

A fusion category $\2C$ is called {\bf nilpotent} \cite{GN1} when its upper central series
\[ \2C\supset\2C_{ad}\supset(\2C_{ad})_{ad}\supset\cdots \]
leads to the trivial category after finitely many steps.

Example: If $G$ is a finite group then $\Rep\,G$ is nilpotent if and only if $G$ is nilpotent.

\item We call a square $n\times n$-matrix $A$ indecomposable if there is no proper subset
  $S\subset\{1,\ldots,n\}$ such that $A$ maps the coordinate subspace span$\{ e_s\ | \ s\in S\}$ into
  itself. Let $A$ be an indecomposable square matrix $A$ with non-negative entries and eigenvalues
$\lambda_i$. Then the theorem of Perron and Frobenius states that there is a unique non-negative
  eigenvalue $\lambda$, the Perron-Frobenius eigenvalue, such that
  $\lambda=\max_i|\lambda_i|$. Furthermore, the associated eigenspace is one-dimensional and
  contains a vector with all components non-negative. Now, given a finite hypergroup 
$(I,\{N_{i,j}^k\},0,i\mapsto\ol{\imath})$ and $i\in I$, define $N_i\in\mathrm{Mat}(|I|\times|I|)$ by
$(N_i)_{jk}=N_{i,j}^k$. Due to the existence of duals, this matrix is indecomposable. Now the {\bf
  Perron-Frobenius dimension} $d_{FP}(i)$ of $i\in I$ is defined as the Perron-Frobenius eigenvalue
of $N_i$. Cf.\ e.g.\ \cite[Section 3.2]{FK}. Then: 
\[ d_{FP}(i)d_{FP}(j)=\sum_k N_{i,j}^k d_{FP}(k). \]
Also the hypergroup $I$ has a Perron-Frobenius dimension: $FP-\dim(I)=\sum_i d_{FP}(i)^2$. This also
defines the PF-dimension of a fusion category, cf.\ \cite{eno}

\item {\bf Ocneanu rigidity}: Up to equivalence there are only finitely many fusion   categories
with given fusion hypergroup. The general statement was announced by Blanchard/A.\ Wassermann, and a
proof is given in \cite{eno}, using the deformation cohomology theory of Davydov \cite{dav0} and
Yetter \cite{yetter}. 

\item Ocneanu rigidity was preceded and motivated by several related results on Hopf algebras:
  Stefan \cite{stefan} proved that the number of isomorphism classes of semisimple and co-semisimple
  Hopf algebras of given finite dimension is finite. For Hopf $*$-algebras, Blanchard
  \cite{blanchard} even proved a bound on the number of iso-classes in terms of the dimension. There
  also is an upper bound on the number of iso-classes of semisimple Hopf algebras with given number
  of irreducible representations, cf.\ Etingof's appendix to \cite{ostrik3}.

\item There is an enormous literature on hypergroups. Much of this concerns harmonic analysis on the
latter and is not too relevant to tensor categories. But the notion of amenability of hypergroups
does have such applications, cf.\ e.g.\ \cite{hiaiizumi}. For a review of some aspects of hypergroups, in
particular the discrete ones relevant here, cf.\ \cite{wild}.

\item A considerable fraction of the literature on tensor categories is devoted to categories that are
$k$-linear over a field $k$ with finite dimensional Hom-spaces. Clearly this a rather restrictive
condition. It is therefore very remarkable that $k$-linearity can actually be deduced under suitable
assumptions, cf.\ \cite{kup3}.

\item {\bf $*$-categories}: A `$*$-operation' on a $\7C$-linear category $\2C$ is a contravariant
  functor $*:\2C\rarr\2C$ which acts trivially on the objects,  is antilinear, involutive ($s^{**}=s$)
 and monoidal $(s\otimes t)^*=s^*\otimes t^*$ (when $\2C$ is monoidal).
A $*$-operation is called positive if $s^*\circ s=0$ implies $s=0$.
Categories with (positive) $*$-operation are also called hermitian (unitary). We will use
`$*$-category' as a synonym for `unitary category'.) Example: The category of Hilbert spaces
$\2H\2I\2L\2B$ with bounded linear maps and $*$ given by the adjoint.

\item It is easy to prove that a finite dimensional $\7C$-algebra with positive $*$-operation is
  semisimple. Therefore, a unitary category with finite dimensional hom-sets has semisimple
  endomorphism algebras. If it has direct sums and splitting idempotents then it is semisimple.
\item  Banach-, $C^*$- and von Neumann categories: A Banach category \cite{karoubi} is a
  $\7C$-linear additive category, where each $\Hom(X,Y)$ is a Banach space, and the norms satisfy 
\[  \|s\circ t\|\le \|s\|\,\|t\|, \quad\quad\| s^*\circ s\|=\|s\|^2. \]
(They were introduced by Karoubi with a view to applications in K-theory, cf.\ \cite{karoubi}.)
A Banach $*$-category is a Banach category with a positive $*$-operation. A $C^*$-category is a
Banach $*$-category satisfying $\|s^*\circ s\|=\|s\|$ for any morphism $s$ (not only endomorphisms).
In a $C^*$-category, each $\End(X)$ is a $C^*$-algebra. Just as an additive category is a `ring with
several objects', a $C^*$-category is a ``$C^*$-algebra with several objects''. Von Neumann
categories are defined similarly, cf. \cite{glr}. They turned out to have applications to
$L^2$-cohomology (cf.\ Farber \cite{farber}), representation theory of quantum groups (Woronowicz
\cite{woro2}), subfactors \cite{lro}, etc.

Remark: A $*$-category with finite dimensional hom-spaces and $\End\,\11=\7C$ automatically is a
$C^*$-category in a unique way. (Cf.\ \cite{mue06}.)

\item If $\2C$ is a $C^*$-tensor category, $\End\,\11$ is a commutative $C^*$-algebra, thus
  $\cong C(S)$ for some compact Hausdorff space $S$. Under certain technical conditions, the spaces
  $\Hom(X,Y)$ can be considered as vector bundles over $S$, or at least as (semi)continuous fields
  of vector spaces. (Work by Zito \cite{zito} and Vasselli \cite{vasselli1}.)
In the case where $\End\,\11$ is finite dimensional, this boils down to a direct
  sum decomposition of $\2C=\oplus_i \2C_i$, where each $\2C_i$ is a tensor category with
  $\End_{\2C_i}(\11_{\2C_i})=\7C$. (In this connection, cf.\ Baez' comments a Doplicher-Roberts type
theorem for finite groupoids \cite{baez2}.)

\end{itemize}


\section{Symmetric tensor categories}\label{sec-2}
\begin{itemize}
\item Many of the obvious examples of tensor categories encountered in Section \ref{sec-1}, like the
  categories $\2S\2E\2T$, Vect$_k$, representation categories of groups and Cartesian categories
  (tensor product $\otimes$ given by the categorical product), have an additional piece of
  structure, to which this section is dedicated.

\item A {\bf symmetry} on a tensor category $(\2C,\otimes,\11,\alpha,\rho,\lambda)$ is a natural
  isomorphism $c: \otimes\rarr\otimes\circ\sigma$, where $\sigma:\2C\times\2C\rarr\2C\times\2C$ is
  the flip automorphism of $\2C\times\2C$, such that $c^2=\id$. (I.e., for any two objects $X,Y$
  there is an isomorphism $c_{X,Y}:X\otimes Y\rarr Y\otimes X$, natural w.r.t.\ $X,Y$ such that 
$c_{Y,X}\circ c_{X,Y}=\id_{X\otimes Y}$.), where ``all properly built diagrams commute'', i.e.\ the
category is coherent. A {\bf symmetric tensor category} (STC) is a tensor category equipped with a symmetry.

We represent the symmetry graphically by 
\[ c_{X,Y}=\quad\begin{tangle} \object{Y}\Step\object{X}\\ \X\\\object{X}\Step\object{Y}\end{tangle}\]
\item As for tensor categories, there are two versions of the Coherence Theorem. Version I (Mac Lane
  \cite{macl2}): Let $(\2C,\otimes,\11,\alpha,\rho,\lambda)$ be a tensor category. Then a natural
  isomorphism $c: \otimes\rarr\otimes\circ\sigma$ satisfying $c^2=\id$ is a symmetry if and only if 
\begin{diagram}
 (X\otimes Y)\otimes Z & \rTo^{c_{X,Y}\otimes\id_Z} & (Y\otimes X)\otimes Z & \rTo^{\alpha_{Y,X,Z}} &
   Y\otimes (X\otimes Z) & \rTo^{\id_Y\otimes c_{X,Z}} & Y\otimes (Z\otimes X) \\
\dTo^{\alpha_{X,Y,Z}} &&&&&& \uTo_{\alpha_{Y,Z,X}} \\
X\otimes(Y\otimes Z) &&& \rTo_{c_{X,Y\otimes Z}} &&& (Y\otimes Z)\otimes X
\end{diagram}
commutes. (In the strict case, this reduces to 
$c_{X,Y\otimes Z}=\id_Y\otimes c_{X,Z}\circ c_{X,Y}\otimes\id_Z$.)

A {\bf symmetric tensor functor} is a tensor functor $F$ such that 
$F(c_{X,Y})=c'_{F(X),F(Y)}$. Notice that a  natural transformation between symmetric tensor functors is just
a monoidal natural transformation, i.e.\ there is no new condition.
\item Now we can state version II of the Coherence theorem: Every symmetric tensor category is
  equivalent (by a symmetric tensor functor) to a strict one.
\item Examples of symmetric tensor categories:
\begin{itemize}
\item The category $\7S$ defined earlier, when $c_{n,m}:n+m\rarr n+m$ is taken to be the element of
  $S_{n+m}$  defined by $(1,\ldots,n+m)\mapsto(n+1,\ldots,n+m,1,\ldots n)$. It is the free symmetric
  monoidal   category generated by one object.
\item Non-strict symmetric categorical groups were classified by Sinh \cite{sinh}. We postpone our
  discussion to Section \ref{sec-4}, where we will also consider the braided case.
\item Vect$_k$, representation categories of groups: We have the canonical symmetry
$c_{X,Y}: X\otimes Y\rarr Y\otimes X,\ x\otimes  y\mapsto y\otimes   x$. 
\item The tensor categories obtained using products or coproducts are symmetric.
\end{itemize}

\item Let $\2C$ be a strict STC, $X\in\2C$ and $n\in\7N$. Then there is a unique homomorphism
\[ \Pi_n^X: S_n\rarr\Aut\,X^{\otimes n}\quad\quad\mbox{such that}\quad\quad
  \sigma_i\mapsto\id_{X^{\otimes (i-1)}}\otimes c_{X,X}\otimes\id_{X^{\otimes (n-i-1)}}. \]
Proof: This is immediate by the definition of STCs and the presentation 
\[ S_n=\{ \sigma_1,\ldots, \sigma_{n-1}\ | \ \sigma_i\sigma_{i+1}\sigma_i=\sigma_{i+1}\sigma_i\sigma_{i+1},\
    \sigma_i\sigma_j=\sigma_j\sigma_i \ \ \mbox{when}\ \ |i-j|>1,\ \sigma_i^2=1\} \]
of the symmetric groups.

These homomorphisms in fact combine to a symmetric tensor functor $F:\7S\rarr\2C$ such that
$F(n)=X^{\otimes n}$. (This is why $\7S$ is called the free symmetric tensor category on one generator.)
\item In the $\otimes$-category $\2C=\Vect_k^{\mathrm{fin}}$, $\Hom(V,W)$ is itself an object of
$\2C$, giving rise to an internal hom-functor: 
$\2C^\op\times\2C\rarr\2C,\ X\times  Y\mapsto[X,Y]=\underline{\Hom}(X,Y)$ 
  satisfying some axioms. In the older literature, a symmetric tensor category with such an
  internal-hom functor is called a {\bf closed} category. There are coherence theorems for closed
  categories. \cite{kelml,kel3}.

Since in $\Vect_k^{\mathrm{fin}}$ we have  $\Hom(V,W)\cong V^*\otimes W$, it is sufficient -- and
more transparent -- to axiomatize duals $V\mapsto V^*$, as is customary in the more recent
literature. We won't mention `closed' categories again. (Which doesn't mean that they have no uses!)

\item We have seen that, even if a tensor category has left and right duals 
${}^\vee\! X,\ X^\vee$ for every object, they don't need to be isomorphic.
But if $\2C$ is symmetric and $X\mapsto ({}^\vee\! X,e_X,d_X)$ is a left duality, then
defining 
\[ X^\vee={}^\vee\! X,\quad e'_X=e_X\circ c_{X,{}^\vee\! X},\quad d'_X=c_{X,{}^\vee\! X}\circ d_X,\]
one easily checks that $X\mapsto(X^\vee,e'_X,d'_X)$ defines a right duality. We can thus take
${}^\vee\! X=X^\vee$ and denote this more symmetrically by $\ol{X}$. 

\item Let $\2C$ be symmetric with given left duals and with right duals as just defined, and let
  $X\in\2C$. Define the (left) {\bf trace} $\Tr_X:\End\,X\rarr\End\,\11$ by
\[ \Tr_X(s)=e_X\mcirc\id_{\ol{X}}\otimes s\mcirc d'_X=
\quad\quad\begin{tangle}
\coev\mobj{e_X}\\
\step[-.7]\obj{\ol{X}}\step[.7]\id\Step\O{s}\\
\ev\step[-.4]\mobj{d'_X}\end{tangle} 
\quad\quad\quad\left(=\quad\quad
\begin{tangle}
\coev\mobj{e_X}\\
\step[-.7]\obj{\ol{X}}\step[.7]\id\Step\O{s}\\
\X\\
\ev\step[-.4]\mobj{d_X}\end{tangle} 
\quad\quad\right)
\]
Without much effort, one can prove the trace property $\Tr_X(ab)=\Tr_X(ba)$ and multiplicativity
under $\otimes: \Tr_{X\otimes Y}(a\otimes b)=\Tr_X(a)\Tr_Y(b)$.   
Finally, $\Tr_X$ equals the right trace defined using $e'_X, d_X$. 
For more on traces in tensor categories cf.\ e.g.\ \cite{js4, mal2}.

\item Using the above, we define the categorical {\bf dimension} of an object $X$ by
$d(X)=\Tr_X(\id_X)\ \in\End\,\11$. If $\End\,\11=k\id_\11$, we can use this identification
  to obtain $d(X)\in k$. 

With this dimension and the usual symmetry and duality on $\Vect_k^{\mathrm{fin}}$, one verifies
$d(V)=\dim_kV\cdot 1_k$.  

However, in the category SVect$_k$ of super vector spaces (which coincides with the representation
category $\Rep_k\7Z_2$, but has the symmetry modified by the Koszul rule) it gives the
super-dimension, which can be negative, while one might prefer the total dimension. Such situations
can be taken care of (without changing the symmetry) by introducing twists.

\item If $(\2C,\otimes,\11)$ is strict symmetric, we define a {\bf twist} to be natural family $\{
  \Theta_X \in\End\,X,\ X\in\2C\}$ of isomorphisms  satisfying 
\be \Theta_{X\otimes Y}=\Theta_X\otimes\Theta_Y, \quad\quad\Theta_\11=\id_\11\label{3}\ee
i.e., $\Theta$ is a monoidal natural isomorphism of the functor $\id_\2C$. If $\2C$ has a left
duality, we also require
\[ {}^\vee\!(\Theta_X)=\Theta_{{}^\vee\! X}. \]
The second condition implies $\Theta_X^2=\id$. Notice that $\Theta_X=\id_X\ \forall X$ is a legal 
choice. This will not remain true in braided tensor categories! 

Example: If $G$ is a compact group and $\2C=\Rep\,G$, then the twists $\Theta$ satisfying only
(\ref{3}) are in bijection with $Z(G)$. The second condition reduces this to central elements of order
two. (Cf.\ e.g.\ \cite{mue16}.)

\item Given a strict symmetric tensor category with left duality and a twist, we can define a right
  duality by $X^\vee={}^\vee\! X$, writing $\ol{X}={}^\vee\! X=X^\vee$, but now
\be e'_X=e_X\mcirc c_{X,\ol{X}}\mcirc\Theta_X\otimes\id_{\ol{X}},
   \quad\quad d'_X=\id_{\ol{X}}\otimes\Theta_X\mcirc c_{X,\ol{X}}\circ d_X,\label{eq-e'}\ee
still defining a right duality and the maps $Tr_X:\End\,X\rarr\End\,\11$ still are traces.

\item Conversely, the twist can be recovered from $X\mapsto(\ol{X},e_X,d_X,e'_X,d'_X)$ by
\[ \Theta_X=(\Tr_X\otimes\id)(c_{X,X})=
\quad\quad\begin{tangle}
\Step\object{X}\\
\hh\step[-1]\obj{e_X}\step\hcoev\step\id\\
\step[-.7]\obj{\ol{X}}\step[.7]\id\step\hX\\
\hh\step[-1]\obj{d'_X}\step\hev\step\id\\
\Step\object{X}\end{tangle}
\]
Thus: Given a symmetric tensor category with fixed left duality, every twist gives rise to a right
duality, and every right duality that is `compatible' with the left duality gives a twist. (The
trivial twist $\Theta\equiv\id$ corresponds to the original definition of right duality. The latter
does not work in proper braided categories!)
This compatibility makes sense even for categories without symmetry (or braiding) and will be
discussed later ($\leadsto$ pivotal categories).

\item The symmetric categories with  $\Theta\equiv\id$ are now called {\bf even}.

\item The category SVect$_k$ of super vector spaces with $\Theta$ defined in terms of the
  $\7Z_2$-grading now satisfies  $\dim(V)\ge 0$ for all $V$.

\item The standard examples for STCs are $\Vect_k, S\Vect_k, \Rep\,G$ and the representation
  categories of supergroups.
In fact, rigid STCs are not far from being representation categories of (super)groups. However, they
not always are, cf.\ \cite{GK2} for examples of non-Tannakian symmetric categories.)

\item A category $\2C$ is called {\bf concrete} if its objects are sets and
$\Hom_\2C(X,Y)\subset\Hom_{\mathrm{Sets}}(X,Y)$. A $k$-linear category is called concrete if the
objects are fin.dim.\ vector spaces over $k$ and $\Hom_\2C(X,Y)\subset\Hom_{\Vect_k}(X,Y)$. 
However, a better way of thinking of a concrete category is as a (abstract) category
  $\2C$ equipped with a {\bf fiber functor}, i.e.\ a faithful functor $E:\2C\rarr Sets$,
  respectively $E:\2C\rarr\Vect_k$. The latter is required to be monoidal when $\2C$ is monoidal. 
\item Example: $G$ a group. Then $\2C:=\Rep_kG$ should be considered as an abstract $k$-linear
$\otimes$-category together with a faithful $\otimes$-functor $E:\2C\rarr\Vect_k$.
\item The point of this that a category $\2C$  may have inequivalent fiber functors!!
\item But: If $k$ is algebraically closed of characteristic zero, $\2C$ is rigid {\it symmetric}
  $k$-linear with $\End\,\11=k$ and $F,F'$ are {\it symmetric} fiber functors then $F\cong F'$ (as
  $\otimes$-functors). (Saavedra Rivano \cite{sr,DM}). 
\item The first non-trivial application of (symmetric) tensor categories probably were the
  reconstruction theorems of Tannaka \cite{tannaka} (1939!) and Saavedra Rivano \cite{sr,DM}.

Let $k$ be algebraically closed.
Let $\2C$ be rigid symmetric $k$-linear with $\End\,\11=k$ and $E:\2C\rarr\Vect_k$ faithful tensor
functor.  (Tannaka did this for $k=\7C$, $\2C$ a $*$-category and $E$ $*$-preserving.)  Let
$G=\Aut_\otimes E$ be the group of natural monoidal [unitary] automorphisms of $E$. Define a functor
$F:\2C\rarr\Rep\,G$ [unitary representations] by 
\[ F(X)=(E(X),\pi_X),\quad\quad \pi_X(g)=g_X\quad\quad(g\in G). \]
Then $G$ is pro-algebraic [compact] and $F$ is an equivalence of symmetric tensor [$*$-]categories. 

Proof: The idea is the following (Grothendieck, Saavedra Rivano \cite{sr}, cf.\ Bichon
\cite{bichon}): Let $E_1,E_2:\2C\rarr\Vect_k$ be fiber functors. Define a unital $k$-algebra
$A_0(E_1,E_2)$ by  
\[ A_0(E_1,E_2)= \bigoplus_{X\in\2C}\Hom_{\mathrm{Vect}}(E_2(X),E_1(X)), \]
spanned by elements $[X,s],\ X\in\2C, s\in\Hom(E_2(X),E_1(X))$, with
$[X,s]\cdot[Y,t]=[X\otimes Y, u]$, where $u$ is the composite
\[ \begin{diagram}          
 E_2(X\otimes Y) & \rTo^{(d^2_{X,Y})^{-1}} & E_2(X)\otimes E_2(Y) & \rTo^{s\otimes t} &
   E_1(X)\otimes E_1(Y) & \rTo^{d^1_{X,Y}} & E_1(X\otimes Y). \end{diagram} \]
This is a unital associative algebra, and $A(E_1,E_2)$ is defined as the quotient by the ideal
generated by the elements $[X,a\circ E_2(s)]-[Y,E_1(s)\circ a]$, where 
$s\in\Hom_\2C(X,Y),\ a\in\Hom_{\mathrm{Vect}}(E_2(Y),E_1(X))$.

\item Remark: Let $E_1,E_2:\2C\rarr\Vect_k$ be fiber functors as above. Then 
the map 
\[ X\times Y\ \mapsto \ \Hom_{\Vect_k}(E_2(X),E_1(Y))\]
extends to a functor $F: \2C^\op\times\2C\rarr\Vect_k$. Now the algebra $A(E_1,E_2)$ is just the
{\bf coend} $\int^X F(X,X)$ of $F$, a universal object. Coends are a categorical, non-linear
version of traces, but we refrain from going into them since it takes some time to appreciate the
concept. (Cf.\ \cite{cwm}.) 

\item Now one proves \cite{bichon,mue16}:
\begin{itemize}
\item If $E_1,E_2$ are {\it symmetric} tensor functors then $A(E_1,E_2)$ is commutative.
\item If $\2C$ is $*$-category and $E_1,E_2$ are *-preserving then $A(E_1,E_2)$ is a $*$-algebra and
has a $C^*$-completion. 
\item If $\2C$ is finitely generated (i.e.\ there exists a monoidal generator $Z\in\2C$ such that
  every $X\in\2C$ is direct summand of some $Z^{\otimes N}$) then $A(E_1,E_2)$ is finitely generated.
\item There is a bijection between natural monoidal (unitary) isomorphisms $\alpha:E_1\rarr E_2$ and
  ($*$-)characters on $A(E_1,E_2)$.
\end{itemize}

Thus: If $E_1,E_2$ are symmetric and either $\2C$ is finitely generated or a $*$-category, the
algebra $A(E_1,E_2)$ admits characters (by the Nullstellensatz or by Gelfand's theory), thus 
$E_1\cong E_2$. One also finds that
$G=\mathrm{Aut}_\otimes E\cong(*$-$)\mathrm{Char}(A(E,E))$ and $A(E)=\mathrm{Fun}(G)$
(representative respectively continuous functions). This is used to prove that $F:\2C\rarr\Rep\,G$
is an equivalence. 

\item Remarks: 1. While it has become customary to speak of Tannakian categories, the work of Kre\u\i n,
  cf.\ \cite{krein}, \cite[Section 30]{HR}, should also be mentioned since it can be considered as a
  precursor of the later generalizations to non-symmetric categories, in particular in Woronowicz's
  approach.

2. The uniqueness of the symmetric fiber functor $E$ implies that $G$ is unique up to  isomorphism.

3. For the above construction, we need to have a fiber functor. Around 1989, Doplicher and Roberts
  \cite{DR}, and independently Deligne \cite{del} construct such a functor under weak assumptions on
  $\2C$. See below.  

4. The uniqueness proof fails if either of $E_1,E_2$ is not symmetric (or $\2C$ is not
  symmetric). Given a group $G$, there is a tautological fiber functor $E$. The fact that there may
  be (non-symmetric) fiber functors that are not naturally isomorphic to $E$ reflects the fact that
  there can be groups $G'$ such that $\Rep\,G\simeq\Rep\,G'$ as tensor categories, but not as
  symmetric tensor categories! This phenomenon was independently discovered by Etingof/Gelaki
  \cite{EG}, who called such $G,G'$ {\bf isocategorical} and produced examples of isocategorical but
  non-isomorphic finite groups, by Davydov \cite{dav2} and by Izumi and Kosaki \cite{IK}.
The treatment in \cite{EG} relies on the fact that if $G,G'$ are isocategorical then
$\7CG'\cong\7C^J$ for some Drinfeld twist $J$. A more categorical approach, allowing also an
extension to compact groups, will be given in \cite{mue21}.
A group $G$ is called categorically rigid if every $G'$ isocategorical to $G$ is actually isomorphic
to $G$. (Compact groups that are abelian or connected are categorically rigid in a strong sense
since they are determined already by their fusion hypergroups.)

\item Consider the free rigid symmetric tensor $*$-category $\2C$ with $\End\,\11=\7C$ generated by
  one object $X$ of dimension $d$. If $d\in\7N$ then $\2C$ is equivalent to $\Rep\,U(d)$ or
  $\Rep\,O(d)$ or $\Rep\,Sp(d)$, depending on whether $X$ is non-selfdual or orthogonal or
  symplectic. The proof \cite{baez2} is straightforward once one has the Doplicher-Roberts theorem. 

\item The free rigid symmetric categories just mentioned can be constructed in a topological way,
  in a fashion very similar to the construction of the Temperley-Lieb categories TL$(\tau)$. The
  main difference is that one allows the lines in the pictures defining the morphisms to cross. (But
  they still live in a plane.) Now one quotients out the negligible morphisms and completes
  w.r.t.\ direct sums and splitting idempotents. (In the non-self dual case, the objects are words
  over the alphabet $\{+,-\}$ and the lines in the morphisms are directed.) All this is noted in
  passing by Deligne in a paper \cite{del2} dedicated to the exceptional groups!
 Notice that when $d\not\in\7N$, these categories are examples of rigid symmetric categories that
 are not Tannakian.

\item The above results already establish strong connections between tensor categories and representation
theory, but there is much more to say.

\end{itemize}


\section{Back to general tensor categories}\label{sec-3}
\begin{itemize}
\item In a general tensor category, left and right duals need not coincide. This can already be seen
for the left module category $H-\Mod$ of a Hopf algebra $H$. This category has left and right duals,
related to $S$ and $S^{-1}$. ($S$ must be invertible, but 
can be aperiodic!) They coincide when $S^2(x)=uxu^{-1}$ with $u\in H$.  

\item We only consider tensor categories that have isomorphic left and right duals, i.e.\
two-sided duals, which we denote $\ol{X}$.
\item If $\2C$ is $k$-linear with $\End\,\11=k\,\id$ and $\End\,X=k\,\id$ ($X$ is
  simple/irreducible), one can canonically define the squared dimensions $d^2(X)\in k$ by
\[ d^2(X)=(e_X\circ d'_X)\cdot(e'_X\circ d_X)\ \in\End\,\11. \]
(Since $X$ is simple, the morphisms $d,d',e,e'$ are unique up to scalars, and well-definedness of
$d^2$ follows from the equations involving $(d,e),(d',e')$ bilinearly.) Cf.\ \cite{mue09}.
\item If $\2C$ is a fusion category, we define its {\bf dimension} by $\dim\2C=\sum_i d^2(X_i)$. 
\item If $H$ is a finite dimensional semisimple and co-semisimple Hopf algebra then
  $\dim\,H\mathrm{-}\Mod=\dim_kH$. (A finite dimensional Hopf algebra is co-semisimple if and only
  if the dual Hopf algebra $\widehat{H}$ is semisimple.) 
\item Even if $\2C$ is semisimple, it is not clear whether one can choose roots $d(X)$ of the above
numbers  $d^2(X)$ in such a way that $d$ is additive and multiplicative! 

\item In pivotal categories this can be done. A strict {\bf pivotal} category \cite{FY1,FY2} is a
strict left rigid category with a monoidal structure on the functor $X\mapsto{}^\vee\! X$ and a
monoidal equivalence of the functors $\id_\2C$ and $X\mapsto {}^{\vee\vee}X$. As a consequence, one
can define a right duality satisfying $X^\vee={}^\vee\! X$.

\item In a strict pivotal categories we can define left and right traces for every endomorphism:
\be\label{eq-traces} \Tr_X^L(s)=\quad\quad
\begin{tangle}
\coev\step[.3]\obj{e_X}\\
\step[-1]\obj{\ol{X}}\step\id\Step\O{s}\\
\ev\step[-.4]\obj{d'_X}
\end{tangle}
\quad\quad\quad\quad\quad
\Tr_X^R(s)=\quad\quad
\begin{tangle}
\coev\step[.3]\mobj{e'_X}\\
\O{s}\Step\id\step[.3]\obj{\ol{X}}\\
\ev\step[-.4]\mobj{d_X}
\end{tangle}
\ee

Notice: In general $\Tr^L_X(s)\ne\Tr^R_X(s)$.

\item We now define dimensions by $d(X)=\Tr^L_X(\id_X)\in\End\,\11$. One then automatically has 
  $d(\ol{X})=\Tr_X^R(\id_X)$, which can differ from $d(X)$. But for simple $X$ we have
  $d(X)d(\ol{X})=d^2(X)$ with $d^2(X)$ as above. 

\item In a pivotal category, we can use the trace to define pairings
  $\Hom(X,Y)\times\Hom(Y,X)\rarr\End\,\11$ by $(s,t)\mapsto\Tr^L_X(t\circ s)$. In the semisimple
  $k$-linear case with $\End\,\11$, these pairings are non-degenerate for all
  $X,Y$. Cf.\ e.g.\ \cite{GK1}. In general, a morphism $s:X\rarr Y$ is called {\bf negligible} if
  $\Tr(t\circ s)=0$ for all $t:Y\rarr X$. We call an Ab-category {\it non-degenerate} if only
the zero morphisms are negligible. The negligible morphisms form a monoidal ideal,
  i.e.\ composing or tensoring a negligible morphism with any morphism yields a negligible
  morphism. It follows that one can quotient out the negligible morphisms in a straightforward way,
  obtaining a non-degenerate category. A non-degenerate abelian category is semisimple \cite{del4},
but a counterexample given there shows that non-degeneracy plus pseudo-abelianness do not imply
semisimplicity! 

\item A {\bf spherical} category \cite{BW1} is a pivotal category where the left and right
traces 
coincide. Equivalently, it is a strict autonomous category (i.e.\ a tensor category equipped with a
left and a right duality) for which the resulting functors $X\mapsto X^\vee$ and $X\mapsto{}^\vee\! X$
coincide. 

Sphericity implies $d(X)=d(\ol{X})$, and if $\2C$ is semisimple, the converse implication holds.
\item The Temperley-Lieb categories $\2T\2L(\tau)$ are spherical.
\item A finite dimensional Hopf algebra that is involutive, i.e. satisfies $S^2=\id$, gives rise to
  a spherical category. (It is known that every semisimple and co-semisimple Hopf algebra is
  involutive.) More generally, `spherical Hopf algebras', defined as satisfying $S^2(x)=wxw^{-1}$,
  where $w\in H$ is invertible with $\Delta(w)=w\otimes w$ and $\Tr(\theta w)=\Tr(\theta w^{-1})$ for
  any finitely generated projective left $H$-module $V$, give rise to spherical categories \cite{BW1}.

\item In a $*$-category with conjugates, traces of endomorphisms, in particular dimensions of
objects, can be defined uniquely without choosing a spherical structure, cf.\ \cite{DR,lro}. The
dimension satisfies $d(X)\ge 1$ for every non-zero $X$, and $d(X)=1$ holds if and only if $X$ is
invertible. Furthermore, one has \cite{lro} a $*$-categorical version of the quantization of the
Jones index \cite{vfr0}: 
\[ d(X)\in \left\{ 2\cos\frac{\pi}{n},\ n=3,4,\ldots\right\} \cup[2,\infty). \]

On the other hand, every tensor $*$-category can be equipped \cite{yama4} with an (essentially)
unique spherical structure such the traces and dimension defined using the latter coincide with
those of \cite{lro}. 

\item In a $\7C$-linear fusion category (no $*$-operation required!) one has $d^2(X)>0$ for all
  $X$, cf.\ \cite{eno}. 

The following is a very useful application: If $\2A\subset\2B$ is a full inclusion of $\7C$-linear
fusion category then $\dim\2A\le\dim\2B$, and equality holds if and only if $\2A\simeq\2B$.

\item In a unitary category, $\dim\2C=FP-\dim\2C$. Categories with the latter property are called
{\bf pseudo-unitary} in \cite{eno}, where it is shown that every pseudo-unitary category admits a
unique spherical structure such that $FP-d(X)=d(X)$ for all $X$.

\item There are Tannaka-style theorem for not necessarily symmetric categories (Ulbrich \cite{ulb2},
  Yetter \cite{yet1}, Schauenburg \cite{schau1}): 
Let $\2C$ be a $k$-linear pivotal category with $\End\,\11=k\id_\11$ and let $E:\2C\rarr\Vect_k$
a fiber functor. Then the algebra $A(E)$ defined as above admits a coproduct and an antipode, thus
the structure of a Hopf algebra $H$, and an equivalence $F: \2C\rarr\mathrm{Comod}\,H$ such that
$E=K\circ F$, where $K:\mathrm{Comod}\,H\rarr\Vect_k$ is the forgetful functor.
(If $\2C$ and $E$ are symmetric, this $H$ is a commutative Hopf algebra of functions on the group 
obtained earlier $G$.)
Woronowicz proved a similar result \cite{woro2} for $*$-categories, obtaining a compact quantum
group (as defined by him \cite{woro1,woro3}). Commutative compact quantum groups are just algebras
$C(G)$ for a compact group, thus one recovers Tannaka's theorem. Cf.\ \cite{js3} for an excellent
introduction to the area of Tannaka-Krein reconstruction.

\item Given a fiber functor, can one find an algebraic structure whose {\it representations} (rather
than corepresentations) are equivalent to $\2C$? The answer is positive, provided one uses a slight
generalization of Hopf algebras, to wit A.\ van Daele's `Algebraic Quantum Groups' \cite{VD2,VD} (or
`Multiplier Hopf algebras with Haar functional'). They are not necessarily unital algebras equipped
with a coproduct $\Delta$  that takes values in the multiplier algebra $M(A\otimes A)$ and with a 
left-invariant Haar-functional $\mu\in A^*$. A nice feature of algebraic quantum groups is that they
admit a nice version of Pontryagin duality (which is not the case for infinite dimensional ordinary
Hopf algebras). 

In \cite{MRT} the following was shown: If $\2C$ is a semisimple spherical ($*$-)category and $E$ a
($*$-)fiber functor then there is a discrete   multiplier Hopf ($*$-)algebra $(A,\Delta)$ and an
equivalence $F:\2C\rarr\Rep(A,\Delta)$ such that   $K\circ F=E$, where $K:\Rep(A,\Delta)\rarr\Vect$
is the forgetful functor. (This $(A,\Delta)$ is 
  the Pontrjagin dual of the $A(E)$ above.) This theory exploits the semisimplicity from the very
  beginning, which makes it quite transparent: One defines
\[ A=\bigoplus_{i\in I}\End\,E(X_i)
\quad\quad\mbox{and}\quad\quad M(A)=\prod_{i\in I}\End\,E(X_i)\,\cong\,\mathrm{Nat}\,E,\]
where the summation is over the equivalence classes of simple objects in $\2C$.
Now the tensor structures of $\2C$ and $E$ give rise to a coproduct $\Delta:A\rarr M(A\otimes A)$ in
a very direct way.

Notice: This reconstruction is related to the preceding one as follows. Since
$H-\mathrm{comod}\simeq\2C$ is semisimple, the Hopf algebra $H$ has a left-invariant integral $\mu$,
thus $(H,\mu)$ is a compact algebraic quantum group, and the discrete algebraic quantum group
$(A,\Delta)$ is just the Pontrjagin dual of the latter.

\item In this situation, there is a bijection between braidings on $\2C$ and R-matrices (in
$M(A\otimes A)$), cf.\ \cite{MRT}. But: The braiding on $\2C$ plays no essential r\^ole in the
  reconstruction. (Since \cite{MRT} works with the category of finite dimensional representations,
  which in general does not contain the left regular representation, this is more work than e.g.\ in
  \cite{kassel} and requires the use of semisimplicity.) 
\item Summing up: Linear [braided] tensor categories admitting a fiber functor are
  (co)representa\-tion categories of [(co)quasi-triangular] discrete (compact) quantum groups. 

Notice that here `Quantum groups' refers to Hopf algebras and suitable generalizations thereof, but
not necessarily to q-deformations of some structure arising from groups!

\item WARNING: The non-uniqueness of fiber functors means that there can be non-isomorphic quantum
  groups whose (co)representation categories are equivalent to the given $\2C$!

The study of this phenomenon leads to Hopf-Galois theory and is connected (in the $*$-case) to the
study of ergodic actions of quantum groups on  $C^*$-algebras. (Cf.\ e.g.\ Bichon, de Rijdt, Vaes
\cite{BV}).

\item Despite this non-uniqueness, one may ask whether one can {\it intrinsically}
  characterize the tensor categories admitting a fiber functor, thus being related to quantum
  groups. (Existence of a fiber functor is an extrinsic criterion.) The few known results to this
  questions are of two types. On the one hand there are some recognition theorems for certain classes of
 representation categories of quantized enveloping algebras, which will be discussed somewhat
 later. On the other hand, there are results based on the regular representation, to which we turn
 now. However, it is only in the symmetric case that this leads to really satisfactory results.

\item The left regular representation $\pi_l$ of a compact group $G$ (living on $L^2(G)$) has the
following  well known properties:
\bean \pi_l\cong\bigoplus_{\pi\in\widehat{G}} d(\pi)\cdot\pi, && \mbox{(Peter-Weyl theorem)} \\ 
  \pi_l\otimes\pi\cong d(\pi)\cdot\pi_l\quad\forall\pi\in\Rep\,G. &&\mbox{(absorbing property)}.\eean

\item The second property generalizes to {\it any} algebraic quantum group' $(A,\Delta)$, cf.\
\cite{MT}: 

1. Let $\Gamma=\pi_l$ be the left regular representation. If $(A,\Delta)$ is discrete, then $\Gamma$
carries a monoid structure $(\Gamma,m,\eta)$ with $\dim\Hom(\11,\Gamma)=\11$, which we call the {\bf
regular monoid}. (Algebras in $k$-linear tensor categories satisfying $\dim\Hom(\11,\Gamma)=\11$
have been called `simple' or `haploid'.) If $(A,\Delta)$ is compact, $\Gamma$ has a comonoid
structure. (And in the finite (=compact + discrete) case, the algebra and coalgebra structures
combine to a Frobenius algebra, cf.\ \cite{mue09}, discussed below.)  

2. If $(A,\Delta)$ is a discrete algebraic quantum group, one has a monoid version of the absorbing
property: For every $X\in\Rep(A,\Delta)$ one has an isomorphism 
\be (\Gamma\otimes X,m\otimes\id_X)\cong n(X)\cdot (\Gamma,m) \label{1}\ee
of $(\Gamma,m,\eta)$-modules in $\Rep(A,\Delta)$. (Here $n(X)\in\7N$ is the dimension of the vector
space of the representation $X$, which in general differs from the categorical dimension.)

\item The following theorem from \cite{MT} is motivated by Deligne's \cite{del}:  Let $\2C$ be a
$k$-linear  category and $(\Gamma,m,\eta)$ a monoid in $\2C$ (more generally, in the associated
category Ind$\,\2C$ of inductive limits) satisfying $\dim\Hom(\11,\Gamma)=1$ and (\ref{1}) for some
function $n:\Obj\,\2C\rarr\7N$. Then 
\[ E(X)=\Hom_{\Vect_k}(\11,\Gamma\otimes X) \]
defines a faithful $\otimes$-functor  $E:\2C\rarr\Vect_k$, i.e.\ a fiber functor. (One has $\dim
E(X)=n(X)\ \forall X$ and $\Gamma\cong\oplus_i n(X_i)X_i$.) If $\2C$ is symmetric and
$(\Gamma,m,\eta)$ commutative (i.e.\ $m\circ c_{\Gamma,\Gamma}=m$), then $E$ is symmetric.

Remark: Deligne considered this only in the symmetric case, but did not make the requirement
$\dim\Hom(\11,\Gamma)=1$. This leads to a tensor functor $E:\2C\rarr A-\Mod$, where
$A=\Hom(\11,\Gamma)$ is the commutative $k$-algebra of `elements of $\Gamma$' encountered earlier.

\item This gives rise to the following implications:
\[\begin{array}{c}\begin{picture}(300,150)(80,0)\thicklines

\put(160,130){There is a discrete AQG $(A,\Delta)$}

\put(170,115){such that $\2C\simeq\Rep(A,\Delta)$}


\put(290,20){$\2C$ admits an absorbing monoid}


\put(20,20){There is a fiber functor}
\put(50,05){$E:\2C\rarr\2H$}


\put(130,50){\vector(1,1){45}}

\put(290,95){\vector(1,-1){45}}

\put(265,23){\vector(-1,0){40}}

\end{picture} \end{array}\]

Remarks: 1. This can be considered as an intrinsic characterization of quantum group categories. (Or
rather semi-intrinsic, since the regular monoid lives in the Ind-category of $\2C$ rather than $\2C$ 
itself.) 

2. The case of finite $*$-categories had been treated in \cite{lo3}, using subfactor theory
and a functional analysis.

3. This result is quite unsatisfactory, but I doubt that a better result can be obtained without
restriction to {\it special classes of categories} or adopting a wide {\it generalization of the
  notion of quantum groups}. Examples for both will be given below.

4. For a different approach, also in terms of the regular representation, cf.\ \cite{dopr}.

\item Notice that having an absorbing monoid in $\2C$ (or rather Ind$(\2C)$) means having an
  $\7N$-valued dimension function $n$ on the hypergroup $I(\2C)$ and an associative product on
  the object $\Gamma=\oplus_{i\in I}n_iX_i$. The latter is a cohomological condition. 

If  $\2C$ is finite, one can show using Perron-Frobenius theory that there is only one dimension
function, namely the intrinsic one $i\mapsto d(X_i)$. Thus a finite category with non-integer
intrinsic dimensions cannot be Tannakian (in the above sense). 

\item We now turn to a very beautiful result of Deligne \cite{del} (simplified considerably by
  Bichon \cite{bichon}): 

Let $\2C$ be a semisimple $k$-linear rigid {\bf even symmetric} category satisfying $\End\,\11=k$,
where $k$ is algebraically closed of characteristic zero. Then there is an absorbing commutative
monoid as above. (Thus we have a symmetric fiber functor, implying $\2C\simeq\Rep\,G$.) 

Sketch: 
The homomorphisms $\Pi_n^X:S_n\rarr\Aut\,X^{\otimes  n}$ allow to define the idempotents
   \[ P_\pm(X,n)=\frac{1}{n!}\sum_{\sigma\in S_n}\mathrm{sgn}(\sigma)\,\Pi_n^X(\sigma) 
   \ \in \End(X^{\otimes  n})\]
and their images $S^n(X),A^n(X)$, which are direct summands of $X^{\otimes n}$. Making crucial use
of the evenness assumption on $\2C$, one proves
\[ d(A^n(X))=\frac{d(X)(d(X)-1)\cdots(d(X)-n+1)}{n!}\quad\forall n\in\7N. \]
In a $*$-category, this must be non-negative $\forall n$, implying $d(X)\in\7N$, cf.\ \cite{DR}.
Using this -- or assuming it as in \cite{del} -- one has $d(A^{d(X)}(X))=1$, and $A^{d(X)}(X)$ is called the
{\bf determinant} of $X$. On the other hand, one can define a commutative monoid structure on
\[ S(X)=\bigoplus_{n=0}^\infty S^n(X), \]
obtaining the {\bf symmetric algebra} $(S(X),m,\eta)$ of $X$. Let $Z$ be a $\otimes$-generator $Z$
of $\2C$ satisfying $\det\,Z=\11$. Then the `interaction' between symmetrization (symmetric algebra)
and antisymmetrization (determinants) allows to construct a maximal ideal $I$ in the commutative
algebra $S(Z)$ such that the quotient algebra $A=S(Z)/I$ has all desired properties: it is
commutative, absorbing and satisfies $\dim\Hom(\11,A)=1$. QED.

Remarks: 1. The absorbing monoid $A$ constructed in \cite{del,bichon} did {\it not} satisfy
$\dim\Hom(\11,A)=1$. Therefore the construction considered above does not give a fiber functor
to $\Vect_\7C$, but to $\Gamma_A-\Mod$, and one needs to quotient by a maximal ideal in
$\Gamma_A$. Showing that one can achieve $\dim\Hom(\11,A)=1$ was perhaps the main innovation of
\cite{mue16}. This has the advantage that $(A,m,\eta)$ actually is (isomorphic to) the regular
monoid of the group $G=\mathrm{Nat}_\otimes E$. As a consequence, the latter group can be obtained
simply as the automorphism group 
\[ \mathrm{Aut}(\Gamma,m,\eta)\equiv\{ g\in\Aut\,\Gamma\ | \ g\circ m=m\circ g\otimes g,\ g\circ\eta=\eta\}\]
of the monoid -- without even mentioning fiber functors!

2. Combining Tannaka's theorem with those on fiber functors from monoids and with the above, one has
the following beautiful

Theorem \cite{DR,del}: Let $k$ be algebraically closed of characteristic zero and $\2C$ a semisimple
k-linear rigid even symmetric category with $\End\,\11=k$. Assume that all objects have dimension in
$\7N$. Then there is a pro-algebraic group $G_a$, unique up to isomorphism, such that
$\2C\simeq\Rep\,G_a$ (finite dimensional rational representations). If $\2C$ is a $*$-category then
semisimplicity and the dimension condition are redundant, and there is a unique compact group $G_c$
such that $\2C\simeq\Rep\,G_c$ (continuous unitary finite-dimensional representations). In this case,
$G_a$ is the complexification of $G_c$.

3. If $\2C$ is symmetric but not even, its symmetry can be `bosonized' into an even one, cf.\
\cite{DR}. Then one applies the above result and obtains a group $G$. The $\7Z_2$-grading on $\2C$
given by the twist gives rise to an element $k\in Z(G)$ satisfying $k^2=e$. Thus
$\2C\simeq\Rep(G,k)$ as symmetric category. Cf.\ also \cite{del6}.

\item The above result has several applications in pure mathematics: It plays a big r\^ole in the
  theory of motives \cite{andre,levine} and in differential Galois theory and the related Riemann Hilbert 
problem, cf.\ \cite{put}. It is used for the classification of triangular Hopf algebras in terms of
Drinfeld twists of group algebras (Etingof/Gelaki, cf.\ \cite{gelaki} and references therein) and
for the modularization of braided tensor categories \cite{brug3,mue06}, cf.\ below.

The work of Doplicher and Roberts \cite{DR} was motivated by applications to quantum field theory in
$\ge 2+1$ dimensions \cite{dhr3,DR2}, where it leads to a Galois theory of quantum fields, cf.\ also
\cite{halv}. 

\item Thus, at least in characteristic zero (in the absence of a $*$-operation one needs to impose
integrality of all dimensions) rigid symmetric categories with $\End\,\11=k\id_\11$ are reasonably
well understood in terms of compact or pro-affine groups.  What about relaxing the last  
  condition? The category of a representations (on continuous fields of Hilbert spaces) of a compact
  {\bf groupoid} $\2G$ is a symmetric $C^*$-tensor   category. Since a lot of information is lost in
  passing from $\2G$ to $\Rep\,\2G$, there is no hope of reconstructing $\2G$ up to isomorphism, but
  one may hope to find a compact group bundle giving rise to the given category and proving that it
  is Morita equivalent to $\2G$. However, there seem to be  topological   obstructions to this being
  always the case, cf.\ \cite{vasselli2}.  

\item In this context, we mention related work by Brugui\`eres/Maltsiniotis
  \cite{mal1,brugM,brug1} on Tannaka theory for {\bf quasi quantum groupoids} in a purely algebraic
  setting. 

\item We now turn to the characterization of certain {\it special classes} of tensor categories:

\item Combining Doplicher-Roberts reconstruction with the mentioned result of McMullen and Handelman
  one obtains a simple prototype: If $\2C$ is an even symmetric tensor $*$-category with conjugates
 and $\End\,\11=\7C$ whose fusion hypergroup is isomorphic to that of a connected compact Lie group $G$,
  then $\2C\simeq\Rep\,G$.  

\item Kazhdan/Wenzl \cite{Wenzl2}: Let $\2C$ be a semisimple $\7C$-linear spherical
$\otimes$-category with $\End\,\11=\7C$, whose fusion hypergroup is isomorphic to that of
$\6s\6l(N)$. Then there is a $q\in\7C^*$ such that $\2C$ is equivalent (as a tensor category)  to
the representation category of the Drinfeld/Jimbo   quantum group $SL_q(N)$ (or one of finitely
many twisted versions of it). Here $q$ is either $1$   or not a root of unity and unique up to
$q\rarr q^{-1}$. (For another approach to a characterization of the $SL_q(N)$-categories, excluding
the root of unity case, cf.\ \cite{pinz}.)

Furthermore: If $\2C$ is a semisimple $\7C$-linear rigid $\otimes$-category with $\End\,\11=\7C$,
whose fusion hypergroup is isomorphic to that of the (finite!) representation category of $SL_q(N)$,
where $q$ is a primitive root of unity of order $\ell>N$, then $\2C$ is equivalent to
$\Rep\,SL_q(N)$ (or one of finitely many twisted versions).

We will say (a bit) more on quantum groups later. The reason that we mention the Kazhdan/Wenzl
result already here is that it does not require $\2C$ to come with a braiding.
Unfortunately, the proof is not independent of quantum group theory, nor does it provide a {\it
construction} of the categories.

Beginning of proof: The assumption on the fusion rules implies that $\2C$ has a multiplicative generator
$Z$. Consider the full monoidal subcategory $\2C_0$ with objects $\{ Z^{\otimes n},\ n\in\7Z_+\}$. 
Now $\2C$ is equivalent to the idempotent completion (`Karoubification') of $\2C_0$. (Aside: Tensor
categories with objects $\7N_+$ and $\otimes=+$ for objects appear quite often: The symmetric category
$\7S$, the braid category $\7B$, PROPs \cite{macl3}.) A semisimple $k$-linear category with objects
$\7Z_+$ is called a {\bf monoidal algebra}, and is equivalent to having a family $\2A=\{A_{n,m}\}$
of vector spaces together with semisimple algebra structures on $A_n=A_{n,n}$ and bilinear
operations $\circ: A_{n,m}\times A_{m,p}\rarr A_{n,p}$ and $\otimes: A_{n,m}\times A_{p,q}\rarr A_{n+p,m+q}$
satisfying obvious axioms. A monoidal algebra is  {\bf diagonal} if $A_{n,m}=0$ for $n\ne m$ and
of {\bf type N} if $\dim A(0,n)=\dim A(n,0)=1$ and $A_{n,m}=0$ unless $n\equiv m (\mod\, N)$. If $\2A$
is of type $N$, there are exactly $N$ monoidal algebras with the same diagonal. The possible
diagonals arising from type $N$ monoidal algebras can be classified, using Hecke algebras $H_n(q)$
(defined later). 

\item There is an analogous result (Tuba/Wenzl \cite{tuW}) for categories with the other classical
(BCD) fusion rings, but that {\it does} require the categories to come with a braiding.

\item For fusion categories, there are a number of classification results in the case of low rank
(number of simple objects) (Ostrik: fusion categories of rank 2 \cite{ostrik3}, braided fusion
  categories of rank 3 \cite{ostrik4}) or special dimensions, like $p$ or $pq$ (Etingof/Gelaki/Ostrik
\cite{EGO}). Furthermore, one can classify {\bf near group categories}, i.e.\ fusion categories with
all simple objects but one invertible (Tambara/Yamagami \cite{TY}, Siehler \cite{siehler}).

\item In another direction one may try to represent more tensor categories as module categories by
  {\it generalizing the notion of Hopf algebras}. We have already encountered a very modest (but
  useful) generalization, to wit Van Daele's multiplier Hopf algebras. (But the main rationale
  for the latter was to repair the breakdown of Pontrjagin duality for infinite dimensional Hopf
  algebras, which works so nicely for finite dimensional Hopf algebras.)

\item Drinfeld's {\bf quasi-Hopf algebras} \cite{drin3} go in a different direction: One considers
  an associative unital algebra $H$ with a unital algebra homomorphism $\Delta:H\rarr H\otimes H$, where
  coassociativity holds only up to conjugation with an invertible element $\phi\in H\otimes H\otimes H$:
\[ \id\otimes\Delta\mcirc\Delta(x)=\phi(\Delta\otimes\id\mcirc\Delta(x))\phi^{-1}, \]
where $(\Delta,\phi)$ must satisfy some identity in order for $\Rep\,H$ with the tensor product
defined in terms of $\Delta$ to be (non-strict) monoidal. Unfortunately, duals of quasi-Hopf
algebras are not quasi-Hopf algebras. They are useful nevertheless, even for the proof of results
concerning ordinary Hopf algebras, like the Kohno-Drinfeld theorem for $U_q(\6g)$,
cf.\ \cite{drin3,drin4} and \cite{kassel}.

Examples: Given a finite group $G$ and $\omega\in Z^3(G,k^*)$, there is a finite dimensional 
quasi Hopf algebra $D^\omega(G)$, the twisted quantum double of Dijkgraaf/Pasquier/Roche
\cite{dpr}. (We will later define its representation category in a purely categorical way.) Recently,
Naidu/Nikshych \cite{NN1} have given necessary and sufficient conditions on pairs 
$(G,[\omega]), (G',[\omega]')$ for $D^\omega(G)-\Mod, D^{\omega'}(G')-\Mod$ to be equivalent as
braided tensor categories. But the question for which pairs $(G,[\omega])$ $D^\omega(G)-\Mod$ is
Tannakian (i.e.\ admits a fiber functor and therefore is equivalent to the representation category
of an ordinary Hopf algebra) seems to be still open.

\item There have been various attempts at proving generalized Tannaka reconstruction theorems in
  terms of quasi-Hopf algebras \cite{majid2} and ``weak quasi-Hopf
  algebras''. (Cf.\ e.g.\ \cite{MS,HO}.) As it turned out, it is sufficient to consider `weak', but 
  `non-quasi' Hopf algebras:   

\item Preceded by Hayashi's `face algebras' \cite{haya1}, which largely went unnoticed, B\"ohm and
  Szlach\'anyi \cite{szl0} and then Nikshych, Vainerman, L. Kadison introduced {\bf weak Hopf
    algebras}, which may be considered as finite-dimensional quantum groupoids: They are associative
  unital algebras $A$ with {\it coassociative} algebra homomorphism  $\Delta:A\rarr A\otimes A$, but
  the axioms $\Delta(\11)=\11\otimes\11$ and $\ve(\11)=1$ are weakened. 
 
Weak Hopf algebras are closely related to Hopf algebroids and have various desirable properties:
Their duals are weak Hopf algebras, and Pontrjagin duality holds. The categorical dimensions of
their representations can be non-integer. And they are general enough to `explain'
finite-index depth-two inclusions of von Neumann factors, cf.\ \cite{NV}.

\item Furthermore, Ostrik \cite{ostrik} proved that every fusion category is the module category of
  a semisimple weak Hopf algebra. (Again, there was related earlier work by Hayashi \cite{haya2} in
  the context of his face algebras  \cite{haya1}.)  

Proof idea: An {\bf $R$-fiber functor} on a fusion category $\2C$ is a faithful tensor functor 
$\2C\rarr\mathrm{Bimod}\,R$, where $R$ is a finite direct sum of matrix algebras. Szlach\'anyi
\cite{szl1}: An R-fiber functor on $\2C$ gives rise to an equivalence $\2C\simeq A-\Mod$ for a weak
Hopf algebra (with base $R$). (Cf.\ also \cite{phh}.) How to construct an R-fiber functor?

Since $\2C$ is semisimple, we can choose an algebra $R$ such that $\2C\simeq R-\Mod$ (as abelian
categories). Since $\2C$ is a module category over itself, we have a $\2C$-module structure on
$R-\Mod$. Now use that, for $\2C$ and $R$ as above, there is a bijection between R-fiber functors and
$\2C$-module category structures on $R-\Mod$ (i.e.\ tensor functors $\2C\rarr\End(R-\Mod)$.

Remarks: 1. $R$ is highly non-unique: The only requirement was that the number of simple direct summands
equals the number of simple objects of $\2C$. (Thus there is a unique commutative such $R$, but even
for that, there is no uniqueness of $R$-fiber functors.)

2. The above proof uses semisimplicity. (A non-semisimple generalization was announced by
Brugui\`eres and Virezilier in 2008.)

\item Let $\2C$ be fusion category and $A$ a weak Hopf algebra such that $\2C\simeq A-\Mod$. Since
  there is a dual weak Hopf algebra $\widehat{A}$, it is natural to ask how
  $\widehat{\2C}=\widehat{A}-\Mod$ is related to $\2C$. (One may call such a category dual to $\2C$,
  but must keep in mind that there is one for every weak Hopf algebra $A$ such that $\2C\simeq A-\Mod$.)

\item Answer: $\widehat{A}-\Mod$ is (weakly monoidally) Morita equivalent to $\2C$. This notion
  (M\"uger \cite{mue09}) was inspired by subfactor theory, in particular ideas of Ocneanu, cf.\
  \cite{ocn1,ocn2}. For this we need the following: 

\item A {\bf Frobenius algebra} in a strict tensor category is a quintuple $(A,m,\eta,\Delta,\ve)$,
where $(A,m,\eta)$ is an algebra, $(A,\Delta,\ve)$ is a coalgebra and the Frobenius identity
\[ m\otimes\id_A\mcirc\id_A\otimes\Delta=\Delta\mcirc m=\id_A\otimes m\mcirc\Delta\otimes\id_A\]
holds. Diagrammatically:
\[ \begin{tangle} \hh\hcd\step\id\\ \hh\id\step\hcu\end{tangle}
\quad=\quad\begin{tangle}\hh\hcu\\ \hh\hcd\end{tangle}
\quad=\quad\begin{tangle}\hh\id\step\hcd \\ \hh\hcu\step\id\end{tangle}.\]
A Frobenius algebra in a $k$-linear category is called {\bf strongly separable} if 
\[ \ve\circ\eta=\alpha\,\id_\11, \quad\quad m\circ\Delta=\beta\,\id_\Gamma, \quad\quad\alpha\beta\in
k^*.\]
The roots of this definition go quite far back. F.\ Quinn \cite{quinn} discussed them under the name
`ambialgebras', and L.\ Abrams \cite{abrams} proved that Frobenius algebras in
$\Vect_k^{\mathrm{fin}}$ are the usual Frobenius algebras, i.e.\ $k$-algebras 
$V$ equipped with a $\phi\in V^*$ such that $(x,y)\mapsto\phi(xy)$ is non-degenerate. Frobenius
algebras play a central r\^ole for topological quantum field theories in $1+1$ dimensions,
cf.\ e.e.\ \cite{kock}. 

\item Frobenius algebras arise from two-sided duals in tensor categories: Let $X\in\2C$ with 
  two-sided dual $\ol{X}$, and define $\Gamma=X\otimes\ol{X}$. Then $\Gamma$ carries a Frobenius algebra
  structure, cf.\ \cite{mue09}:
\[ m=\quad
\begin{tangle}
\object{X}\step[4]\object{\ol{X}}\\
\id\step\coev\obj{e_X}\step\id\\
\object{X}\step\object{\ol{X}}\Step\object{X}\step\object{\ol{X}}
\end{tangle}
\quad\quad\quad
\Delta=\quad
\begin{tangle}
\object{X}\step\object{\ol{X}}\Step\object{X}\step\object{\ol{X}}\\
\id\step\ev\step[-.5]\obj{d'_X}\step[1.5]\id\\
\object{X}\step[4]\object{\ol{X}}
\end{tangle}
\quad\quad\quad
\eta=\quad
\begin{tangle}
\object{X}\Step\object{\ol{X}}\\
\ev\step[-.5]\obj{d_X}
\end{tangle}
\quad\quad\quad
\ve=\quad
\begin{tangle}
\coev\obj{e'_X}\\
\object{X}\Step\object{\ol{X}}
\end{tangle}
\]
Verifying the Frobenius identities and strong separability is a trivial exercise. In view of 
$\End(V)\cong V\otimes V*$ in the category of finite dimensional vector spaces, the above Frobenius
algebra is called an `endomorphism (Frobenius) algebra'.

\item This leads to the question whether every (strongly separable) Frobenius algebra in a
  $\otimes$-category arise in this way. The answer is, not quite, but: If $\Gamma$ is a strongly
  separable Frobenius algebra in a $k$-linear spherical tensor  category $\2A$ then there exist 
\begin{itemize}
\item a spherical $k$-linear 2-category $\2E$ with two objects $\{\6A,\6B\}$,
\item a 1-morphism $X\in\Hom_\2E(\6B,\6A)$ with 2-sided dual $\ol{X}\in\Hom_\2E(\6A,\6B)$, and
  therefore a Frobenius algebra $X\circ\ol{X}$ in the $\otimes$-category $\End_\2E(\6A)$,
\item a monoidal equivalence $\End_\2E(\6A)\stackrel{\simeq}{\rarr}\2A$ mapping the the Frobenius
  algebra  $X\circ\ol{X}$ to $\Gamma$. 
\end{itemize}
Thus every Frobenius algebra in $\2A$ arises from a 1-morphism in a bicategory $\2E$ containing
$\2A$ as a corner. In this situation, the tensor category $\2B=\End_\2E(\6B)$ is called weakly
monoidally Morita equivalent to $\2A$ and the bicategory $\2E$ is called a Morita context. 

\item The original proof in \cite{mue09} was tedious. Assuming mild technical conditions on $\2A$
and strong separability of $\Gamma$, the bicategory  $\2E$  can simply be obtained as follows:
\bean \Hom_\2E(\6A,\6A) &=& \2A, \\
\Hom_\2E(\6A,\6B) &=& \Gamma-\Mod_\2A,\\
\Hom_\2E(\6B,\6A) &=& \Mod_\2A-\Gamma, \\
\Hom_\2E(\6B,\6B) &=& \Gamma-\Mod_\2A-\Gamma, \eean
with the composition of 1-morphisms given by the usual tensor products of (left and right)
$\Gamma$-modules. Cf.\ \cite{yama3}. (A discussion free of any technical assumptions on $\2A$ was
recently given in \cite{lauda}.)

\item Weak monoidal Morita equivalence of tensor categories also admits an interpretation in terms
  of module categories: If $\6A,\6B$ are objects in a bicategory $\2E$ as above, the category
  $\Hom_\2E(\6A,\6B)$ is a left module category over the tensor category $\End_\2E(\6B)$ and a right
  module category over $\2A=\End_\2E(\6A)$. In fact, the whole structure can be formulated in terms
  of module categories, thereby getting rid of the Frobenius algebras, cf.\  \cite{EO1,eno}: 
Writing $\2M=\Hom_\2E(\6A,\6B)$, the dual category $\2B=\End_\2E(\6B)$ can be obtained as the tensor
category HOM$_\2A(\2M,\2M)$, denoted $\2A^*_\2M$ in \cite{EO1}, of right $\2A$-module functors from
$\2M$ to itself. 

Since the two pictures are essentially equivalent, the choice is a matter of taste. The picture with
Frobenius algebras and the bicategory $\2E$ is closer to subfactor theory. What speaks in favor of
the module category picture is the fact that non-isomorphic algebras in $\2A$ can have equivalent
module categories, thus give rise to the same $\2A$-module category. (But not in the case of
commutative algebras!)

\item Morita equivalence of tensor categories indeed is an equivalence relation, denoted $\approx$. (In
  particular, $\2B$  contains a strongly   separable Frobenius algebra $\widehat{\Gamma}$ such that
  $\widehat{\Gamma}-\Mod_\2B-\widehat{\Gamma}\simeq\2A$.)

\item As mentioned earlier, the left regular representation of a finite dimensional Hopf algebra $H$
  gives rise to a Frobenius algebra $\Gamma$ in $H-\Mod$. $\Gamma$ is strongly separable if and only if $H$ is
  semisimple and cosemisimple. In this case, one finds for the ensuing Morita equivalent category:
\[ \2B=\Gamma-\Mod_{H-\Mod}-\Gamma\ \simeq\ \widehat{H}-\Mod. \]
(This is a situation encountered earlier in subfactor theory.) Actually, in this case the Morita
context $\2E$ had been defined independently by Tambara \cite{tambara}.

The same works for weak Hopf algebras, thus for any semisimple and co-semisimple weak Hopf algebra we have
$A-\Mod\approx\widehat{A}-\Mod$, provided the weak Hopf algebra is Frobenius, i.e.\ has a
non-degenerate integral. (It is unknown whether every weak Hopf algebra is Frobenius.)

\item The above concept of Morita equivalence has important applications: If $\2C_1,\2C_2$ are
  Morita equivalent (spherical) fusion  categories then
\begin{enumerate}
\item $\dim\2C_1=\dim\2C_2$.
\item $\2C_1$ and $\2C_2$ give rise to the same triangulation TQFT in 2+1 dimensions (as defined by
  Barrett/Westbury \cite{BW2} and S.\ Gelfand/Kazhdan \cite{GK1}, generalizing the Turaev/Viro TQFT
\cite{TV,turaev} to non-braided categories. Cf.\ also Ocneanu \cite{ocn3}.)

This fits nicely with the known fact (Kuperberg \cite{kup1}, Barrett/Westbury \cite{BW3}) that, the
spherical categories $H-\Mod$ and $\widehat{H}-\Mod$ (for a semisimple and co-semisimple Hopf algebra $H$) give
rise to the same triangulation TQFT.

\item The braided centers $Z_1(\2C_1), Z_1(\2C_2)$ (to be discussed in the next section) are
  equivalent as braided tensor categories. This is quite immediate by a result of Schauenburg
  \cite{schau2}.
\end{enumerate}

\item We emphasize that (just like Vect$_k$) a fusion category can contain many (strongly separable)
  Frobenius algebras, thus it can be Morita equivalent to many other tensor categories. In view of
  this, studying (Frobenius) algebras in fusion categories is an important and interesting
  subject. (Even more so in the braided case.) 

\item Example: Commutative algebras in a representation category $\Rep\,G$ (for $G$ finite) are the
  same as commutative algebras carrying a   $G$-action   by algebra automorphisms. The condition
  $\dim\Hom(\11,\Gamma)=1$ means that the 
  $G$-action is   ergodic. Such algebras correspond to closed subgroups $H\subset G$ via
  $\Gamma_H=C(G/H)$. Cf.\ \cite{ko}.

\item Algebras in and module categories over the category $\2C_k(G,\omega)$ defined in Section
  \ref{sec-1} were studied in \cite{ostrik2}. 

\item A {\bf group theoretical category} is a fusion category that is weakly Morita equivalent
(or `dual') to a pointed fusion category, i.e.\ one of the form $\2C_k(G,\omega)$ (with $G$ finite
  and $[\omega]\in H^3(G,\7T)$). (The original definition \cite{ostrik} was in terms of quadruples
$(G,H,\omega,\psi)$ with $H\subset G$ finite groups, $\omega\in Z^3(G,\7C^*)$ and 
$\psi\in C^2(H,\7C^*)$ such that $d\psi=\omega_{|H}$, but the two notions are equivalent by Ostrik's
analysis of module categories of $\2C_k(G,\omega)$ \cite{ostrik}.) For more on group theoretical
categories cf.\ \cite{naidu1,GN2}.

\item The above considerations are closely related to subfactor theory (at finite Jones index): A
  {\bf factor} is a von Neumann algebra with center $\7C\11$. For 
  an inclusion $N\subset M$ of factors, there is a notion of index $[M:N]\in[1,+\infty]$ (not
  necessarily integer!!), cf.\ \cite{vfr0,lo1}. One has $[M:N]<\infty$ if and only if the canonical
N-M-bimodule $X$ has a dual 1-morphism $\ol{X}$ in the bicategory of von Neumann algebras, bimodules
and their intertwiners. Motivated by Ocneanu's bimodule picture of subfactors \cite{ocn1,ocn2} one
observes that the bicategory with the objects $\{N,M\}$ and bimodules generated by $X,\ol{X}$ is a
Morita context. On the other hand, a single factor $M$ gives rise to a certain tensor $*$-category
$\2C$ (consisting of $M-M$-bimodules or the endomorphisms $\End\,M$) such that, by Longo's work
\cite{lo3}, the   Frobenius algebras (``Q-systems'' \cite{lo3}) in $\2C$ are (roughly) in bijection
with the subfactors $N\subset M$ with $[M:N]<\infty$. (Cf.\ also the introduction of \cite{mue09}.)


\end{itemize}


\section{Braided tensor categories}\label{sec-4}
\begin{itemize}
\item The symmetric groups have the well known presentation
\[ S_n=\{ \sigma_1,\ldots, \sigma_{n-1}\ | \ \sigma_i\sigma_{i+1}\sigma_i=\sigma_{i+1}\sigma_i\sigma_{i+1},\
    \sigma_i\sigma_j=\sigma_j\sigma_i \ \ \mbox{when}\ \ |i-j|>1,\ \sigma_i^2=1\}. \]
Dropping the last relation, one obtains  the {\bf Braid groups}:
\[ B_n=\{ \sigma_1,\ldots, \sigma_{n-1}\ | \ \sigma_i\sigma_{i+1}\sigma_i=\sigma_{i+1}\sigma_i\sigma_{i+1},\
    \sigma_i\sigma_j=\sigma_j\sigma_i \ \ \mbox{when}\ \ |i-j|>1\}. \]
They were introduced by Artin in 1928, but had appeared implicitly in much earlier work by Hurwitz,
cf.\ \cite{kt3}. They have a natural geometric interpretation:

{\unitlength=0.6mm
\begin{picture}(200,42)(-10,0)
\put(5,16){\makebox(0,0){\mbox{$\sigma_1=$}}}
\put(20,0){\makebox(0,0){\mbox{$\bullet$}}}
\put(30,0){\makebox(0,0){\mbox{$\bullet$}}}
\put(40,0){\makebox(0,0){\mbox{$\bullet$}}}
\put(52,0){\makebox(0,0){\mbox{$\bullet$}}}
\put(20,30){\makebox(0,0){\mbox{$\bullet$}}}
\put(30,30){\makebox(0,0){\mbox{$\bullet$}}}
\put(40,30){\makebox(0,0){\mbox{$\bullet$}}}
\put(52,30){\makebox(0,0){\mbox{$\bullet$}}}
\put(47,15){\makebox(0,0){\mbox{$\cdots$}}}
\put(20,0){\line(1,3){10}} \put(30,0){\line(-1,3){4}}
\put(20,30){\line(1,-3){4}} \put(40,0){\line(0,1){30}}
\put(52,0){\line(0,1){30}} \put(54,16){\makebox(0,0){\mbox{,}}}
\put(70,16){\makebox(0,0){\mbox{$\sigma_2=$}}}
\put(80,0){\makebox(0,0){\mbox{$\bullet$}}}
\put(90,0){\makebox(0,0){\mbox{$\bullet$}}}
\put(100,0){\makebox(0,0){\mbox{$\bullet$}}}
\put(110,0){\makebox(0,0){\mbox{$\bullet$}}}
\put(80,30){\makebox(0,0){\mbox{$\bullet$}}}
\put(90,30){\makebox(0,0){\mbox{$\bullet$}}}
\put(100,30){\makebox(0,0){\mbox{$\bullet$}}}
\put(110,30){\makebox(0,0){\mbox{$\bullet$}}}
\put(105,15){\makebox(0,0){\mbox{$\cdots$}}}
\put(80,0){\line(0,1){30}} \put(90,0){\line(1,3){10}}
\put(100,0){\line(-1,3){4}} \put(90,30){\line(1,-3){4}}
\put(110,0){\line(0,1){30}} \put(112,16){\makebox(0,0){\mbox{,}}}
\put(129,16){\makebox(0,0){\mbox{$\sigma_{n-1}=$}}}
\put(143,0){\makebox(0,0){\mbox{$\bullet$}}}
\put(155,0){\makebox(0,0){\mbox{$\bullet$}}}
\put(165,0){\makebox(0,0){\mbox{$\bullet$}}}
\put(175,0){\makebox(0,0){\mbox{$\bullet$}}}
\put(143,30){\makebox(0,0){\mbox{$\bullet$}}}
\put(150,15){\makebox(0,0){\mbox{$\cdots$}}}
\put(155,30){\makebox(0,0){\mbox{$\bullet$}}}
\put(165,30){\makebox(0,0){\mbox{$\bullet$}}}
\put(175,30){\makebox(0,0){\mbox{$\bullet$}}}
\put(143,0){\line(0,1){30}}
\put(165,0){\line(1,3){10}}
\put(175,0){\line(-1,3){4}}
\put(165,30){\line(1,-3){4}}
\put(155,0){\line(0,1){30}}
\end{picture}
}

Note: $B_n$ is infinite for all $n\ge 2$, $B_2\cong\7Z$. The representation theory of $B_n, n\ge 3$ is
difficult. It is known that all $B_n$ are linear, i.e.\ they have faithful finite dimensional
representations $B_n\hookrightarrow GL(m,\7C)$ for suitable $m=m(n)$. Cf. Kassel/Turaev \cite{kt3}.

\item Analogously, one can drop the condition $c_{Y,X}\circ c_{X,Y}=\id$ on a symmetric tensor
  category. This leads to the concept of a {\bf braiding}, due to Joyal and Street \cite{js0,js1},
  i.e.\ a family of natural isomorphisms $c_{X,Y}: X\otimes Y\rarr Y\otimes X$ satisfying two hexagon
identities but not necessarily the condition $c^2=\id$. Notice that without the latter condition,
one needs to require two hexagon identities, the second being obtained from the first one by the
replacement $c_{X,Y}\leadsto c_{Y,X}^{-1}$ (which does nothing when $c^2=\id$). (The latter is the
non-strict generalization of $c_{X\otimes Y,Z}=c_{X,Z}\otimes\id_Y\circ\id_X\otimes c_{Y,Z}$.)
A {\bf braided tensor category} (BTC) now is a tensor category equipped with a braiding.

\item In analogy to the symmetric case, given a BTC $\2C$ and $X\in\7N,\ n\in\7Z_+$, one has a
homomorphism $\Pi_n^X:B_n\rarr\mathrm{Aut}(X^{\otimes n})$.

\item The most obvious example of a BTC that is not symmetric is provided by the braid
  category $\7B$. In analogy to the symmetric category $\7S$, it is defined by
  $\mathrm{Obj}\,\7B=\7Z_+$, $\End(n)=B_n$, $n\otimes m=n+m$,   while on the morphisms $\otimes $ is
  defined by juxtaposition of braid diagrams. The definition of the braiding
  $c_{n,m}\in\End(n+m)=B_{n+m}$ is illustrated by the example $(n,m)=(3,2)$: 
\[ c_{n,m}=  \quad
\begin{picture}(100,50)(0,30)
\put(0,0){\line(1,1){70}}
\put(30,0){\line(1,1){70}}
\put(60,0){\line(1,1){70}}

\put(85,0){\line(-1,1){10}}
\put(70,15){\line(-1,1){10}}
\put(55,30){\line(-1,1){10}}
\put(40,45){\line(-1,1){25}}

\put(115,00){\line(-1,1){20}}
\put(85,30){\line(-1,1){10}}
\put(70,45){\line(-1,1){10}}
\put(55,60){\line(-1,1){10}}
\end{picture}
\]

\vspace{1cm}

\item If $\2C$ is a strict BTC and $X\in\2C$, there is a unique braided tensor functor
$F:\7B\rarr\2C$ such that $F(1)=X$ and $F(c_{2,2})=c_{X,X}$. Thus $\7B$ is the {\bf free braided
tensor category} generated by one object.  

\item Centralizer and center $Z_2$:

If $\2C$ is a BTC, we say that two objects $X,Y$ {\bf commute} if $c_{Y,X}\circ c_{X,Y}=\id_{X\otimes Y}$.
If $\2D\subset\2C$ is subcategory (or just subset of $\Obj\,\2C)$, we define the {\bf centralizer}
$\2C\cap\2D'\subset\2C$ as the full subcategory defined by  
\[ \Obj(\2C\cap\2D')=\{ X\in\2C\ | \ c_{Y,X}\circ c_{X,Y}=\id_{X\otimes Y}\ \ \forall Y\in\2D\}. \]
Now, the {\bf center} $Z_2(\2C)$ is
\[ Z_2(\2C)=\2C\cap\2C'. \]

Notice that $\2C\cap\2D'$ is monoidal and $Z_2(\2C)$ is symmetric! In fact, a BTC $\2C$ is symmetric
if and only if   $\2C=Z_2(\2C)$. Apart from `central', the objects of $Z_2(\2C)$ have been called `degenerate'
\cite{khr} or `transparent' \cite{brug3}.
\item We thus see that STC are maximally commutative BTCs. Does it make sense to speak of maximally
  non-commutative BTCs?  $\7B$ is an example since $\obj\,Z_2(\7B)=\{0\}$. Braided fusion categories
  with `trivial' center will turn out to be just Turaev's modular categories, cf.\ Section \ref{sec-5}.

\item Since the definition of BTCs is quite natural if one knows the braid groups, one may wonder
  why they appeared more than 20 years after symmetric categories. Most likely, this was a
  consequence of a lack of really interesting examples. When they finally appeared in \cite{js0},
  this was mainly motivated by developments internal to category theory (and homotopy theory). It is
  a remarkable historical accident that this happened at the same time as (and independently from)
  the development of quantum groups, which dramatically gained in popularity in the wake of
  Drinfeld's talk  \cite{drin1}. 
\item In 1971 it was shown \cite{dhr3} that certain representation theoretic considerations for
  quantum field theories in spacetimes of dimension $\ge 2+1$ lead to symmetric categories. Adapting
  this theory to $1+1$ dimensions inevitably leads to braided categories, as was finally shown in
  1989, cf.\ \cite{frs1}. That this was not done right after the appearance of \cite{dhr3} must be 
  considered as a missed opportunity. 

\item As promised, we will briefly look at braided categorical groups. Consider $\2C(G)$ for $G$
  abelian. As shown in \cite{js1} -- and in much more detail in the preprints \cite{js0} -- the
  braided categorical groups $\2C$ with $\pi_0(\2C)\cong G$ (isomorphism classes of objects) and
$\pi_1(\2C)\cong A$ ($\End\,\11$) are classified by the group $H^3_{ab}(G,A)$,  where
$H^n_{ab}(G,A)$ refers to the Eilenberg-Mac Lane cohomology theory for {\it abelian} groups, 
  cf. \cite{macl1}. (Whereas $H^3(G,A)$ can be defined in terms of topological cohomology theory as
  $H^3(K(G,1),A)$ of the Eilenberg-Mac Lane space $K(G,1)$, one has
$H^3_{ab}(G,A):=H^4(K(G,2),A)$. This group also has a description in terms of quadratic functions 
  $q:G\rarr A$. The subgroup of $H^3_{ab}(G,A)$ corresponding to symmetric braidings is isomorphic
to $H^5(K(G,3),A)$, cf.\ \cite{CK}.) 

\item Duality: Contrary to the symmetric case, in the presence of a (non-symmetric) braiding, having a
  left duality is not sufficient for a nice theory: If we define a right duality in terms of a left
  duality and the braiding, the left and right traces will fail to have all the properties they do
  have in the symmetric case. Therefore, some additional concepts are needed:

\item A {\bf twist} for a braided category with left duality is a natural family
$\{ \Theta_X  \in\End\,X,\ X\in\2C\}$ 
of isomorphisms (i.e.\ a natural isomorphism of the functor $\id_\2C$) satisfying 
\[ \Theta_{X\otimes Y}=\Theta_X\otimes\Theta_Y\mcirc c_{Y,X}\circ c_{X,Y}, \quad\quad\quad\Theta_\11=\id_\11,
 \quad\quad\quad{}^\vee\!(\Theta_X)=\Theta_{{}^\vee\! X}. \]
Notice: If $c_{Y,X}\circ c_{X,Y}\not\equiv\id$ then the natural isomorphism $\Theta$ is {\it not
  monoidal} and $\Theta=\id$ is not a legal twist!

\item A {\bf ribbon category} is a strict braided tensor category equipped with a left duality and a
twist. 

\item Let $\2C$ be a ribbon category with left duality $X\mapsto({}^\vee\!\!X,e_X,d_X)$ and twist
  $\Theta$. We   define a right duality $X\mapsto (X^\vee, e'_X,d'_X)$ by $X^\vee={}^\vee\!\! X$ and
(\ref{eq-e'}). Now one can show, cf.\ e.g.\ \cite{kassel}, that the maps $\End\,X\rarr\End\,\11$
  defined as in (\ref{eq-traces}) coincide and that $\Tr(s):=\Tr_L(s)=\Tr_R(s)$ has the trace property
  and behaves well under tensor products, as previously in the symmetric case. Writing
  $\ol{X}={}^\vee\! X=X^\vee$, one finds that $\2C$ is a spherical category in the sense of
  \cite{BW1}. Conversely, if $\2C$ is spherical and braided, then {\it defining} 
\[ \Theta_X=(\Tr_X\otimes\id_X)(c_{X,X}), \]
$\{ \Theta_X,\ X\in\2C\}$ satisfies the axioms of a twist and thus forms a ribbon structure together
with the left duality. (Cf.\ Yetter \cite{y}, based on ideas of Deligne, and Barrett/Westbury \cite{BW1}.) 

(Personally, I prefer to consider the twist as a derived structure, thus talking about spherical
categories with a braiding, rather than about ribbon categories. In some situations, e.g.\ when the center
$Z_1(\2C)$ is involved, this is advantageous. This also is the approach of the Rome school \cite{DR2,lro}.)

\item So far, our only example of a non-symmetric braided category is the free braided category
  $\7B$, which is not rigid. In the remainder of this section, we will consider three main `routes'
  to braided categories: (A) the topological route, (B) the ``non-perturbative approach'' via
  quantum doubles and categorical centers, and (C) the ``perturbative approach'' via deformation (or
  `quantization') of symmetric categories.

\item We briefly mention one construction of an interesting braided category that doesn't seem to
  fit nicely into one of our routes: While the usual representation category of a group is symmetric,
  the category of representations of the general linear group $GL_n(\7F_q)$ over a finite field with
  the  {\it external} tensor product of representations turns out to be braided and non-symmetric, cf.\
\cite{js5}. 
\end{itemize}

\subsection{Route A: Free braided categories (tangles) and their quotients}
\begin{itemize}

\item Combining the ideas behind the Temperley-Lieb categories TL$(\tau)$ (which have duals) and the
  braid category $\7B$ (which is braided but has no duals), one arrives at the categories of {\bf
    tangles} (Turaev \cite{turaev0}, Yetter \cite{yet0}. See also \cite{turaev,kassel}.) One must
distinguish between categories of {\bf unoriented tangles} having  $\Obj\,U-\2T\2A\2N=\7Z_+$ with tensor
  product (of objects) given by addition and {\bf oriented tangles}, based on
  $\Obj\,O-\2T\2A\2N=\{+,-\}^*$ (i.e.\ finite words in $\pm,\ \11=\emptyset$) with concatenation as
  tensor product. In either case, the morphisms are given as sets of pictures as in Figure 1, or
  else by linear combinations of such pictures with coefficients in a commutative ring or
  field. All this is just as in the discussion of the free symmetric categories at the end of
  Section \ref{sec-2}. The only difference is that one must distinguish between over- and
  undercrossings of the lines; for technical reasons it is more convenient to do this in terms of
  pictures embedded in 3-space. 

\begin{figure}[h]
         \centerline{
           \scalebox{1}{
             \epsfig{file=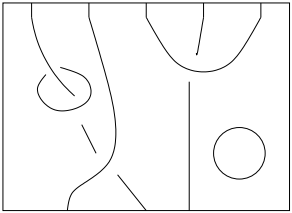, scale=0.6}
             }
           }
         \caption{An unoriented 3-5 tangle}
         \label{fig1}
       \end{figure}

There also is a category $O-\2T\2A\2N$ of oriented tangles, where the objects are finite words in
$\pm,\ \11=\emptyset$ and the lines in the morphisms are directed, in a way that is compatible with
the signs of the objects. It is clear that the morphisms in $\End(\11)$ in $U-\2T\2A\2N$
($O-\2T\2A\2N$) are just the unoriented (oriented) links. 

While the definition is intuitively natural, the details are tedious and we refer to the textbooks
\cite{turaev,kassel,yetter}. In particular, we omit discussing {\it  ribbon} tangles. 

\item The tangle categories are pivotal, in fact spherical, thus ribbon categories. $O-\2T\2A\2N$ is the
free ribbon category generated by one element, cf.\ \cite{shum}.

\item Let $\2C$ be a ribbon category. Then one can define a category $\2C-\2T\2A\2N$ of $\2C$-labeled
oriented tangles and a ribbon tensor functor $F_\2C: \2C-\2T\2A\2N\rarr\2C$. (This is the rigorous
rationale behind the diagrammatic calculus for braided tensor categories!)

Let $\2C$ be a ribbon category and $X$ a self-dual object. Given an unoriented tangle, we can label
every edge by $X$. This gives a composite map
\[ \{ \mbox{links} \} \stackrel{\cong}{\longrightarrow} \Hom_{U-\2T\2A\2N}(0,0) \longrightarrow
   \Hom_{\2C-\2T\2A\2N}(0,0) \stackrel{F_\2C}{\longrightarrow} \End_\2C\11. \]
In particular, if $\2C$ is $k$-linear with $\End\,\11=k\id$, we obtain a map from \{ links \} to $k$,
which is easily seen to be a knot invariant. If $\2C=U_q(\6s\6l(2))-\Mod$ and $X$ is the fundamental
object, one essentially obtains the Jones polynomial. Cf.\ \cite{turaev0,rt1}.
(The other objects of $\2C$ give rise to the colored Jones polynomials, which are much studied in
the context of the volume conjecture for hyperbolic knots.)

\item So far, all our examples of braided categories have come from topology. In a sense, they are quite
  trivial, since they are just the universal braided (ribbon) categories freely generated by one
  object. Furthermore, we are primarily interested in linear categories. Of course, we can apply the
  $k$-linearization functor $\2C\2A\2T\rarr k$-$\mathrm{lin.}$-$\2C\2A\2T$. But the categories we
  obtain have infinite dimensional hom-sets and are not more interesting than the original
  ones. (This should be contrasted to the symmetric case, where this construction produces the
  representation categories of the classical groups, cf.\ Section \ref{sec-2}.)

\item Thus in order to obtain interesting $k$-linear ribbon categories from the tangle categories,
  we must reduce the infinite dimensional hom-spaces to finite dimensional ones.

We consider the following analogous situation in the context of associative algebras: The braid group
$B_n\ (n>1)$ is infinite, thus the group algebra $\7CB_n$ is infinite 
dimensional. But this algebra has finite dimensional quotients, e.g.\ the {\bf  Hecke algebra}
$H_n(q)$, the unital $\7C$-algebra generated by $\sigma_1,\ldots,\sigma_{n-1}$, modulo the relations 
\[ \sigma_i\sigma_{i+1}\sigma_i=\sigma_{i+1}\sigma_i\sigma_{i+1},\quad\quad
    \sigma_i\sigma_j=\sigma_j\sigma_i \ \ \mbox{when}\ \ |i-j|>1, \quad\quad
   \sigma_i^2=(q-1)\sigma_i+q\11. \]
This algebra is finite dimensional for any $q$, and for $q=1$ we have $H_n(q)\cong\7CS_n$. In fact,
$H_n(q)$ is isomorphic to $\7CS_n$, thus semisimple, whenever $q$ is not a root of unity, but this
isomorphism is highly non-trivial. Cf.\ e.g.\ \cite{lehrer}. 

The idea now is to do a similar thing on the level of categories, or to `categorify' the Hecke
algebras or other quotients of $\7CB_n$ like the Birman/Murakami/Wenzl- (BMW-)-algebras \cite{BW}. 

\item We have seen that ribbon categories give rise to knot invariants. One can go the other way
  and construct $k$-linear ribbon categories from link invariants. This approach was initiated in
  \cite[Chapter XII]{turaev}, where a topological construction of the representation category of
  $U_q(sl(2))$ was given. A more general approach was studied in \cite{tw2}. A $k$-valued link
  invariant $G$ is said to {\bf admit functorial extension to tangles} if there exists a tensor
  functor $F: U-\2T\2A\2N\rarr k-\Mod$ whose restriction to
  $\End_{U-\2T\2A\2N}(0)\cong\{\mbox{links}\}$ equals $G$. 

For any $X\in U-\2T\2A\2N,\ f\in\End(X)$, let $L_f$ be the link obtained by closing $f$ on the
right, and define $\Tr_G(f)=G(L_f)$. If $\2C$ is the $k$-linearization of $U-\2T\2A\2N$, it is shown
in \cite{tw2}, under weak assumptions on $G$, that the idempotent and direct sum completion of the
quotient of $\2C$ by the ideal of negligible morphisms is a semisimple ribbon category with finite
dimensional hom-sets. Cf.\ \cite{tw2}. 

Example: Applying the above procedure $G=V_t$, the Jones polynomial, one obtains a Temperley-Lieb
category $\2T\2L_\tau$, which in turn is equivalent to a category $U_q(\6s\6l(2))-\Mod$. 
Cf.\ \cite[Chapter XII]{turaev}. Applying it to the Kauffman polynomial \cite{kauff}, one obtains
the quantized BTCs of types BCD, cf.\ \cite{tw2}. The general theory in \cite{tw2} is quite nice,
but it should be noted that the assumption of functorial extendability to tangles is rather strong:
It implies that the resulting semisimple category admits a fiber functor and therefore is the
representation category of a discrete quantum group. Furthermore, the application of the general
formalism of \cite{tw2} to the Kauffman polynomial used input from (q-deformed) quantum group theory
for the proof of functorial extension to tangles and of modularity. This drawback was repaired by
Beliakova/Blanchet, cf.\ \cite{bebl1,bebl2}. 

Blanchet \cite{blanchet} gave a similar construction with HOMFLY polynomial \cite{HOMFLY}, obtaining
the type A categories. (The HOMFLY polynomial is an invariant for oriented links, thus one must work
with oriented tangles.)

Remark: The ribbon categories of BCD type arising from the Kauffman polynomial give rise to
topological quantum field theories. The latter can even be constructed directly from the Kauffman
bracket, bypassing the categories, cf.\ \cite{BHMV}. This construction actually preceded those
mentioned above. 

\item The preceding constructions reinforce the close connection between braided categories and knot
invariants. It is important to realize that this reasoning is not circular, since the polynomials
of Jones, HOMFLY, Kauffman can (nowadays) be constructed in rather elementary ways, independently
of categories and quantum groups, cf.\ e.g.\ \cite{lick}. Since the knot polynomials are
defined in terms of skein relations, we speak of the skein construction of the quantum categories, 
which arguably is the simplest known so far.

\item In the case $q=1$, the skein constructions of the ABCD categories reduce to the construction
  of the categories arising from classical groups mentioned in Section \ref{sec-2}. (This happens
  since $q=1$ corresponds to parameters in the knot polynomials for which they fail to distinguish
  over- from under-crossings. Then one can replace the tangle categories by symmetric categories of
  non-embedded cobordisms (oriented or not) as in \cite{del2}.)

\item Concerning the exceptional Lie algebras and their quantum categories, inspired by work of
Cvitanovic, cf.\ \cite{cvi} for a book-length treatment, and by Vogel \cite{vogel}, Deligne
conjectured \cite{del2} that there is a one parameter family of symmetric tensor categories $\2C_t$
specializing to $\Rep\,G$ for the exceptional Lie groups at certain values of $t$. This is still
unproven, but see \cite{CdM,del3,del5} for work resulting from this conjecture. (For the
$E_n$-categories, including the $q$-deformed ones, cf.\ \cite{wen2}.) 
\item In a similar vein, Deligne defined \cite{del4} a one parameter family of rigid symmetric tensor
categories $\2C_t$ such that $\2C_t\simeq\Rep\,S_t$ for $t\in\7N$. These categories were studied
further in \cite{comes}. (Recall that $S_n$ is considered as the $GL_n(\7F_1)$ where $\7F_1$ is the
`field with one element', cf.\ \cite{soule}.)

\item More generally, one can define linear categories by generators and relations, cf.\ e.g.\ \cite{kup2}.
\end{itemize}

\subsection{Route B: Doubles and centers}
We begin with a brief look at Hopf algebras.
\begin{itemize}
\item {\bf Quasi-triangular} Hopf algebras (Drinfeld, 1986 \cite{drin1}): If $H$ is a Hopf algebra
  and $R$ an invertible element of (possibly a completion of) $H\otimes H$, satisfying
\[ R\Delta(\cdot)R^{-1}=\sigma\circ\Delta(\cdot),\quad\quad\sigma(x\otimes y)=y\otimes x,\]
\[ (\Delta\otimes\id)(R)=   R_{13}R_{23}, \quad (\id\otimes\Delta)(R)=R_{13}R_{12}.\]
\[ (\ve\otimes\id)(R)=(\id\otimes\ve)(R)=\11. \]
If $(V,\pi),(V',\pi')\in H-\Mod$, the definition
$c_{(V,\pi),(V',\pi')}=\Sigma_{V,V'}(\pi\otimes\pi')(R)$ produces a braiding for $H-\Mod$. 

\item But this has only shifted the problem: How to get quasi-triangular Hopf algebras? To this
  purpose,  Drinfeld \cite{drin1} gave the quantum double construction $H\leadsto D(H)$, which
  associates a quasi-triangular Hopf algebra $D(H)$ to a Hopf algebra $H$. Cf.\ also \cite{kassel}. 
\item Soon after, an analogous categorical construction was given by  Drinfeld (unpublished),
  Joyal/Street \cite{js2} and Majid \cite{majid1}): The (braided) center $Z_1(\2C)$, defined as follows. 

Let $\2C$ be a strict tensor category and let $X\in\2C$. A half braiding
$e_X$ for $X$ is a family $\{e_X(Y)\in\Hom_\2C(X\otimes Y,Y\otimes X),\ Y\in\2C\}$ of isomorphisms,
natural w.r.t.\ $Y$, satisfying  $e_X(\11)=\id_X$ and
\[ e_X(Y\otimes Z)=\id_{Y}\otimes e_X(Z) \mcirc e_X(Y)\otimes\id_{Z} \quad \forall Y,Z\in\2C. \]

Now, the {\bf center $Z_1(\2C)$} of $\2C$ has as objects pairs $(X, e_X)$, where $X\in\2C$ and $ e_X$ is a
half braiding for $X$. The morphisms are given by  
\[ \Hom_{Z_1(\2C)}((X, e_X),(Y, e_Y))= \{ t\in\Hom_\2C(X,Y) \ | \  \id_X\otimes t \mcirc e_X(Z)=
   e_Y(Z)\mcirc t\otimes\id_X \quad \forall Z\in\2C \}.\] 
The tensor product of objects is given by 
$(X, e_X)\otimes(Y, e_Y)=(X\otimes Y, e_{X\otimes Y})$, where
\[ e_{X\otimes Y}(Z)= e_X(Z)\otimes\id_Y\mcirc\id_X\otimes e_Y(Z). \]
The tensor unit is $(\11,e_\11)$ where $e_\11(X)=\id_X$. The composition and tensor
product of morphisms are inherited from $\2C$. Finally, the braiding is given by 
\[  c((X, e_X),(Y, e_Y))=  e_X(Y). \]
(The author finds this definition is much more transparent than that of $D(H)$ even though {\it a
  priori} little is known about $Z_1(\2C)$.)

\item Just as the centralizer $\2C\cap\2D'$ generalizes $Z_2(\2C)=\2C\cap\2C'$, there is a version
  of $Z_1$ relative to a subcategory $\2D\subset\2C$, cf.\ \cite{majid1}.

\item $Z_1(\2C)$ is categorical version (generalization) of Hopf algebra quantum double in the
  following sense: If $H$ is a finite dimensional Hopf algebra, there is an equivalence
\be Z_1(H-\Mod)\simeq D(H)-\Mod \label{eq}\ee
of braided tensor categories, cf.\ e.g.\ \cite{kassel}. (If $H$ is infinite dimensional, one still
has an equivalence between $Z_1(H-\Mod)$ and the category of Yetter-Drinfeld modules over $H$.)
\item If $\2C$ is a category and $\2D:=Z_0(\2C)=\End(\2C)$ is its tensor category of endofunctors, then
  $Z_1(\2D)$ is trivial. (This may be considered as the categorification of the fact that the center
  (in the usual sense) of the endomorphism monoid $\End(S)$ of a set $S$ is trivial, i.e.\ equal to
  $\{\id_S\}$.) But in general, the braided center of a tensor category is a non-trivial braided
  category that is not symmetric. Unfortunately, this doesn't seem to have been studied
  thoroughly. Presently, strong results on $Z_1(\2C)$ exist only in the case where $\2C$ is a fusion
  category.

\item There are abstract categorical considerations, quite unrelated to topology and quantum groups,
  that provide rationales for studying BTCs: 

(A): A second, compatible, multiplication functor on a tensor
  category gives rise to a braiding, and conversely, cf. \cite{js1}. (This is a higher
dimensional version of the Eckmann-Hilton argument mentioned earlier.) 

(B): Recall that tensor categories are bicategories with one object. Now, braided tensor categories
turn out to be monoidal bicategories with one object, which in turn are weak 3-categories with one
object and one 1-morphism. Thus braided (and symmetric) categories really are a manifestation of the
existence of $n$-categories for $n>1$!

\item Baez-Dolan \cite{baez3} conjectured the following `periodic table' of `k-tuply monoidal n-categories':

\begin{center}
{\small
\begin{tabular}{|c|c|c|c|c|c|}  \hline
         & $n = 0$   & $n = 1$    & $n = 2$            & $n=3$  & $n=4$ \\     \hline
$k = 0$  & sets      & categories & 2-categories   &  3-categories & \ldots \\     \hline
$k = 1$  & monoids   & monoidal   & monoidal       & monoidal & \ldots\\ 
         &           & categories & 2-categories   & 3-categories & \\     \hline
$k = 2$  &commutative& braided    & braided        & braided & \ldots \\
         & monoids   & monoidal   & monoidal       &  monoidal  & \\
         &           & categories & 2-categories   & 3-categories &  \\     \hline
$k = 3$  &           & symmetric  & `sylleptic' &  &\\
         &     "      & monoidal   & monoidal        & ?  & \ldots \\
         &           & categories & 2-categories    &  &\\    \hline
$k = 4$  &           &          & symmetric & & \\
         &    "      &    "        & monoidal        & ? & \ldots \\
         &           &            & 2-categories    &   &\\     \hline
$k = 5$  &         &          &               & symmetric & \\
         &    "       &   "       &    "             & monoidal & \ldots\\ 
         &           &            &                 & 3-categories & \\     \hline
$k=6$ & " & " & " & "  & \ldots \\ \hline
\end{tabular}} \vskip 1em
\end{center}
In particular, one expects to find `center constructions' from each structure in the table to the
one underneath it. For the column $n=1$ these are the centers $Z_0,Z_1,Z_2$ discussed above. For
$n=0$ they are given by the endomorphism monoid of a set and the ordinary center of a monoid. The
column $n=2$ is also relatively well understood, cf.\ Crans \cite{crans}. There is an accepted
notion of a non-strict 3-category (i.e.\ $n=3, k=0$) (Gordon/Power/Street \cite{GPS}), but there are
many competing definitions of weak higher categories. We refrain from moving any further into this 
subject. See e.g.\ \cite{baezM}.

\item With this heuristic preparation, one can give a high-brow interpretation of $Z_1(\2C)$,
  cf.\ \cite{js1,str2}: Let $\2C$ be tensor category and $\Sigma\2C$ the corresponding bicategory
  with one object. Then the 
  category $\End(\Sigma\2C)$ of endofunctors of $\Sigma\2C$ is a monoidal bicategory (with natural
  transformations as 1-morphisms and `modifications' as 2-morphisms). Now,
  $\2D=\End_{\End(\Sigma\2C)}(\11)$ is a tensor category with two compatible $\otimes$-structures
  (categorifying $\End\,\11$ in a tensor category), thus braided, and it is equivalent to $Z_1(\2C)$. 

\item For further abstract considerations on the center $Z_1$, consider the work of Street
  \cite{str2,str3} and of Brugui\`eres and Virelizier \cite{brug9,brug8}.

\item If $\2C$ is braided there is a braided embedding $\iota_1:\2C\hookrightarrow Z_1(\2C)$, given by
$X\mapsto   (X,e_X)$, where $e_X(Y)=c(X,Y)$. Defining $\widetilde{\2C}$ to be the tensor category
$\2C$ with `opposite' braiding   $\widetilde{c}_{X,Y}=c_{Y,X}^{-1}$, there is an analogous embedding 
$\widetilde{\iota}:\2C\hookrightarrow Z_1(\2C)$. In fact, one finds that the images of $\iota,\iota'$
  are each others' centralizers:
\[ Z_1(\2C)\cap\iota(\2C)'=\widetilde{\iota}(\widetilde{\2C}), \quad\quad
  Z_1(\2C)\cap\widetilde{\iota}(\widetilde{\2C})'=\iota(\2C). \]
Cf.\ \cite{mue10}. On the one hand, this is an instance of the double commutant principle, and on
the other hand, this establishes one connection  
\[ \iota(\2C)\cap\widetilde{\iota}(\widetilde{\2C})=\iota(Z_2(\2C))=\widetilde{\iota}(Z_2(\widetilde{\2C})),\]
between $Z_1$ and $Z_2$ which suggests that
``$Z_1(\2C)\simeq\2C\times\widetilde{\2C}$'' when $Z_2(\2C)$ is ``trivial''.
We will return to both points in the next section.
\end{itemize}

\subsection{Route C: Deformation of groups or symmetric categories}
\begin{itemize}
\item As for Route B, there is a more traditional approach via deformation of Hopf algebras and a
  somewhat more recent one focusing directly on deformation of tensor categories. 

\item (C$_1$): The earlier approach to braided categories relies on deformation of Hopf algebras
  related to groups. For lack of space we will limit ourselves to providing just enough information as
  needed for the discussion of the more categorical approach. For more, we refer to the textbooks, in
  particular \cite{kassel, cp, J, lusztig}. In any case, one chooses a simple (usually compact) Lie
  group $G$ and considers either the enveloping algebra $U(\6g)$ of its Lie algebra $\6g$ in  terms
  of Serre's generators and relations   \cite{serre}, or one departs from the algebra Fun$(G)$ of
  regular functions on $G$, which can also be described in terms of finitely many relations,
  cf.\ e.g.\ \cite{woro1}. In a nutshell, one inserts factors of a `deformation parameter' $q$ into
  the presentation of $U(\6g)$ or Fun$(G)$ in such a way that for $q\ne 1$ one still obtains a
  (non-trivial) Hopf algebra. Quantum group theory began with the discovery that this is possible
  at all.

\item Obviously, this `definition' is a farcical caricature. But there is some truth in it: In
the mathematical literature on quantum groups, cf.\ e.g.\ \cite{kassel, cp, J, lusztig}, it is
all but impossible to find a comment on the origin of the presentation of the quantum group under
study and of the underlying motivation. 
While the initiators of quantum group theory from the Leningrad school (Faddeev, Kulish,
Semenov-Tian-Shansky, Sklyanin, Reshetikhin, Drinfeld and others) were very well aware of these
origins, this knowledge has now almost faded into obscurity. (This certainly has to do with the fact
that the applications to theoretical physics for which quantum groups were invented in the first
place are still exclusively pursued by physicists, cf.\ e.g.\ \cite{gomez}.) One point of this
section will be that -- quite independently of the original physical motivation -- the categorical
approach to quantum deformation is mathematically better motivated.

\item In what  follows, we will concentrate on the  enveloping algebra approach. The usual
  Drinfeld-Jimbo presentation of the quantized enveloping algebra is as follows, Consider the
  algebra $U_q(\6g)$  generated by elements $E_i$, $F_i$, $K_i$, 
$K_i^{-1}$, $1\le i\le r$, satisfying the relations
$$
K_iK_i^{-1}=K_i^{-1}K_i=1,\ \ K_iK_j=K_jK_i,\ \
K_iE_jK_i^{-1}=q_i^{a_{ij}}E_j,\ \
K_iF_jK_i^{-1}=q_i^{-a_{ij}}F_j,
$$
$$
E_iF_j-F_jE_i=\delta_{ij}\frac{K_i-K_i^{-1}}{q_i-q_i^{-1}},
$$
$$
\sum^{1-a_{ij}}_{k=0}(-1)^k\begin{bmatrix}1-a_{ij}\\
k\end{bmatrix}_{q_i} E^k_iE_jE^{1-a_{ij}-k}_i=0,\ \
\sum^{1-a_{ij}}_{k=0}(-1)^k\begin{bmatrix}1-a_{ij}\\
k\end{bmatrix}_{q_i} F^k_iF_jF^{1-a_{ij}-k}_i=0,
$$
where $\displaystyle\begin{bmatrix}m\\
k\end{bmatrix}_{q_i}=\frac{[m]_{q_i}!}{[k]_{q_i}![m-k]_{q_i}!}$,
$[m]_{q_i}!=[m]_{q_i}[m-1]_{q_i}\dots [1]_{q_i}$,
$\displaystyle[n]_{q_i}=\frac{q_i^n-q_i^{-n}}{q_i-q_i^{-1}}$ and
$q_i=q^{d_i}$. This is a Hopf algebra with coproduct $\Delta$ and
counit $\varepsilon$ defined by
$$
\Delta(K_i)=K_i\otimes K_i,\ \
\Delta(E_i)=E_i\otimes1+ K_i\otimes E_i,\ \
\Delta(F_i)=F_i\otimes K_i^{-1}+1\otimes F_i,
$$
$$
\varepsilon(E_i)=\varepsilon(F_i)=0,\ \ \varepsilon(K_i)=1.
$$

One should distinguish between  Drinfeld's \cite{drin1} formal approach, where one constructs a Hopf
algebra $H$ over the ring $\7C[[h]]$ of formal power series in such a way that $H/hH$ is isomorphic
to the enveloping algebra $U(\6g)$, and the non-formal deformation of Jimbo \cite{jimbo}, who
obtains an honest quasi-triangular Hopf algebra $U_q(\6g)$ (over $\7C$) for any value $q\in\7C$ of a
deformation parameter. (In this approach, the properties of the resulting Hopf algebra depend 
  heavily on whether $q$ is a root of unity or not. In the formal approach, this distinction
  obviously does not arise.) The relation between both approaches becomes clear by inserting $q=e^h$
  in Jimbo's definition and considering the result as a Hopf algebra over $\7C[[h]]$.

\item (C$_2$): As mentioned, one can obtain non-symmetric braided categories directly by `deforming'
  symmetric categories. This approach was initiated by Cartier \cite{cartier} and worked out in more
  detail in \cite[Appendix]{kassel} and \cite{kt2}. (These works were all motivated by applications
  to Vassiliev link invariants, which we cannot discuss here.)

Let $\2S$ be a strict symmetric Ab-category. Now an {\bf infinitesimal braiding} on $\2S$ is a
natural family  of endomorphisms $t_{X,Y}: X\otimes Y\rarr X\otimes Y$ satisfying
\[ c_{X,Y}\circ t_{X,Y}=t_{Y,X}\circ c_{X,Y}\quad\forall X,Y, \]
\[ t_{X,Y\otimes Z}=t_{X,Y}\otimes\id_Z+c_{X,Y}^{-1}\otimes\id_Z\mcirc\id_Y\otimes t_{X,Z}\mcirc
   c_{X,Y}\otimes\id_Z\quad \forall X,Y,Z. \] 
Strict symmetric Ab-categories equipped with an infinitesimal braiding were called {\bf
  infinitesimal symmetric}. (We would prefer to call them symmetric categories equipped with an
infinitesimal braiding.)

\item Example: If $H$ is a Hopf algebra, there is a bijection between infinitesimal braidings $t$ on 
$\2S=H-\Mod$ and elements $t\in\mbox{Prim}(H)\otimes\mbox{Prim}(H)$ (where Prim$(H)=\{x\in H\ |
\ \Delta(x)=x\otimes\11+\11\otimes x\}$) satisfying $t_{21}=t$ and
$[t,\Delta(H)]=0$, given by $t_{X,Y}=(\pi_X\otimes\pi_Y)(t)$.

\item Now we can define the formal deformation of a symmetric category associated to an
  infinitesimal braiding: Let $\2S$ be a strict $\7C$-linear symmetric category with finite
  dimensional hom-sets and let $t$ be an infinitesimal braiding for $\2S$. We write $\2S[[h]]$ for
  the $\7C[[h]]$-linear category obtained by extension of
  scalars. (I.e.\ Obj$\,\2S[[h]]=\mbox{Obj}\,\2S$ and  
$\Hom_{\2S[[h]]}(X,Y)=\Hom_\2S(X,Y)\otimes_\7C \7C[[h]]$.) Also the functor
$\otimes:\2S\times\2S\rarr\2S$ lifts to $\2S[[h]]$. For objects $X,Y,Z$, define
\[ \alpha_{X,Y,Z}=\Theta_{KZ}(h\,t_{X,Y}\otimes\id_Z,h\,\id_X\otimes t_{Y,Z}), \quad\quad
  \widetilde{c}_{X,Y}=c_{X,Y}\circ e^{ht_{X,Y}/2}. \]
Here $\Theta_{KZ}$ is a Drinfeld associator \cite{drin3}, i.e.\ a formal power series
\[ \Theta_{KZ}(A,B)=\sum_{w\in\{A,B\}^*} c_w\,w \]
in two non-commuting variables $A,B$, where $c_w\in\7C$, satisfying certain identities. (Cf.\
\cite[Chapter XIX, (8.27)-(8.29)]{kassel}.)
Then $(\2S[[h]],\otimes,\11,\alpha)$ is a (non-strict) tensor category with associativity
constraint $\alpha$, trivial unit constraints and $\widetilde{c}$ a braiding. If $\2S$ is rigid,
then $(\2S[[h]],\otimes,\11,\alpha,\widetilde{c})$ admits a ribbon structure. 

\item Application: Let $\6g$ be a simple Lie algebra$/\7C$. Let $\2S=\6g-\Mod$ and define
  $\{t_{X,Y}\}$ be as 
in the example, corresponding to $t=(\sum_i x_i\otimes x^i + x^i\otimes x_i)/2$, where $x_i,x^i$ are
dual bases of $\6g$ w.r.t.\ the Killing form. Then $[t,\Delta(\cdot)]=0$ and one can prove
\be (\2S[[h]],\otimes,\11,\alpha,\widetilde{c})\ \simeq\ U_h(\6g)-\Mod \label{4}\ee
as $\7C[[h]]$-linear ribbon categories. (The proof is a corollary of the proof of the Kohno-Drinfeld 
theorem \cite{drin3,drin4}, cf.\ also \cite{kassel}.) 

Remark: 1. Obviously, we have cheated: The main difficulty resides in the definition of $\Theta_{KZ}$!
Giving the latter and proving its properties requires ca.\ 10-15 pages of rather technical material
(but no Lie theory). Le and Murakami explicitly wrote down an associator; cf.\ e.g.\ 
\cite[Remark XIX.8.3]{kassel}. Drinfeld also gave a non-constructive proof of existence of an
associator defined over $\7Q$, cf.\ \cite{drin4}.

2. The above is relevant for a more conceptual approach to the theory of finite-type knot
invariants (Vassiliev invariants), cf.\ \cite{cartier,kt2}.

3. A disadvantage of the above is that we obtain only a formal deformation of $\2S$. If $\6g$ is a
simple Lie algebra and $\2S=\6g-\Mod$, we know by (\ref{4}), that we obtain the $\7C[[h]]$-category
$U_h(\6g)-\Mod$. On the other hand, thanks to the work of Jimbo \cite{jimbo} and others
\cite{lusztig,J} we know that there is a non-formal version $U_q(\6g)$ of the quantum group with
$\7C$-linear representation category. One would therefore hope that the $\7C$-linear categories
$U_q(\6g)-\Mod$ can be obtained directly as deformations of the module categories
$U(\6g)-\Mod$. Indeed, for numerical $q\in\7C\backslash\7Q$, with some more analytical effort one
can make sense of $\alpha_q=\Theta_{KZ}(h\,t_{X,Y}\otimes\id_Z,h\,\id_X\otimes t_{Y,Z})$ as an
element of $\End(X\otimes Y\otimes Z)$ and define a non-formal, $\7C$-linear category $\2C(\6g,q)$
and prove an equivalence 
\[ \2C(\6g,q)=(\2S,\otimes,\11,\alpha_q,\widetilde{c}_q)\ \simeq\ U_q(\6g)-\Mod \]
of $\7C$-linear ribbon categories. This was done by Kazhdan and Lusztig \cite{kl}, but see also the
nice recent exposition by Neshveyev/Tuset \cite{NT}.


\item Fact: If $q\in\7C^*$ is generic, i.e.\ {\it not} a root of unity, then
  $\2C(\6g,q):=U_q(\6g)-\Mod$ is a semisimple braided ribbon category whose fusion hypergroup is
  isomorphic to that of $U(\6g)$, thus of the category of $\6g$-modules, cf.\ \cite{J,lusztig}. But
  it is not symmetric for $q\ne 1$, thus certainly not equivalent to the latter. In fact,
  $U_q(\6g)-\Mod$ and  $U(\6g)-\Mod$ are already inequivalent as $\otimes$-categories. (Recall that
  associativity constraints $\alpha$ can be considered as generalized 3-cocycles, and the
  $\alpha_q$ for different $q$ are not cohomologous.)

\item We have briefly discussed the Cartier/Kassel/Turaev formal deformation quantization of
  symmetric categories equip\-ped with an infinitesimal braiding. There is a cohomology theory for
  Ab- tensor categories and tensor functors that classifies deformations due to Davydov \cite{dav0}
  and Yetter \cite{yetter}.

Definition: Let $F:\2C\rarr\2C'$ a tensor functor. Define $T_n:\2C^n\rarr\2C$ by 
$X_1\times\cdots\times X_n\mapsto X_1\otimes\cdots\otimes X_n$. ($T_0(\emptyset)=\11,\ T_1=\id$.)
Let $C^n_F(\2C)=\End(T_n\circ F^{\otimes n})$. ($C^0_F(\2C)=\End\,\11'$.) For a fusion category,
this is finite dimensional. Define $d: C^n_F(\2C)\rarr C^{n+1}_F(\2C)$ by
\[ df= \id\otimes f_{2,\ldots,n+1}-f_{12,\cdots,n+1}+f_{1,23,\ldots,n+1}- \cdots
   +(-1)^nf_{1,\ldots,n(n+1)}+(-1)^{n+1}f_{1,\ldots,n}\otimes\id, \]
where, e.g., $f_{12,3,\ldots,n+1}$ is defined in terms of $f$ using the isomorphism 
$d^F_{X_1,X_2}: F(X_1)\otimes F(X_2)\rarr F(X_1\otimes F_2)$ coming with the tensor functor $F$.

One has $d^2=0$, thus $(C^i,d)$ is a complex. Now $H^i_F(\2C)$ is the cohomology of this complex,
and $H^i(\2C)=H^i_F(\2C)$ for $F=\id_\2C$. 

In low dimensions one finds that  $H^1_F$ classifies derivations of the tensor functor $F$, $H^2_F$
classifies deformations of the tensor structure $\{d^F_{X,Y}\}$ of $F$. $H^3(\2C)$ classifies
deformations of the associativity constraint $\alpha$ of $\2C$.

Examples: 1. If $\2C$ is fusion then $H^i(\2C)=0\ \forall i>0$. This implies Ocneanu rigidity,
cf.\ \cite{eno}. 

2. If $\6g$ is a reductive algebraic group with Lie algebra $\6g$ and $\2C=\Rep\,G$ (algebraic
representations). Then $H^i(\2C)\cong (\Lambda^i\6g)^G \ \forall i$. If $\6g$ is simple then
$H^1(\2C)=H^2(\2C)=0$, but $H^3(\2C)$ is one-dimensional, corresponding to a one-parameter family
of deformations $\2C$. According to \cite{eno} ``it is easy to guess that this deformation comes
from an actual deformation, namely the deformation of $O(G)$ to the quantum group $O_q(G)$''. 
It is not clear to this author whether this suggestion should be considered as proven. If so,
together with the one-dimensionality of $H^3(\6g-\Mod)$ it provides a very satisfactory
`explanation' for the existence of the quantized categories $\2C(\6g,q)\simeq U_q(\6g)-\Mod$. 

\item In analogy to the result of Kazhdan and Wenzl mentioned in Section \ref{sec-3}, Tuba and Wenzl
  \cite{tuW} proved that a semisimple ribbon category with the fusion  hypergroup isomorphic to that of a simple 
classical Lie algebra $\6g$ of BCD type (i.e.\ orthogonal or symplectic) is equivalent to the
category $\2C(\6g,q)$, with $q=1$ or not a root of unity, or one of finitely many twisted versions
thereof. Notice that in contrast to the  Kazhdan/Wenzl result \cite{Wenzl2}, this result needs the
category to be braided! (Again, this is a characterization, not a construction of the categories.)

\item Finkelberg \cite{fink} proved a braided equivalence between $\2C(\6g,q)$, $q=e^{i\pi/m\kappa}$,
  where $m=1$ for ADE, $m=2$ for BCD and $m=3$ for $G_2$, and the ribbon category
  $\tilde{\2O}_\kappa$ of integrable representations of the affine Lie algebra $\hat{\6g}$ of
  central charge $c=\kappa-\check{h}$, where $\check{h}$ is the dual Coxeter number of $\6g$.

The category $\tilde{\2O}_\kappa$ plays an important r\^ole in conformal field theory, either
in terms of vertex operator algebras or via the representation theory of loop groups (Wassermann
\cite{wass}, Toledano-Laredo \cite{Tol1}). This is the main reason for the relevance of quantum
groups to CFT. 

\item Finally, we briefly discuss the connection between routes (B) and (C) to BTCs: In order to find
an R-matrix for the Hopf algebra $U_q(\6g)$ one traditionally uses the quantum double, appealing to 
an isomorphism $U_q(\6g)\cong D(B_q(\6g))/I$, where $B_g(\6g)$ is the q-deformation of a Borel
subalgebra of $\6g$ 
and $I$ an ideal in $D(B_q(\6g))$. Now $R_{U_q(\6g)}=(\phi\otimes\phi)(R_{D(B_q(\6g))})$, where
$\phi$ is the quotient map. Since a surjective Hopf algebra homomorphism $H_1\rarr H_2$ corresponds
to a full monoidal inclusion $H_2-\Mod\hookrightarrow H_1-\Mod$, and recalling the connection (\ref{eq})
between Drinfeld's double construction and the braided center $Z_1$, we conclude that the BTC
$U_q(\6g)-\Mod$ is a full monoidal subcategory of $Z_1(B_q(g)-\Mod)$ (with the inherited
braiding). Therefore, also in the deformation approach, the braiding can be understood as ultimately
arising from the $Z_1$ center construction.  

\item Question: It is natural to ask whether a similar observation also holds for $q$ a root of unity,
i.e., whether the modular categories $\2C(\6g,q)$, for $q$ a root of unity, can be understood as full
$\otimes$-subcategories of $Z_1(\2D)$, where $\2D$ is a fusion category corresponding to the
deformed Borel subalgebra $B_q(\6g)$. Very recently, Etingof and Gelaki \cite{EG2} gave an
affirmative answer in some cases.

Remark: In the next section, we will discuss a criterion that allows to recognize the quantum doubles
$Z_1(\2C)$ of fusion categories.
\end{itemize}


\section{Modular categories}\label{sec-5}
\begin{itemize}
\item Turaev \cite{turaev1,turaev}: A {\bf modular} category is a fusion category that is ribbon
(alternatively, spherical and braided) such that the matrix $S=(S_{i,j})$
\[ S_{i,j}=\Tr_{X\otimes Y}(c_{Y,X}\circ c_{X,Y}),\ \ i,j\in I(\2C), \]
where $I(\2C)$ is the set of simple objects modulo isomorphism, is invertible. 
\item A fusion category that is ribbon is modular if and only if $\dim\2C\ne 0$ and the center
$Z_2(\2C)$ is trivial. (In the sense of consisting only of the objects $\11\oplus\cdots\oplus\11$.)
(This was proven by Rehren \cite{khr} for $*$-categories and by 
  Beliakova/Blanchet \cite{bebl2} in general. Cf.\ also \cite{brug3} and \cite{AB}.) 

Thus: Modular categories are \underline{braided fusion categories with trivial center}, i.e.\
the maximally non-symmetric ones. (This definition seems more conceptual than the original one in
terms of invertibility of $S$.)

\item Why are these categories called `modular'? Let $S$ as above and $T=\mathrm{diag}(\omega_i)$,
where $\Theta_{X_i}=\omega_i\id_{X_i},\ i\in I$. Then 
\[ S^2=\alpha\,C,\quad (ST)^3=\beta\,C, \quad\quad\quad(\alpha\beta\ne 0) \]
where $C_{i,j}=\delta_{i,\ol{\jmath}}$, thus $S,T$ give rise to a projective representation of the
modular group $SL(2,\7Z)$ (which has a presentation $\{ s,t\ | \ (st)^3=s^2=c, \ c^2=e\}$). Cf.\
\cite{khr,turaev}. 

\item At first sight, this is somewhat mysterious. Notice: $SL(2,\7Z)$ is the mapping class group of
the 2-torus $S^1\times S^1$. Now, by work of Reshetikhin/Turaev \cite{rt2,turaev}, providing a
rigorous version of ideas of Witten, every modular category gives rise to a {\bf topological quantum 
field theory} in $2+1$ dimensions. Every such TQFT in turn gives rise to
projective representation of the mapping class groups of all closed 
surfaces, and for the torus one obtains just the above representation of $SL(2,\7Z)$.  Cf.\
\cite{turaev,BK}. We don't have the time to say more about TQFTs.

\item Turaev's motivation came from conformal field theory (CFT). (Cf.\ e.g.\ Moore-Seiberg
\cite{ms2}). In fact, there is a (rigorous) definition of {\bf rational chiral CFT}s (using von
Neumann algebras) and their representations, for which one can prove that the latter are unitary
modular (Kawahigashi, Longo, M\"uger \cite{klm}). Most of the examples considered in the (heuristic)
physics literature fit into this scheme. (E.g.\ the loop group models: \cite{wass,xu} and the minimal
Virasoro models with $c<1$ \cite{loke}.)  

In the context of vertex operator algebras, similar results were proven by Huang \cite{huang}.

\item It is natural to ask whether there are less complicated ways to produce modular categories?
The answer is positive; we will reconsider our three routes to braided categories.

\item Route A: Recall that the classical categories can be obtained from the linearized tangle
categories (type A: oriented tangles, types BCD: unoriented tangles), dividing by ideals defined in
terms of the knot polynomials of HOMFLY and Kauffman. At roots of unity, this leads to modular
categories, cf.\ \cite{tw2,blanchet,bebl2}.

\item Route C$_1$: H.\ Andersen et al. \cite{andersen}, Turaev/Wenzl \cite{tw1} (and others): Let
  ${\6g}$ be a simple Lie algebra and $q$ a primitive root of unity. Then $U_q(\6g)-\Mod$ gives rise
  to a modular category $\2C(\6g,q)$. (Using tilting modules, dividing by negligible morphisms, etc.)


\item Let $q$ be primitive root of unity of order $\ell$. Then $\2C(\6g,q)$ has a positive
$*$-operation (i.e.\ is unitary) if $\ell$ is even (Kirillov Jr.\ \cite{kir0}, Wenzl \cite{wen1}) and is
not unitarizable for odd $\ell$ (Rowell \cite{row1}). 

\item Characterization theorem: A braided fusion category with the fusion hypergroup of
$\2C(\6g,q)$, where $\6g$ is a simple Lie algebra of BCD type and $q$ a root of unity, is equivalent
to $\2C(\6g,q)$ or one of finitely many twisted versions. (Tuba/Wenzl \cite{tuW})

\item Before we reconsider Route B, we assume that we already have a braided fusion category, or
{\bf pre-modular category}.

As we have seen, failure of modularity is due to non-trivial center $Z_2(\2C)$. Idea: Given a
braided (but not symmetric) category with {\it even} center $Z_2(\2C)$, kill the latter, using the 
Deligne / Doplicher-Roberts theorem: 
$Z_2(\2C)\simeq\Rep\,G$. The latter contains a commutative (Frobenius) algebra $\Gamma$
corresponding to the regular representation of $G$. Now $\Gamma-\Mod_\2C$ is modular. (Brugui\`eres
\cite{brug3}, M\"uger \cite{mue06}). This construction can be interpreted as Galois closure in a Galois
theory for BTCs, cf.\ \cite{mue06}.

\item Route B to braided categories: Quantum doubles: If $G$ is a finite group then $D(G)-\Mod$ 
and $D^\omega(G)-\Mod$ are modular (Bantay \cite{bantay1}, Altschuler/ Coste \cite{ac}). If $H$ is a
finite-dimensional semisimple and cosemisimple Hopf algebra then $D(H)-\Mod$ is modular
(Etingof/Gelaki \cite{EG1}). If $A$ is a finite-dimensional weak Hopf algebra then $D(A)-\Mod$
modular (Nikshych/ Turaev/ Vainerman \cite{ntv}).  

\item The center $Z_1$ of a left/right rigid, pivotal, spherical category has the same properties. In
particular, the center of a spherical category is spherical and braided, thus a ribbon
category. (Under weaker assumptions, this is not true, and existence of a twist for the center, if
desired, must be enforced by a categorical version of the ribbonization of a Hopf algebra,
cf.\ \cite{kt1}.) 

\item The braided center $Z_1$: If $\2C$ is spherical fusion category and $\dim\2C\ne 0$ then
$Z_1(\2C)$ is modular and  $\dim Z_1(\2C)=(\dim\2C)^2$. (M\"uger
\cite{mue10}.)  

Comments on the proof: Semisimplicity not difficult. Next, one finds a Frobenius algebra $\Gamma$ in 
$\2D=\2C\boxtimes\2C^\op$ such that the dual category $\Gamma-\Mod_\2D-\Gamma$ is equivalent to
$Z_1(\2C)$, implying $\dim Z_2(\2C)=(\dim\2C)^2$. Here $\Gamma=\oplus_i X_i\boxtimes X_i^\op$, which
is again a coend and can exist also in non-semisimple categories. 

\item This contains all the earlier modularity results on $D(G)-\Mod$ and $D(H)-\Mod$, but also for 
$D^\omega(G)-\Mod$ since:
\[ D^\omega(G)-\Mod\simeq Z_1(\2C_k(G,\omega)). \]
(Using work by Hausser/Nill \cite{HN} or Panaite \cite{pvo} on quantum double of quasi Hopf-algebras.)

\item Modularity of $Z_1(\2C)$ also follows by combination of Ostrik's result that every fusion
category arises from a weak Hopf algebra $A$, combined with modularity of $D(A)-\Mod$ \cite{ntv},
provided one proves $D(A)-\Mod\simeq Z_1(A-\Mod)$, generalizing the known result for Hopf algebras.
But the purely categorical proof avoiding weak Hopf algebras seems preferable, not least since it
probably extends to finite non-semisimple categories.

\item In the Morita context having $\2C\boxtimes\2C^\op$ and $Z_1(\2C)$ as its corners, the two 
off-diagonal categories are equivalent to $\2C$ and $\2C^\op$, and their structures as
$\2C\boxtimes\2C^\op$-module categories are the obvious ones. Therefore, the center can also be
understood as (using the notation of EO):
\[ Z_1(\2C)\simeq(\2C\boxtimes\2C^\op)^*_\2C. \]
A (somewhat sketchy) proof of this equivalence can be found in \cite[Prop.\ 2.5]{ostrik2}.

\item We give another example for a purely categorical result that can be proven using weak Hopf
algebras: Radford's formula for $S^4$ has a generalization to weak Hopf algebras \cite{nik}, and
this can be used to prove that in every fusion category, there exists an isomorphism of tensor
functors $\id\rarr ****$, cf.\ \cite{eno2}. (Notice that in every pivotal category we have $\id\cong
**$, thus here it is important that we understand `fusion' just to mean existence of two-sided
duals. But in \cite{eno} it is conjectured that every fusion category admits a pivotal structure.) 

\item If $\2C$ is already modular then there is a braided equivalence
$Z_1(\2C)\simeq\2C\boxtimes\2C^\op$, cf.\ \cite{mue10}. 
Thus, every modular category $\2M$ is full subcategory of $Z_1(\2C)$ for some fusion category. (This 
probably is not very useful for the classification of modular categories, since there are `more
fusion categories than modular categories': Recall from Section \ref{sec-3} that
$\2C_1\approx\2C_2\ \impl\ Z_1(\2C_1)\simeq Z_1(\2C_2)$. (For converse, see below.) 

\item There is a ``Double commutant theorem'' for modular categories (M\"uger \cite{mue11}, inspired
  by Ocneanu   \cite{ocn4}):
Let $\2M$ a modular category and a $\2C\subset\2M$ a replete full tensor subcategory. Then:
\begin{enumerate}
\item $(\2M\cap(\2M\cap\2C')')=\2C.$
\item $\DS \dim\2C\cdot\dim(\2M\cap\2C')=\dim\2M$,
\item If, in addition $\2C$ is modular, then also $\2D=\2M\cap\2C'$ is modular and
$\2M\simeq\2C\boxtimes\2D$. (Thus every full inclusion of modular categories arises from a direct
product.)
\end{enumerate}
These results indicate that `modular categories are better behaved than finite groups'.

\item Corollary: If $\2M$ is modular and $\2S\subset\2M$ symmetric then
$\2S\subset\2M\cap\2S'$. Thus 
\[ (\dim\2S)^2\le\dim\2S\cdot\dim(\2M\cap\2S')=\dim\2M, \] 
implying $\dim\2S\le\sqrt{\dim\2M}$. Notice that the bound is satisfied by
$\Rep\,G\subset D^{(\omega)}(G)-\Mod$. In fact, existence of a symmetric subcategory attaining the
bound characterizes the representation categories of twisted doubles, cf.\ below.

\item On the other hand, consider $\2C\subset\2M$ with $\2M$ modular. We have 
$\2M\cap\2C'\supset Z_2(\2C)$, implying $\dim\2M\ge\dim\2C\cdot\dim Z_2(\2C)$.
This provides a lower bound on the dimension of a modular category containing a given pre-modular
subcategory as a full tensor subcategory. In \cite{mue11} it was conjectured that this bound can
always be attained.

\item It is natural to ask how primality of $D(G)-\Mod$ is related to simplicity of $G$. It turns
out that the two properties are independent. On the one hand, there are non-simple finite groups for
which $D(G)-\Mod$ is prime. (This is a corollary of the classification of the full fusion subcategories
of $D(G)-\Mod$ given in \cite{NN2}.) On the other hand, for $G=\7Z/p\7Z$ one finds that $D(G)-\Mod$ 
is prime if and only if $p=2$. For $p$ an odd prime, $D(G)-\Mod$ has two prime factors, both of
which are modular categories with $p$ invertible objects, cf.\ \cite{mue11}. But for every finite
simple non-abelian $G$, one finds that $D(G)-\Mod$ is prime. In fact, it has only one replete full
tensor subcategory at all, namely $\Rep\,G$. Thus all these categories are mutually inequivalent:
The classification of prime modular categories contains that of finite simple groups. 

\item If $\2C$ is symmetric and $(\Gamma,m,\eta)$ a commutative algebra in $\2C$, then
$\Gamma-\Mod_\2C$ is again symmetric and
\be \dim\Gamma-\Mod_\2C=\frac{\dim\2C}{d(\Gamma)}. \label{2}\ee
Now, if $\2C$ is only braided, $\Gamma-\Mod_\2C$ is a fusion category satisfying (\ref{2}), but in
general it fails to be braided! (Unless $\Gamma\in Z_2(\2C)$, as was the case in the context of
modularization.) 

\item Example: Given a BTC $\2C\supset\2S\simeq\Rep\,G$, let $\Gamma$ be the regular monoid in
$\2S$ as considered in Section \ref{sec-3}. Then $\2C\rtimes\2S:=\Gamma-\Mod_\2C$ is fusion category, but it
is braided only if $\2S\subset Z_2(\2C)$, as in the discussion of modularization. In general, one
obtains a {\bf braided crossed G-category} as defined by Turaev \cite{turaev6,turaev7} (cf.\ also
Carrasco and Moreno \cite{carr}), i.e.\ a tensor category with $G$-grading $\partial$ on the
objects, a $G$-action $\gamma$ such that $\del(\gamma_g(X))=g\del X g^{-1}$ and a `braiding' 
$c_{X,Y}:X\otimes Y\stackrel{\cong}{\longrightarrow}\gamma_{\del X}(Y)\otimes X$. The
degree zero part is  $\Gamma-\Mod_{\2C\cap\2S'}\simeq\Gamma-\Mod_\2C^0$ (cf.\ below).
(Kirillov Jr.\ \cite{kir1,kir3}, M\"uger \cite{mue13}). This construction has an interesting
connection to conformal orbifold models (\cite{mue15,mue19}).  

\item Even if $\Gamma\not\in Z_2(\2C)$, there is a full tensor subcategory
$\Gamma-\Mod_\2C^0\subset\Gamma-\Mod_\2C$ that is braided. Calling a module $(X,\mu)\in\Gamma-\Mod_\2C$ 
{\bf dyslectic} if  
\[ \mu\circ c_{X,\Gamma}=\mu\circ c_{\Gamma,X}^{-1},  \]
one finds that the full subcategory $\Gamma-\Mod_\2C^0$ of dyslectic modules is not only monoidal,
but also inherits the braiding from $\2C$, cf.\ Pareigis \cite{par3}. This was rediscovered by
Kirillov and Ostrik \cite{ko} who in addition proved that if $\2C$ is modular then
$\Gamma-\Mod_\2C^0$ is modular and the following identity, similar to (\ref{2}) but different,
holds: 
\[ \dim\Gamma-\Mod_\2C^0=\frac{\dim\2C}{d(\Gamma)^2}. \]

Remark: Analogous results were previously obtained by B\"ockenhauer, Evans and Kawa\-higa\-shi
\cite{BEK} in an operator algebraic context. While the transposition of their work to tensor
$*$-categories is immediate, removing the $*$-assumption requires some work.

\item The above implies (for $*$-categories, but also in general over $\7C$ by \cite{eno}) that
$d(\Gamma)\le\sqrt{\dim\2C}$ for commutative Frobenius algebras in modular categories.  (The above
bound on the dimension of full symmetric categories follows from this, since the regular
monoid in $\2S$ is a commutative Frobenius algebra $\Gamma$ with $d(\Gamma)=\dim\2S$.) 

\item All these facts have applications to chiral conformal field theories in the operator algebraic
framework, reviewed in more detail in \cite{mue18}:

Longo/Rehren \cite{lre}: Finite local extensions of a CFT $A$ are classified by the `local
Q-systems' ($\approx$ commutative Frobenius algebras) in $\Rep\,A$, which is a $*$-BTC.

B\"ockenhauer/Evans \cite{BE}, \cite{mue18}: If $B\supset A$ is the finite local extension
corresponding to the commutative Frobenius algebra $\Gamma\in\Rep\,A$, then
$\Rep\,B\simeq\Gamma-\Mod^0_{\Rep\,A}$. 

Analogous results for vertex operator algebras were formulated by Kirillov and Ostrik \cite{ko}.

Remark: It is perhaps not completely absurd to compare these results to local class field theory,
where finite Galois extensions of a local field $k$ are shown to be in bijection to finite index
subgroups of $k^*$.

\item Drinfeld, Gelaki, Nikshych and Ostrik \cite{dgno}, and independently Kitaev and the author,
observed that every commutative Frobenius algebra $\Gamma$ in a modular category $\2M$ gives rise to
a braided equivalence 
\be Z_1(\Gamma-\Mod_\2M)\ \simeq\ \2M\,\boxtimes\, \widetilde{\Gamma-\Mod^0_\2M}. \label{eq5}\ee
Taking $\Gamma=\11$, one recovers the fact $Z_1(\2M)\simeq\2M\boxtimes\widetilde{\2M}$. The latter
raises the question whether one can find a smaller fusion category $\2C$ such that $\2M\subset
Z_1(\2C)$. The answer given by (\ref{eq5}) is that the bigger a commutative algebra one can find in
$\2M$, the smaller one can take $\2C$ to be. In particular, if $\Gamma-\Mod_\2M^0$ is trivial 
(which is equivalent to $d(\Gamma)^2=\dim\2M$ over $\7C$) then $\2M\simeq Z_1(\Gamma-\Mod_\2M)$
is not just contained in a center of a fusion category but is such a center. In fact, this criterion
identifies the modular categories of the form $Z_1(\2C)$ since, conversely, cf.\ \cite{dmno}, one
finds that the center $Z_1(\2C)$ of a fusion category contains a commutative Frobenius algebra 
$\Gamma$ of the maximal dimension $d(\Gamma)=\sqrt{\dim Z_1(\2C)}=\dim\2C$ such that
\[ \Gamma-\Mod_{Z_1(\2C)}^0 \mbox{\ trivial},\quad\quad\quad\Gamma-\Mod_{Z_1(\2C)}\simeq\2C.\]

\item As an application one obtains that if $\2M$ is modular and $\2S\subset\2M$ symmetric and even
such that $\dim\2S=\sqrt{\dim\2M}$ then $\2M\simeq D^\omega(G)-\Mod$, where $\2S\simeq\Rep\,G$ and 
$\omega\in Z^3(G,\7T)$.

This has an application in CFT: If $A$ is a chiral CFT with trivial representation category
$\Rep\,A$ (i.e.\ $A$ is `holomorphic') acted upon by finite group $G$. Then $\Rep\,A^G\simeq
D^\omega(G)-\Mod$. (Together with the results of \cite{klm}, this proves the folk conjecture, having
its roots in \cite{dvvv,dpr}, that the representation category of a `holomorphic chiral orbifold
CFT' is given by a category $D^\omega(G)-\Mod$.)   

\item As shown in \cite{mue09}, a weak monoidal Morita equivalence $\2C_1\approx\2C_1$ of fusion
categories implies $Z_1(\2C_1)\simeq Z_1(\2C_2)$. (This is an immediate corollary of the definition
of $\approx$, combined with \cite{schau2}.)
The converse is true for group theoretical categories (Naidu/Nikshych \cite{NN1}), and a 
general proof is announced by Nikshych.

\item By definition, a group theoretical category $\2C$ is weakly Morita equivalent (dual) to
$\2C_k(G,\omega)$ for a finite group $G$ and $[\omega]\in H^3(G,\7T)$. Thus 
$Z_1(\2C)\simeq Z_1(\2C_k(G,\omega))\simeq D^\omega(G)-\Mod$. The converse is also true.

Therefore, with $\2M$ modular and $\2C$ fusion we have:
\begin{center}
\begin{tabular}{|c|c|c|} \hline 
contains & $\2M$ & $Z_1(\2C)$ \\ \hline
maximal comm.\ FA $\Gamma$  & $\2M\simeq Z_1(\2C)$ & \mbox{always true} \\ \hline
maximal STC $\2S$ & $\2M\simeq D^\omega(G)-\Mod$ & $\2C$ is group theoretical \\ \hline
\end{tabular}
\end{center}

\vspace{1cm}

\item What can we say about non-commutative (Frobenius) algebras in modular categories? We first
  look at the symmetric case. Let thus $\2C$ be a rigid symmetric $k$-linear tensor category and
  $\Gamma$ a strongly separable Frobenius algebra in $\2C$. Define $p\in\End\,\Gamma$ by 
\be\label{eq6} p=(\Tr_\Gamma\otimes\id_\Gamma)(\Delta\circ m\circ c_{\Gamma,\Gamma})=\quad
\begin{tangle}
\Step\object{\Gamma}\\
\hh\hcoev\step\id\\
\hh\id\step\cu\\
\hh\id\step\cd\\
\id\step\hX\\
\hh\hev\step\id\\
\Step\object{\Gamma}
\end{tangle}
\quad=\quad
\begin{tangle}
\step\object{\Gamma}\\
\hh\step\id\\
\dd\id\d\\
\id\step\hX\\
\hh\hev\step\id\\
\Step\object{\Gamma}
\end{tangle}
\ee
(The fourfold vertex in the right diagram represents the morphism $m^{(2)}=m\circ m\otimes\id$.)
Then $p$ is idempotent (up to a scalar) and its kernel is an ideal. Thus the image of $p$ is a
commutative Frobenius subalgebra of $\Gamma$. The latter is called the {\bf center} of $\Gamma$
since it is the ordinary center in the case $\2C=\Vect_k^{\mathrm{fin}}$. 

\item Application to TQFT: Every finite dimensional semisimple $k$-algebra $A$ gives rise to a TQFT in
$1+1$ dimensions via triangulation (Fukuma/Hosono/Kawai \cite{FHK}). By the classification of TQFTs
in $1+1$ dimensions \cite{dijk, abrams,kock}, this TQFT corresponds to a commutative Frobenius
algebra $B$ (in $\Vect_k^{\mathrm{fin}}$), with $A=V(S^1)$ and the product arising from the pants
cobordism. The latter is given by the vector space associated with the circle and the multiplication
is given by the pants cobordism. One finds $B=Z(A)$, and $B$ arises exactly as the image of $A$
under the above projection $p$. (This works since every semisimple algebra is a Frobenius algebra.)

\item If $\2C$ is braided, but not symmetric, we must choose between $c_{\Gamma,\Gamma}$ and
$c_{\Gamma,\Gamma}^{-1}$ in the definition (\ref{eq6}) of the idempotent $p$. This implies that a
non-commutative  Frobenius algebra will typically have two different centers, called the left and
right centers $\Gamma_l,\Gamma_r$. Remarkably, one then obtains an equivalence 
\[ E:\Gamma_l-\Mod_\2C^0\ \stackrel{\simeq}{\longrightarrow}\ \Gamma_r-\Mod_\2C^0 \]
of modular categories, cf.\ B\"ockenhauer, Evans, Kawahigashi \cite{BEK}, Ostrik \cite{ostrik} and
Fr\"ohlich, Fuchs, Runkel, Schweigert \cite{FS,FFRS}.
Conversely, if $\2C$ is modular,  every triple $(\Gamma_l,\Gamma_r,E)$ as above arises from a
  non-commutative algebra in $\2C$, \cite{KR}. (The latter is unique only up to Morita equivalence.)

\item This is relevant for the classification of CFTs in two dimensions: The latter are constructed
from a pair $(A_l,A_r)$ of chiral CFTs and some algebraic datum (`modular invariant') specifying how
the two chiral CFTs are glued together. In the left-right symmetric case, where the two chiral
theories coincide $A_l=A_r=A$, the above result indicates that Frobenius algebras in $\2C=\Rep\,A$
are the structure to use. This is substantiated by a construction, using TQFTs, of a `topological
from a modular category $\2C$ and a Frobenius algebra $\Gamma\in\2C$, cf.\  Fuchs, Runkel,
Schweigert, cf.\ \cite{FuRuS} and sequels. 

\item The Frobenius algebras in / module categories of $SU_q(2)-\Mod$ can be classified in terms of ADE
graphs. (Quantum MacKay correspondence.) Cf.\ B\"ockenhauer, Evans \cite{BE}, Kirillov Jr. and Ostrik
\cite{ko}, Etingof/Ostrik \cite{EO2}.

\item These results should be extended to other Lie groups. If $SU(2)$ already leads to the
ADE graphs (``ubiquitous'' according to \cite{hazew}), the other classical groups should give rise
to very interesting algebraic-combinatorial structures, cf.\ e.g.\ \cite{ocn5,ocn6}. 

\item More generally, when the two chiral theories $A_l, A_r$, and therefore the associated modular
  categories $\2C_l,\2C_r$ differ, it is better to work with triples $(\Gamma_l,\Gamma_r,E)$, where
  $\Gamma_{l/r}\in\2C_{l/r}$ are commutative algebras and
  $E:\Gamma_l-\Mod_{\2C_l}^0\rarr\Gamma_r-\Mod_{\2C_r}^0$ is a braided equivalence. (By the above,
  in the left-right symmetric case $\2C_l=\2C_r=\2C$, this is equivalent to the study of
  non-commutative Frobenius algebras $\Gamma\in\2C$.) Now one finds \cite{mue18} a bijection between
  such triples and {\it commutative} algebras $\Gamma\in\2C_l\boxtimes\widetilde{\2C_r}$ of the maximal
  dimension $d(\Gamma)=\sqrt{\dim\2C_l\cdot\dim\2C_r}$. (This is a categorical version of Rehren's
  approach \cite{khr5} to the classification of modular invariants. It is based on studying local
  extensions $\2A\supset A_l\boxtimes\widetilde{A_r}$, corresponding to commutative algebras 
$\Gamma\in\2C_l\boxtimes\widetilde{\2C_r}$.)

\item There also is a concept of a center of an algebra $A$ in a not-necessarily braided tensor
  category $\2C$, to wit the {\bf full center} defined in \cite{dav3} by a universal property. While
  the full center is a commutative algebra in the  braided center $Z_1(\2C)$ of $\2C$, as apposed to
  in $\2C$ like the above notions of  center, there are connections between these constructions.
\item We close this section giving three more reasons why modular categories are interesting:

1. They have many connections with {\bf number theory}: 
\begin{itemize}
\item Rehren \cite{khr}, Turaev \cite{turaev}: 
\[ \sum_i d_i^2=|\sum_i d_i^2\omega_i|^2.\]
In the pointed case (all simple objects have dimension one) this reduces to
$|\sum_i\omega_i|=\pm\sqrt{|I|}$. For suitable $\2C$, this reproduces Gauss' evaluation of Gauss
sums. (Gauss actually also determined the sign of his sums.)
\item The elements of $T$ matrix are roots of unity, and the elements of $S$ are cyclotomic integers
\cite{BG,et1}. 
\item For related integrality properties in Y=TQFSs, cf.\ Masbaum, Roberts, Wenzl \cite{masR,masW}
and Brugui\`eres \cite{brug4}).
\item The congruence subgroup property: Let $N=\mbox{ord}\,T (<\infty)$. Then
\[ \mathrm{ker}(\pi: SL(2,\7Z)\rarr GL(|I|,\7C))\ \supset \
    \Gamma(N)\equiv\mathrm{ker}(SL(2,\7Z)\rarr SL(2,\7Z/N\7Z)). \]
For the modular categories arising from rational CFTs, this had been known in many cases and
widely believed to be true in general. Considerable progress was made by Bantay \cite{bantay2},
whose arguments were made rigorous by Xu \cite{xu2} using algebraic quantum field theory. Bantay's
work inspired a proof \cite{SZ} by Sommerh\"auser and Zhu for modular Hopf algebras, using the higher
Frobenius-Schur indicators defined by Kashina and Sommerh\"auser \cite{KS}. Finally, Ng and
Schauenburg proved the congruence property for all modular categories along similar lines, cf.\
\cite{ng3}, beginning with a categorical version of the higher Frobenius-Schur indicators \cite{ng2}.
\end{itemize}

2. A modular category $\2M$ gives rise to a surgery TQFT in $2+1$ dimensions
(Reshetikhin, Turaev \cite{rt2,turaev}). In particular, this works for $\2M=Z_1(\2C)$ when $\2C$ is
spherical fusion categories $\2C$ with $\dim\2C\ne 0$. Since such a category $\2C$ also defines a 
TQFT via triangulation \cite{BW2,GK1}, it is natural to expect an isomorphism $RT_\2M=BWGK_\2C$ of
TQFTs. (When $\2C$ is itself modular, this is indeed true by
$Z_1(\2C)\simeq\2C\boxtimes\widetilde{\2C}$ and Turaev's work in \cite{turaev}.) Recently, a general
proof of this result was announced by Turaev and Virelizier, based on the work of Brugui\`eres and
Virelizier \cite{brug9,brug8}, partially joint with S.\ Lack. (Notice in any case that the surgery
construction provides more TQFTs than the triangulation approach, since not all modular categories
are centers.) 

3. We close with the hypothetical application of modular categories to topological quantum computing
\cite{zwang}. There are actually two different approaches to topological quantum computing: The one
initiated by M.\ Freedman, using TQFTs in $2+1$ dimension and the one due to A.\ Kitaev using $d=2$
quantum spin systems. However, in both proposals, the modular representation categories are
central. Cf.\ also Z.\ Wang, E.\ Rowell et al. \cite{HRW,rsw1}.
\end{itemize} 


\section{Some open problems}\label{sec-6}
\begin{enumerate}
\item Characterize the hypergroups arising from a fusion category. (Probably hopeless.) Or at least 
  those corresponding to (connected) compact groups. 
\item Find an algebraic structure whose representation categories give all semisimple pivotal
  categories, generalizing Ostrik's result \cite{ostrik}. Perhaps this will be something like 
the quantum groupoids defined by Lesieur and Enock \cite{LE}?  
\item Classify all prime modular categories. (The next challenge after the classification of finite
  simple groups...)
\item Give a direct construction of the fusion categories associated with the two Haagerup
subfactors \cite{haage,AH,asaeda2}.
\item Prove that every braided fusion category $\2C/\7C$ embeds fully into a modular
  category $\2M$ with $\dim\2M=\dim\2C\cdot\dim Z_2(\2C)$. (This is the optimum allowed by the
  double commutant theorem, cf.\ \cite{mue11}.)
\item Find the most general context in which an analytic (i.e.\ non-formal) version of the
Cartier/ Kassel/ Turaev \cite{cartier,kt2} formal deformation quantization of a symmetric tensor
category $\2S$ with infinitesimal braiding can be given. (I.e.\ give an abstract version of the 
Kazhdan/Lusztig construction of Drinfeld's category \cite{kl} that does not suppose $\2S=\Rep\,G$.) 
\item Generalize the proof of modularity of $Z_1(\2C)$ for semisimple fusion categories to
  not necessarily semisimple finite categories (in the sense of \cite{EO1}), using Lyubashenko's
  definition   \cite{ly2} of modularity.
\item Likewise for the triangulation TQFT \cite{TV,BW2,GK1}. Generalize the relation to surgery TQFT
  to the  non-semisimple case. (For the non-semisimple version of the RT-TQFT in \cite{KL}.)
\item Hard non-commutative analysis: Every countable $C^*$-tensor category with conjugates and
  $\End\,\11=\7C$ embeds fully into the $C^*$-tensor category of bimodules over $L(F_\infty)$ 
and, for any  infinite factor $M$, into $\End(L(F_\infty)\overline{\otimes}M)$.  Here $F_\infty$ is
the free group with countably many generators and $L(F_\infty)$ the type $II_1$ factor associated to
its left regular representation. (This would extend and conceptualize the results of
Popa/Shlyakhtenko \cite{popa} on the universality of the factor $L(F_\infty)$ in subfactor theory.)  
\item Give satisfactory categorical interpretations for various generalizations of quasi-triangular
  Hopf algebras, e.g.\ dynamical quantum groups \cite{et0} and Toledano-Laredo's quasi-Coxeter algebras
  \cite{Tol2}. Soibelman's `meromorphic tensor categories' and the `categories with cylinder
  braiding' of tom Dieck and   H\"aring-Oldenburg \cite{tD} might be relevant -- and in any case
  they deserve further study.
\end{enumerate} 

\vspace{1cm}

\noindent{\it Acknowledgement}: I thank B. Enriques and C. Kassel, the organizers of the 
Rencontre ``Groupes quantiques dynamiques et cat\'egories de fusion'' that took place at CIRM,
Marseille, from April 14-18, 2008, for the invitation to give the lectures that gave rise to these
notes. (No proceedings were published for this meeting.) I am also grateful to N. Andruskiewitsch,
F. Fantino, G. A. Garc\'ia and M. Mombelli for the invitation to the ``Colloquium on Hopf algebras,
quantum groups and tensor categories'', C\'ordoba, Argentina, August 31st to September 4th, 2009, as
well as for their willingness to publish these notes. 

\vspace{1cm}

\noindent{\it Disclaimer}: While the following bibliography is quite extensive, it should be clear
that it has no 
pretense whatsoever at completeness. Therefore the absence of this or that reference should not be
construed as a judgment of its relevance. The choice of references was guided by the principal thrust
of these lectures, namely linear categories. This means that the subjects of quantum groups and low
dimensional topology, but also general categorical algebra are touched upon only tangentially.

\end{document}